\magnification 1200
\input amstex
\documentstyle{amsppt}
\vsize 7.85in
\voffset -1cm
\hoffset 4.5truemm
\topmatter

\rightheadtext{Vector field construction of Segre sets}
\leftheadtext{Jo\"el Merker}
\title Vector field construction of Segre sets
\endtitle
\author Jo\"el Merker
\endauthor
\address Laboratoire d'Analyse, Topologie et Probabilit\'es,
Centre de Math\'ematiques et d'Informatique, UMR 6632, 39 rue Joliot
Curie, F-13453 Marseille Cedex 13, France. Fax: 00 33 (0)4 91 11 35
52\endaddress
\email merker\@cmi.univ-mrs.fr\endemail
\thanks
\endthanks

\subjclass
(revised 2000) 32H02, 32C16\endsubjclass

\loadeufm

\def\L{{\Bbb L}}
\def\SS{{\Bbb S}}
\def\N{{\Bbb N}}

\def\R{{\Bbb R}}
\def\C{{\Bbb C}}
\def\K{{\Bbb K}}
\def\L{{\Bbb L}}

\define \dl{[\![}
\define \dr{]\!]}
\def\v{\vert}
\def\n{\vert\vert}
\def\dim{\text{\rm dim}}
\def\codim{\text{\rm codim}}
\def\1{{\text{\bf 1}}}
\def\genrk{\text{\rm gen-rk}_\C}
\def\sumg{\sum_{\gamma\in \N^m}}

\abstract

\centerline{\bf Table of contents~:}

\smallskip

{\bf \S1.~Geometry of finite type and review of CR-extension theory \dotfill 2.}

{\bf \S2.~Identification of CR vector fields and of Segre varieties 
\dotfill 4.}

{\bf \S3.~Segre varieties and extrinsic complexification \dotfill 10.}

{\bf \S4.~Segre sets and iterated complexifications \dotfill 12.}

{\bf \S5.~Segre varieties and conjugate Segre varieties \dotfill 13.}

{\bf \S6.~Complexification of orbits of CR vector fields on $M$ 
\dotfill 16.}

{\bf \S7.~Complexified Segre $k$-chains and Segre sets \dotfill 19.}

{\bf \S8.~Minimality and holomorphic degeneracy in low codimension 
\dotfill 25.}

{\bf \S9.~Orbits of systems of holomorphic vector fields \dotfill 34.}

{\bf \S10.~Segre geometry of formal CR manifolds \dotfill 37.}

{\bf \S11.~Application of the formalism to the regularity of
CR mappings \dotfill 38.}
\endabstract

\endtopmatter

\document

The metric and differentiable structures called {\it SubRiemannian}
(SR) in the anglo-saxon literature and {\it Carnot-Carath\'eodory} in
France, which are associated to $\Cal C^\infty$-smooth real manifolds
equipped with systems of vector fields satisfying the so-called {\it
Chow's accessibility condition} appear naturally in various domains of
mathematics. The consideration of SR structures finds motivations and
ramifications notably in Control Theory (with concrete applications in
robotics), in the Analysis of Partial Differential Equations (mainly
in the study of hypoelliptic operators and of subelliptic estimates),
in Probability (study of diffusion on manifolds and sublaplacian), in
Hamiltonian Mechanics, in Contact Geometry, in Real Algebraic
Geometry, and finally in {\it Cauchy-Riemann} geometry. Before
becoming a geometric field in its own, the study of subRiemannian
structures has been mainly impulsed by the celebrated theorem of
H\"ormander about the $\Cal C^\infty$ hypoellipticity [H\"o] of the
sums of squares type operators $X=\sum_{j=1}^k L_j^2$, where the 
vector fields $L_j$ with $\Cal C^\infty$
coefficients over $\R^n$ 
satisfy Chow's condition.  Subsequent developements in the analysis of
PDE's have been conducted by M\'etivier [Met], by
Rothschild-Stein [RS] (who introduced the notions of dilatations, of
weight, and the structure of nilpotent Lie group on the tangent space
to a SR structure at a regular point, {\it see} also [Mit]), by
Fefferman-Phong, by Varopoulos, by Jerison-Sanchez-Sanchez-Calle, and
others. The reader may consult the extensive Bourbaki survey by
Kupka [K] for further information and more complete references. The
paradigmatic examples of SR structures are the Gru$\check{\text{\rm
s}}$in plane and the Heisenberg group, which is of contact type. The
interest of contact structures (a particular case of nonholonomy) is
motivated by their link with the symplectic structures introduced by
Lagrange, Souriau and Weinstein as a geometric frame for classical
mechanics, and later highly developed in the Russian school, notably
after Arnold, Gromov and Eliashberg. Quite recently, some fundamentals
of the geometric aspects of SR structures have been developed in a
long article by Gromov [Gro], in which it is addressed an impressive 
number of open questions, conjectures, {\it etc.} 
and among other topics, the
ball-box theorem, the Hausdorff dimension of SR structures, their
imbedding, the inequalities of Sobolev type, the disc theorem, {\it
etc.} The interested reader may consult the introductory article by
Bella\"{\i}che [Bell]. Such SR structures appear also quite naturally
in real analytic (or algebraic) Cauchy-Riemann geometry, a field where
the questions of hypoellipticity are motivated concretely, {\it see}
[Trv], [Trp]. Since 1996, some more algebraic aspects of 
such structures in the real analytic case and especially some new
invariants called {\it Segre sets} have been introduced by 
Baouendi-Ebenfelt-Rothschild [BER1,2], notably in order to establish some
CR regularity theorems about CR mappings between CR
manifolds. The present article is exclusively devoted to a geometric
and detailed exposition of these local CR algebraic invariants and
aims to tend towards a presentation of their properties which is as
elementary as possible. A more general approach of algebraic and
analytic aspects of real analytic complexifiable SR structures
underlies our considerations, but we shall concentrate the exposition 
on real analytic CR manifolds.

\head \S1. Geometry of finite type and review of CR-extension theory 
\endhead

\subhead 1.1.~Real analytic CR structures and CR orbits \endsubhead
A CR-generic real analytic 
($\Cal C^{\omega}$) submanifold $M$ of $\C^n$ carries
two fundamental geometric invariants\,: 
\roster
\item"{\bf 1.}" 
The {\it complex tangent bundle $T^cM=\text{\rm Re} (T^{1,0}M)=\text{\rm
Re}(T^{0,1}M)$}, and\,:
\item"{\bf 2.}" 
The so-called {\it family of Segre varieties}. 
\endroster
Our main goal in this article is to explain in expository style how
these two invariants can be identified, geometrically. But let us begin
first with a quick historical review of these two objects. Let $\K=\R$
or $\C$.  For the discussion which we shall endeavour here, we assume that
the reader is familiar with the following {\it Orbit Theorem}, see
Nagano [N], Sussmann [Sus], Baouendi-Ebenfelt-Rothschild [BER] or \S9
below where we reprove it. Let $\SS=\{X_\alpha\}$, $1\leq \alpha \leq
a$, $a\in \N_*$ be a system of vector fields {\it with $\K$-analytic
coefficients} which is defined over a small $\K$-analytic manifold
$M$. We denote by $\text{\rm Lie} (\SS)$ the Lie algebra generated by
$\SS$, {\it i.e.} consisting of all the multiple Lie brackets of any 
length of elements of $\SS$.  Let $p\in M$ be arbitrary.

\proclaim{Orbit Theorem} \! 
\text{\rm ([N], [Sus])}
There exists a unique germ $\Cal O_\SS(M,p)$
of manifold-piece passing through $p$, called the \text{\rm 
$\SS$-orbit of $p$ in $M$}, satisfying\text{\rm \,:}
\roster
\item"{\bf (1)}"
{\text{\rm Minimality property\,:} Every manifold-piece $\Lambda_p$
through $p$ with the property that each vector field $X_\alpha\in \SS$
is tangent to $\Lambda_p$ must \text{\rm contain} $\Cal O_\SS(M,p)$.
\item"{\bf (2)}"
\text{\rm Reachability\,:}
Each point $q\in \Cal O_\SS(M,p)$ is the endpoint of a piecewise
smooth integral curve of various elements of $\SS$ in a finite number.
\item"{\bf (3)}"
\text{\rm Dimensionality\,:} The dimension of this manifold-piece is
equal to the dimension of the Lie algebra $\text{\rm Lie}(\SS)$ at
$p$. Furthermore, the dimension of $\text{\rm Lie}(\SS)$
is \text{\rm constant} and equal to $\dim_\K \Cal O_\SS(M,p)$
at each point $q\in \Cal O_\SS(M,p)$.
\endroster
\endproclaim

\noindent
By definition,
the CR orbits $\Cal O_{CR}(M,p)$ of $M$ are the $\SS$-orbits with 
$\SS=:T^cM$.

\subhead 1.2.~Smooth CR-extension theory \endsubhead \!
In the last decade, it appeared that the CR-orbits of $M$ adequately
govern the analytic extension of CR functions, in all the regularity
categories $\Cal C^\omega$, $\Cal C^{\infty}$ and $\Cal C^{2,\alpha}$
({\it see} [Trp], [Tu1,2], [J], [M1], and [BER2]). In fact, some important
preliminary results in the $\Cal C^{\omega}$ case have first focused
attention on Lie brackets and finite type conditions, not on CR
orbits. But after the works of Treves [Trv], Tr\'epreau [Trp], Tumanov
[Tu1,2], J\"oricke [J] and the author [M1], it became clear that there
exists a {\it one-to-one correspondence} between the CR-orbits $\Cal
O_{CR}\subset M$ of $M$ and the {\it analytic wedges} attached to the
orbits $\Cal O_{CR}$, to which CR functions on $M$ admit a holomorphic
extension. These analytic wedges are simply complex manifolds
with edge, attached to the CR-orbit and they have dimension $\dim_{\C}
\, \Cal W^{an}=\dim_{\R} \, \Cal O_{CR}-
\dim_{CR} M$ ({\it see} [Tu2] and [J] for the precise definitions). In
case $\dim_\R \, \Cal O_{CR}=\dim_\R \, M$, this statemement is the
theorem on global minimality of [J] and [M1]. A large number of
families of analytic discs of Bishop's type, which are attached to $\Cal
O_{CR}$ and to some subsequent deformations of $\Cal O_{CR}$, fill
in $\Cal W^{an}$ by propagation and deformation [Tu2], [J], [M1]. In
summary, CR-orbits are the good objects for understanding CR-extension
theory.

Also, let us give here a brief review of some of these works on CR
extension. Holomorphic extendability to one side of CR functions
defined over a strongly pseudoconvex hypersurface in $\C^2$ was
discovered by Hans Lewy in 1956 and subsequently worked out by many
mathematicians\,: Hill-Taiani, Boggess-Polking, Bedford,
Fornaess-Rea, Boggess-Pitt, Baouendi-Chang-Treves,
Baouendi-Rothschild, Hanges-Treves and others 
(we refer the reader to the books
[Bo], [BER2] for bibliographical data). In the early
eighties, Treves pointed out the importance of CR orbits and of
Sussmann's construction [Sus]. In 1988, Tumanov proved his celebrated
extension theorem [Tu1]. In 1990, after the work of Hanges and Treves
[HaTr], Tr\'epreau identified the
fundamental correspondence between CR orbits and CR extension by
establishing propagation of wedge extendability using {\smc FBI}
(Fourier-Bros-Iagolnizer) transform. Then Tumanov in [Tu2] pushed forward 
the extension theory by establishing propagation of {\it CR extendability}
by means of deformations of analytic discs.

\subhead 1.3.~Real analytic CR manifolds \endsubhead \!
The second invariant, called {\it Segre varieties}, disappears in the $\Cal
C^{\infty}$ category. Segre varieties indeed arise from
complexification of the defining functions of $M$, which must therefore be
real analytic. Nevertheless, it is well known that the Segre varieties
are important tools in the reflection principle and in the study of CR
mappings between $\Cal C^{\omega}$ CR manifolds. They were introduced
in a short note by Segre in 1931. Motivated by Poincar\'e's
equivalence problem, Segre sought local differential invariants
attached to $M$. Much later in 1977, exhuming Segre varieties, Webster
proved a celebrated theorem\,: a local biholomorphic map between two
strongly pseudoconvex real algebraic hypersurfaces in $\C^n$, $n\geq
2$, must be algebraic [W]. Two other related grounding articles were
written independently by Pinchuk [P] and Lewy [L]. Since then, this
area became an intensive subject of research.  
The popularity of this field is
certainly due to its connections with Algebra, Analysis and Geometry
as well as on the fact that a great number of variations on the theme
leave open many attractive questions ({\it cf.} the 
celebrated conjectures of Diederich and Pinchuk). 
The reader may consult [DP] or
the book [BER2] for the classical presentation and also
[Suk] for a differential Lie geometric
viewpoint, as in Cartan's and Segre's original articles. 
As a matter of fact, Segre varieties have been
developed substantially since 1996, in the works of 
Baouendi-Ebenfelt-Rothschild, especially about CR manifolds
of high codimension.

\head \S2. Identification of CR vector fields and Segre varieties \endhead

\subhead 2.1.~Segre sets and minimality criterion \endsubhead \!
Indeed, in the study of the algebraic regularity mapping problem,
Segre varieties have been used by
Sharipov-Sukhov, who introduced a geometric condition called
``Segre-transversality'', and then by Baouendi-Ebenfelt-Rothschild, 
to whom is due the definition of ``{\it Segre sets}'', a quite
important novelty in the subject. Basically, the Segre sets of $M$ are
simple unions of Segre varieties as follows\,: $Q_p^1:=Q_p$ is the Segre
variety through $p\in M$, with $p\in Q_p^1$\,; $Q_p^2$ is the union of
Segre varieties $Q_q$ as $q$ ranges over $Q_p^1$\,; $Q_p^3$ is the
union of Segre varieties $Q_q$ as $q$ ranges over $Q_p^2$, {\it and so
on.} Clearly, by the (easy) property $r\in Q_q \Longleftrightarrow
q\in Q_r$, one has $Q_p^k\subset Q_p^{k+1}$. Such Segre sets $Q_p^k$
are also biholomorphically invariant, as the $Q_p$'s are. As a matter
of fact, we intend to discuss in this paper this construction of Segre
sets and the celebrated ``{\it minimality criterion}'' proved in
[BER1,2] and slightly modified after the contribution of Zaitsev [Z].
Thus, let us recall that $M$ is called of {\it finite type} in the
sense of Kohn and Bloom Graham at $p$ if the Lie algebra generated by
the complex tangent bundle $T^cM$ spans $T_pM$ at $p$ ({\it cf.}
[BG]). Equivalently, the CR orbit of $p$ contains a neighborhood of
$p$ in $M$, by Chow's accessibility theorem [Cho] ({\it cf.} 
{\bf (2)} of the Orbit Theorem). 
Using a now well established terminology, will say that
$M$ is {\it orbit-minimal} at $p$, or shortly ``minimal'' at $p$, ``in
the sense of Tumanov'', {\it cf.}~[Tu1]. As we have mentioned, the
article [BER1] have provided a preliminary insight about the link
between Segre varieties and finite-typeness of CR vector fields by
establishing the following

\proclaim{Minimality Criterion} 
$M$ is minimal at $p$ if and only if there is an integer $\nu_p \leq
d+1$ such that the $(2\nu_p)$-th Segre set $Q_p^{2\nu_p}$ contains a
neighborhood of $p$ in $\C^n$.
\endproclaim

\remark{Remark}
In this precise form, the statement is in fact due to Zaitsev, {\it
see} [Z] or [BER2]. In [BER1], the authors have established that
$Q_p^{\nu_p}$ contains an open set $\Cal V$ of $\C^n$ with
$\overline{\Cal V}\ni p$ (but $p\not\in \Cal V$).  The least such
integer $\nu_p$ is a biholomorphic invariant of $(M,p)$. Briefly, the
core of the proof in [BER1,2] is to set up normal forms for generic
$\Cal C^{\omega}$ manifolds as follows from the work of Bloom-Graham
[BG], to introduce the classical H\"ormander numbers of $M$ at $p$, to
discuss then the homogeneous algebraic Taylor approximation $M^{app}$,
to show that $Q(M^{app})_p^{2\nu_p}$ contains a neighborhood of $p$ in
$\C^n$ and to conclude by a perturbation argument that it is so for
$Q(M)_p^{2\nu_p}$. We claim that such a proof can be substantially
modified and that the whole theory of Segre sets can be geometrized.
\endremark

\subhead 2.2.~Flows of complexified vector fields and Segre chains 
\endsubhead \!
Indeed, our main re\-sult in this article lies in {\it two canonical
observations} which will fill up an appropriate understanding of some
canonical (and still missing) links between the complex tangent bundle
$T^cM$ {\it and} the Segre varieties. In truth, {\it these links become
clearly visible only after complexifying $M$}. To offer a complete
presentation of these ideas, we need at first to introduce a good deal
of notation ({\it see} also \S3, \S4, \S5, and \S6 below).  
Thus, let $\Cal M\subset\C^{2n}$ denote the extrinsic complexification
of $M\subset
\C^n$. This complexification $\Cal M$ lives in
$\C_t^n\times \C_\tau^n$, 
where $t=(t_1,\ldots,t_n)\in \C^n$, $\tau=(\tau_1,\ldots,\tau_n)\in
\C^n$ and is naturally equipped with two projections $\pi_t:
\Cal M \to \C^n_t$ and $\pi_{\tau} : \Cal M \to \C^n_{\tau}$, the term $t$
denoting some coordinates in $\C^n$ near $(M,p)$ and the term
$\tau:=(\bar{t})^c$ being the ``complexification coordinate'' of
$\bar{t}$ ({\it see}
\S3.5). Set $m:=\dim_{CR} M$, $d:=\codim_\R M$, $m+d=n$. Now, let
${\Cal L}:=({\Cal L}^1,\ldots,{\Cal L}^m)$ and $\underline{\Cal
L}:=(\underline{\Cal L}^1,\ldots,\underline{\Cal L}^m)$ be some two
complexified {\it commuting} bases of the two complexied bundles
$(T^{1,0} M)^c$ and $(T^{0,1} M)^c$ respectively, which have both
holomorphic coefficients and are both tangent to $\Cal M$. A complete
description in local coordinates of $\Cal L$ and of $\underline{\Cal
L}$ is provided in \S5.9 below. As each one of the 
two distributions they define is
Frobenius-integrable, we call $\Cal L$ and $\underline{\Cal L}$ ``{\it
$m$-vector fields}''. We notice that, on the contrary, the distribution
defined by the pair $\{\Cal L,\underline{\Cal L}\}$ is {\it not}
Frobenius-integrable, {\it except} only in the Levi-flat case\,! Now, let
$p\in M$ and set $p^c:=(p,\bar p)\in\Cal M$. Let ${\Cal
L}_w(p^c):={\Cal L}_{w^1}^1 (\ldots {\Cal L}_{w^m}^m(p^c))$ denote the
$m$-flow of ${\Cal L}$, $w\in\C^m$, and {\it idem} for
$\underline{\Cal L}_{\zeta}(p^c)$, $\zeta\in \C^m$. We follow a
notation like in [Sus], but in other texts, such $m$-flow would be
denoted for instance by $\exp (w \Cal L)(p^c)$. The alternated
concatenations of such flows will be called {\it Segre $k$-chains},
for instance, for $k=2j$, it is the holomorphic map
$(w_1,\ldots,w_{2j})\mapsto \underline{\Cal L}_{w_{2j}}( {\Cal
L}_{w_{2j-1}}( \ldots \underline{\Cal L}_{w_2}( {\Cal
L}_{w_1}(p^c))))\in {\Cal M}$, where $w_1\in \C^m, \ldots, w_{2j}\in
\C^m$, or by induction\,: $\Cal L_{w_1}(p^c)$, $\underline{\Cal
L}_{w_2}(\Cal L_{w_1}(p^c))$, $\Cal L_{w_3}(\underline{\Cal
L}_{w_2}(\Cal L_{w_1}(p^c)))$, {\it etc.}  We refer the reader to
\S6.1 below for further explanations. For short, we shall denote by
$w_{(k)}:=(w_1,\ldots,w_k)\in\C^{mk}$ such a $k$-uple and we shall
denote these above concatenated flow maps by the abbreviated
expression $w_{(k)}\mapsto\Gamma_k(w_{(k)})$.
Clearly, in place of the above maps, we could have chosen the
collection of maps $w_1\mapsto \underline{\Cal L}_{w_1}(p^c)$,
$w_{(2)}\mapsto \Cal L_{w_2}(\underline{\Cal L}_{w_1}(p^c))$,
$w_{(3)}\mapsto
\underline{\Cal L}_{w_3}(\Cal L_{w_2}(\underline{\Cal L}_{w_1}(p^c)))$,
{\it etc.} {\it Both the two choices are valuable}. In fact, we see no
reason why there there should
be a pre\-ferred choice here. We therefore introduce a corresponding
notation $\underline{\Gamma}_k(w_{(k)})$ for these second
$k$-concatenated flow maps. Then the two maps 
$\Gamma_k$ and $\underline{\Gamma}_k$
are different, except in the Levi-flat case. In fact, it will
appear that these two maps correspond to each other under the
well known special symmetry $\sigma(t,\tau):=(\bar t,\bar \tau)$ of
$\Cal M$ which is re-defined in \S3.6 below, as follows\,:
$\sigma(\Gamma_k(w_{(k)}))=\underline{\Gamma}_k (\overline{w_{(k)}})$.
Finally, let us define the sets $\Cal S_{\bar p}^k:=\{\Gamma_k(w_{(k)})\:
\n w_{(k)}\n \leq \delta \}$, which will be called {\it complexified
Segre $k$-chains} and $\underline{\Cal S}_p^k:=
\{\underline{\Gamma}_k( w_{(k)})\: \n w_{(k)}\n \leq \delta \}$,
called {\it complexified conjugate Segre $k$-chain}, where 
$\delta>0$. By slight abuse
of terminology in the paper, we shall call {\it Segre $k$-chains}
either the holomorphic maps $\Gamma_k$ or the sets $\Cal S_{\bar
p}^k$.  As we work locally, the value of $\delta>0$ appears to be
unimportant. Notice that in our terminology, the {\it Segre sets} of
[BER1,2] live in ambient space whereas our {\it Segre chains} live in
the complexification $\Cal M$. In fact, these Segre chains can be
identified with the ``orbit-chains'' in the spirit of Sussmann [Sus],
which appear in the step by step construction of the orbit of the pair
of integrable distributions $\{\Cal L,
\underline{\Cal L}\}$ on $\Cal M$, when one constructs their
orbits, as we shall do in \S6 and \S7. 
Accordingly, the complexification $\Cal M$ will be called {\it
orbit-minimal at $p^c$} if $\Cal S_{\bar p}^k$ contains a neighborhood
of $p^c$  in $\Cal M$ for some $k\in\N_*$. This full neighborhood 
covering property holds equivalently for
$\underline{\Cal S}_p^k=\sigma(\Cal S_{\bar p}^k)$. In general, we
denote by $\Cal O_{CR}(M,p)\subset M$ the CR-orbit of $p$ in $M$ and
by $\Cal O_{\Cal L, \underline{\Cal L}} (\Cal M,p^c)$ the $\{\Cal
L,\underline{\Cal L}\}$-orbit of $p^c$ in $\Cal M$.  Our first main
observation can be summarized as follows\,:

\proclaim{Theorem~2.3}
Let $\Cal M:=(M)^c$, where $M$ is $\Cal C^\omega$ as above. Then
\roster
\item"{\bf (a)}"
Segre sets arise from Segre $k$-chains:
$Q_p^{2j}=\pi_t(\underline{\Cal S}_p^{2j})$ and $Q_p^{2j+1} \!
=\pi_\tau(\Cal S_{\bar p}^{2j+1})$.
\item"{\bf (b)}"
The generic ranks $r_k=\underline{r}_k$ of the maps $\Gamma_k$,
$\underline{\Gamma}_k$ are biholomorphic invariants.
\item"{\bf (c)}"
The complexification of CR-orbits satisfies $[\Cal O_{CR}(M,p)]^c=
{\Cal O}_{\Cal L, \underline{\Cal L}}(\Cal M,p^c)$.
\endroster
\endproclaim

\remark{Remarks}
{\bf 1.}~Part {\bf (c)} appears already in [BER1,2], with a
technically very different proof.  We believe that our proof given in
\S6 below is elementary and natural in itself, because by definition
the two generators $L$ and $\bar L$ for $T^{1,0}M$ and $T^{0,1}M$
which produce the CR-orbits admit as complexifications $L^c:=\Cal L$
and $(\bar L)^c=\underline{\Cal L}$. Then $[\Cal O_{L,\bar
L}(M,p)]^c={\Cal O}_{\Cal L,\underline{\Cal L}}(\Cal M,p^c)$ should
hold almost tautologically, if our notations are appropriate.  We will
show that such notations are consistent with the proof that we will
conduct in \S6.1.

{\bf 2.}~Clearly, part {\bf (a)} can be taken as a (new) definition of
Segre sets. One can even forget Segre sets and consider only
Segre chains in $\Cal M$, after working only with 
complexified objects. This is in fact what the author does in
[M2,3], when considering various algebraic or formal CR-mapping problems.
We shall provide in \S11 a summary of what is essential in the study of
these CR-regularity problems.

{\bf 3.}~Clearly also, looking at {\bf (a)} we see that we could have
also defined what could be called ``{\it conjugate Segre sets
$\overline{Q_p^k}$}'' (not explicitely defined in [BER1,2] or [Z]) 
as follows\,: $\overline{Q_p^{2j}}:=\pi_\tau({\Cal S}_{\bar p}^{2j})$ and
$\overline{Q_p^{2j+1}}:=\pi_t(\underline{\Cal S}_p^{2j+1})$. Such a
terminology is explained by the fact that they appear to be simply the
conjugate sets $\{t\in\C^n\:\bar t \in Q_p^k\}=\overline{\{t\in Q_p^k\}}$.
This remark is not innocent.
\endremark

\subhead 2.4.~Necessity of conjugate Segre varieties \endsubhead \!
In fact, inspired by this observation of the existence of
``symmetric'' Segre sets obtained by complex conjugation, we will be
engaged to {\it define a $($new$)$ family of ``conjugate Segre
varieties $\overline{Q_p}$'' inside $\C^n$} (before complexification)
and to provide them with the status of interesting biholomorphically
invariant objects which are ``symmetric'' to the usual Segre
varieties. We know for a fact that, although nowhere published, this
tentative idea is very well known in the folklore\,: if
$M=\{t\in\C^n\:
\rho(t,\bar t)=0\}$, why not put $\rho(p,\bar t)=0$ instead of
$\rho(t,\bar p)=0$ to define Segre varieties\,? This would yield
$Q_p=\{t\in\C^n\: \rho(t,\bar p)=0\}$ and $\overline{Q_p}=\{\bar
t\in\C^n\: \rho(p,\bar t)=0\}$ (warning\,: by Lemma~3.4 {\bf (c)}
below, we have in fact, $\{t\in\C^n\: \rho(p,\bar t)=0\}=Q_p$\,--nothing 
new\,!). In fact, because all the properties of the
$\overline{Q_p}$ are obviously the same ({\it modulo} complex
conjugation) as those of the $Q_p$, it may perharps
seem to be superfluous to endeavour a task of
duplication. However, inspired by the vision that the pair of
$m$-vector fields $\{\Cal L,\underline{\Cal L}\}$ induce two different 
symmetric intrinsic complex analytic foliation of $\Cal M$ ({\it see}
Theorem~2.6 and \S5.13 below), the author feels that the correct point
of view has to go further. To support these
declarations, he claims two general ideas\,:

\smallskip

\roster
\item"{\bf I.}" \
{\it The good
geometric objects appear after complexification.}
\item"{\bf II.}" \
{\it It is important to consider simultaneously the pairs
$\{Q_p,\overline{Q_p}\}$.} 
\endroster

\demo{Argumentation of Claim I}
Claim I bears on an imprecise but well established general philosophy about
the study of real analytic objects that has been carried out
intensively since the grounding works of Webster.
\qed
\enddemo

\demo{Argumentation of Claim II}
It is fairly well known that in principle, there is no preferred
choice between holomorphic and antiholomorphic functions\,: both enjoy
the same properties and are ``the same'' type of objects.
Analogously, there is no preferred choice between the $Q_p$ and the
$\overline{Q_p}$, between $T^{1,0}M$ and $T^{0,1}M$, between
$L$-derivations and $\overline{L}$-derivations in all the CR-mapping
problems. But contrary to the holomorphic category, where a choice
between $t$ and $\bar t$ can be decided {\it once for all} at the
beginning, in the CR category, {\it no choice can be really made once
for all}.  In fact, the simultaneous consideration of both $t$ and
$\bar t$ is crucial and any of the two choices would be in a certain
sense erroneous by incompleteness.  
Indeed, since some real structure is added, the
members of the above pairs are not ``one and the same object''\,: the
$Q_p$ really differ from the $\overline{Q_p}$ and they {\it come
together in the Cauchy-Riemann geometry of $M$}. Any biholomorphic or
anti-biholomorphic change of coordinates leaves invariant {\it the two
families $Q_p$ and $\overline{Q_p}$}. Therefore, there would perharps
exist a sort of general ``{\it principle of symmetry}'' in analytic CR
geometry, represented by the bar symmetry $\overline{(\bullet)}$.
Correspondingly, we will also observe that {\it taking complex
conjugates has a real geometric signification}, since it permutes the
two foliations $\Cal F_{\Cal L}$ and $\Cal F_{\underline{\Cal L}}$ of
Theorem~2.6 below.
\qed
\enddemo

\remark{Remarks}
{\bf 1.}~We have wished to insist on this matter because these
``philosophical'' ideas about complex conjugation in analytic CR
geometry are not yet currently argued in the literature and because
conjugate Segre varieties appear nowhere with an explicit status.

{\bf 2.}~However, we mention that in a {\it non-CR context}, Moser and
Webster introduced a pair of involutions $\{\tau_1,\tau_2\}$ attached
to a $2$-surface with isolated complex tangency in $\C^2$ (defined in
terms of the two projections $\pi_t$ and $\pi_\tau$ of the
complexification $\Sigma^c$ of $\Sigma$), which are intertwined by the
complexification $\sigma(t,\tau):=(\bar\tau,\bar t)$ of the complex
conjugation operator. These involutions exchange the (generically of
cardinal two) fibers of $\pi_t\v_{\Sigma^c}$ and
$\pi_\tau\v_{\Sigma^c}$.  In the CR context, one remark in [MW] p.~262
reminds something corresponding to {\bf (k)} and {\bf (l)} of
Theorem~2.10 below, although it is not formulated there {\it with
pairs} as we shall do. The two involutions of Moser and Webster should
correspond to our pair of flows maps $w\mapsto \Cal L_w(p^c)$ and
$\zeta\mapsto \Cal L_\zeta(p^c)$.  On the other hand, these two
involutions are considered {\it as symmetric pairs} in [MW], as we aim
to do in this article. Although the (more elementary) CR case is not
studied by Moser and Webster, there is a strong structural analogy
between their constructions and ours below. In conclusion, we would
like to say that our geometric constructions are deeply related 
to the ones in [MW].

{\bf 3.}~Furthermore in his recent works [M3], the author has observed
that working out simultaneously {\it reflection identities} and (new)
{\it conjugate reflection identities} in the formal CR-mapping problem
enables to establish a well known conjecture, the convergence of
formal equivalences between two real analytic minimal and
holomorphically nondegenerate CR manifolds. This illustrates the
interest of studying real analytic CR objects {\it as conjugate
pairs}. The content of such pairs of reflection identities will be
exposed briefly in \S11 below.
\endremark

\smallskip

In summary, in real analytic CR geometry, the complex conjugation
operator has an important signification and it deserves to be 
thematized as an independent geometric object linking the geometric and
the algebraic objects {\it as conjugate pairs}.

\subhead 2.5.~Segre type \endsubhead \!
From Theorem~2.3, we will easily derive as a particular case\,:

\proclaim{Theorem~2.6}
The following properties are equivalent\,:
\roster
\item"{\bf (d)}" 
The CR-generic $\Cal C^\omega$ manifold $(M,p)$ is 
$\{L,\bar L\}$-orbit-minimal at $p$.
\item"{\bf (e)}"
Its complexification 
$(\Cal M,p^c)$ is $\{\Cal L,\underline{\Cal L}\}$-orbit-minimal
at $p^c=(p,\bar p)$.
\item"{\bf (f)}"
There exists an integer $\mu_p\leq d+2$ such that the 
$(2\mu_p-1)$-th Segre chain
$\Cal S_{\bar p}^{2\mu_p-1}$ $($or equivalently, its conjugate
$\underline{\Cal S}_p^{2\mu_p-1})$
contains a small open set 
$\Cal U\subset \Cal M$ with $\Cal U \ni p^c$.
\item"{\bf (g)}"
There exist $w_{(2\mu_p-1)}^*\in\C^{m(2\mu_p-1)}$
arbitrarily close to the origin such that
$\Gamma_{2\mu_p-1}\: (\C^{m(2\mu_p-1)},w_{(2\mu_p-1)}^*)\to (\Cal M,p^c)$
is a submersion $($same property for the
map $\underline{\Gamma}_{2\mu_p-1}\:
(\C^{m(2\mu_p-1)},\overline{w_{(2\mu_p-1)}^*})\to (\Cal M,p^c))$.
\item"{\bf (h)}"
With this integer $\mu_p$, the Segre sets $Q_p^{2\mu_p-2}$ and 
$\overline{Q_p^{2\mu_p-2}}$ both contain a neighborhood of $p$.
\endroster
Furthermore, we have $\nu_p=\mu_p-1$ $($where $\nu_p$ is as in the 
Minimality Criterion$)$.
\endproclaim

\remark{Remarks}
{\bf 1.}
As $\nu_p$, the least such integer $\mu_p$ is a biholomorphic
invariant of $\Cal M$, called in this article the {\it Segre type of
$\Cal M$ at $p$}.

{\bf 2.} 
Notice that we state that $Q_p^{2\nu_p}$ contains
a neighborhood of $p$, so we recover the same 
minimality criterion as in [Z], [BER2]. This integer
$2\nu_p$ cannot be improved ({\it e.g.} as an {\it odd integer} 
$2\nu_p-1$) as Example~7.21 will show.
\endremark

\smallskip

Correspondingly, we can reformulate a (new) minimality
criterion in the complexification $\Cal M$ by using Segre chains instead
of Segre sets as follows.

\proclaim{Minimality criterion~2.7}
$(M,p)$ is orbit-minimal if and only if there is an integer $\mu_p\leq
d+2$ such that the $(2\mu_p-1)$-th Segre chain $\Cal S_{\bar p}^{2\mu_p-1}$
contains a neighborhood of $p^c$ in $\Cal M$.
\endproclaim

\noindent
We shall give a complete proof of it in \S7. For the moment,
as a matter of fact, and quite naturally, we claim 
that this (new) minimality criterion
becomes (almost) immediate (hence better than the previous one),
because we can check the following elementary
arguments\,:
\roster
\item"{\bf 1\,:}" \
{\bf (d)} $\Longleftrightarrow$ $\C \otimes T_pM=\text{\rm Lie
algebra}_p(\{L,\bar L\})$. We assume that this equivalence is well
known, {\it i.e.} we assume that the reader is familiar with Chow's
and Nagano's classical theorems ({\it cf.} the Orbit Theorem above).
\item"{\bf 2\,:}" \
The complexification $[\text{\rm Lie algebra}_p(\{L,\bar L\})]^c=
\text{\rm Lie algebra}_p(\{\Cal L,\underline{\Cal L}\})$.
\item"{\bf 3\,:}" \
Hence {\bf (d)} $\Longleftrightarrow$ {\bf (e)}, because
$(T_pM)^c=T_{p^c}\Cal M$.
\item"{\bf 4\,:}" \
Finally, the equivalence {\bf (e)} $\Longleftrightarrow$ {\bf (f)} is
by definition, {\it q.e.d.}
\endroster

\remark{Remarks}
{\bf 1.}~The equivalence {\bf (f)} $\Longleftrightarrow$ {\bf (g)} follows
from Sussmann's construction of orbits. We will recover it in \S7-8
below, where we refine the constructions of [Sus].

{\bf 2.}~In fact, it is also easy to check {\bf (f)}
$\Longleftrightarrow$ {\bf (h)}, hence to recover the (first)
minimality criterion of [BER1,2], [Z] from the Minimality
Criterion~2.7. This equivalence can be quickly
established by observing that $\pi_t\: (\Cal M,p^c)\to (\C_t^n,p)$ and
$\pi_\tau\: (\Cal M,p^c)\to (\C_\tau^n,\bar p)$ are submersions with
fibers satisfying $\pi_t^{-1}(p)\cap \pi_\tau^{-1}(\bar p)=p^c$.
\endremark

\subhead 2.8.~The geometry of the complexification \endsubhead \!
Because the minimality criterion becomes canonical in the
complexification, it seems to be more natural to consider and to work out
Segre chains in the complexification $\Cal M$ instead of considering
Segre sets in $\C^n$. Recently, it also became clear that the geometry
of the pair $\{\Cal L,\underline{\Cal L}\}$ is well adapted to the CR
mapping problems, {\it see} \S11 below. A final argumentation to give
evidence to this conclusion will be provided by the next fundamental
observation.

\subhead 2.9.~Double foliation by complexified Segre varieties 
\endsubhead \!
In fact, the Segre varieties $Q_p$ and the conjugate Segre varieties
$\overline{Q_p}$ also admit canonical extrinsic complexifications,
which are $m$-dimen\-sional complex manifolds contained in $\Cal M$. 
By definition, these complexifications will be $\Cal S_{\tau_p}
:=\{(t,\tau_p)\: \rho(t,\tau_p)=0\}$, where $\tau_p\in\C^n$ is {\it
fixed}, and $\underline{\Cal S}_{t_p}:=\{(t_p,\tau)\:
\rho(t_p,\tau)=0\}$, where $t_p\in\C^n$ is {\it fixed} 
({\it see} \S5.12 for further precisions).  Now, we can state our
second main geometrical observation ({\bf (l)} below).  A very quick
proof of this ``observational'' theorem is provided in
\S5.6, \S5.9, \S5.12 below and can be read independently of the rest
of the article.

\proclaim{Theorem~2.10}
Let $\delta>0$ be small. Recall  that $\Cal L$ and $\underline{\Cal L}$
are Frobenius-integrable.
\roster
\item"{\bf (i)}"
$\Cal L$ and $\underline{\Cal L}$ induce naturally two local flow
foliations $\Cal F_{\Cal L}$ and ${\Cal F}_{\underline{\Cal L}}$ of
$\Cal M$.
\item"{\bf (j)}" 
$\sigma(\Cal F_{\Cal L})=\Cal F_{\underline{\Cal L}}$ and 
their two leaves through $p^c$
satisfy $\Cal F_{\Cal L}(p^c)\cap \Cal F_{\underline{\Cal
L}}(p^c)=p^c$.
\item"{\bf (k)}"
The fibers of the projections $\pi_t$ and $\pi_{\tau}$ also coincide
with the leaves of the flow foliations $\Cal F_{\Cal L}$ and $\Cal
F_{\underline{\Cal L}}$, respectively.
\item"{\bf (l)}"
The leaves of the foliation $\Cal F_{\Cal L}$ are the Segre varieties
$\Cal S_{\tau_p}$ and the leaves of the foliation $\Cal
F_{\underline{\Cal L}}$ are the conjugate Segre varieties
$\underline{\Cal S}_{t_p}$\text{\rm :}
$$
\Cal F_{\Cal L}=\bigcup_{\tau_p\in\C^n, \n \tau_p\n <\delta}
{\Cal S}_{\tau_p} \ \ \ \ \text{\rm and} \ \ \ \ \
\Cal F_{\underline{\Cal L}}=\bigcup_{t_p\in\C^n, \n t_p\n <\delta}
\underline{\Cal S}_{t_p}
\tag 2.11
$$
\endroster

\noindent
In other words, the leaves of these two flow foliations are the two families
of complexified $($conjugate$)$ Segre varieties. In symbolic 
representation, for these two foliations, we have the 
correspondence\text{\rm \,:}
$$
\text{\it CR-flow foliations of} \ \  \Cal M
\ \ \Longleftrightarrow \ \ \text{\it Foliations by complexified
Segre varieties}. 
\tag 2.12
$$
\endproclaim

\remark{Remarks}
{\bf 1.}
Observe that in the foliated unions of eqs.~\thetag{2.11}, the
dimensions correspond\,: the dimension of $\Cal M$ is equal to the
leaf dimension $m$ plus $n$, which is equal to the dimension of the
parameter space $\n \tau_p\n < \delta$ (or $\n t_p\n <\delta$).

{\bf 2.}  
It is well known however that classical Segre varieties {\it do not}
in general foliate a neighborhood of $p$ in $\C^n$. For instance,
think of the simplest quadric $z_2=\bar{z}_2+iz_1\bar{z}_1$ in
$\C^2$. However, {\it when we pass to the complexification, we blow-up \,
$\bigcup_{q\in \Cal V_{\C^n}(p)} Q_q$ and $\bigcup_{q\in \Cal
V_{\C^n}(p)}\overline{Q_q}$ in a double foliation of $\Cal M$ by complex
$m$-dimensio\-nal Segre surfaces.}
\endremark

\subhead 2.13.~Summary \endsubhead \!
Our second observation, valuable only in $\Cal C^\omega$ category, 
is clear\,:
$$ 
\text{\it Complexified \, Segre \, varieties} \
\Longleftrightarrow \
\text{\it Complexified \, CR \, vector \, fields \, foliation.}
\tag 2.14
$$

\remark{Remark} 
Again, we should say at the informal level that
the above equivalence makes very transparent and canonical
the Minimality Criterion of [BER1,2], [Z] expressed by means of Segre
sets. To the knowledge of the author, although very preliminary in the
subject, this interpretation of the two complexified Segre varieties
as $\underline{\text{\rm foliation}}$ leaves of the two flow
$\underline{\text{\rm foliations}}$ by complexified CR $m$-vector
fields ({\it i.e.} $m$-dimensional Frobenius-integrable distributions)
is also new.
\endremark

\subhead 2.15.~Propagation philosophy \endsubhead \!
We notice furthermore that in the $\Cal C^{\infty}$ case, while Segre
varieties disappear, it was known that the remaining $\Cal C^{\infty}$
sections of $T^cM$ propagate some structural analytic properties of CR
distributions\,: vanishing in a neighborhood of a point, being
extendable to a wedge, are propagating properties along CR orbits
[Trv], [Trp], [Tu2], [J], [Me1]. {\it It is therefore natural that in
the $\Cal C^{\omega}$ minimal case, Segre varieties become the support
of propagating properties for CR mappings}, like algebraicity, jet
solvability, analytic regularity, which are propagating properties
along Segre chains, {\it see} [BER1,2], [M2,3], [Z]. The chronology of
this discovery is amazingly inversed\,: the propagation phenomena were
first better understood in the $\Cal C^{\infty}$ category than in the
$\Cal C^{\omega}$ one, perharps because the ideas and tools came from
Analysis and more specifically from Partial Differential Equations. By
the way, the author wonders whether there exist smooth geometrical
objects which fill in the disappearing of Segre varieties in the
various CR-mapping problems between smooth CR manifolds, in analogy
with Bishop's discs in the smooth CR-extension theory.

\subhead 2.16.~Organisation of the remainder of the article \endsubhead \!
In \S3 and in \S4, we summarize the classical presentation of Segre
sets due to Baouendi-Ebenfelt-Rothschild.  The core of the
article begins in \S5.4 where the geometry in the complexification
$\Cal M$ starts. In \S5, we establish Theorem~2.10 by an inspection of
some equations of $\Cal M$ in coordinates. Then in \S6 and in \S7, we
study the (extrinsic or intrinsic) complexifications of CR orbits and
we establish Theorem~2.6. Numerous examples illustrating the
(rather dry) general theory are provided in \S8, taking inspiration
from works of Freeman, of Loboda and of Ebenfelt. Finally, we
construct a slight refinement of Sussmann's constructions in \S9 to
prove Nagano's theorem using flows of vector fields instead of Lie
algebras.

\definition{2.18.~Closing remark}
None of the result presented here is really new and most of the
theorems that we (re)prove are in fact originally due to
Baouendi-Ebenfelt-Rothschild. Only our geometrical viewpoint makes a
difference.
\enddefinition

\definition{2.19.~Acknowledgement}
The author is grateful to an anonymous referee who suggested to
rewrite the manuscript in an expository style and to render the
topic the more attractive and the less technical possible.
\enddefinition

\head \S3. Segre varieties and extrinsic complexification \endhead

\subhead 3.1.~Real analytic CR-generic manifolds \endsubhead \!
To begin with, we need a good deal of preliminary material.
It is well known that real analytic CR manifolds are generic in their
(complex analytic) semi-local intrinsic complexification. 
Accordingly, the local study of $\Cal C^\omega$ CR
manifolds reduces to the study of the CR-generic ones. Thus, here and
in the sequel, let $M$ be a piece through the origin of a $\Cal
C^{\omega}$ CR-{\it generic} manifold in $\C^n$ and set $m:=\text{\rm
dim}_{CR} M$, $d:=\text{\rm codim}_{\R} \ M$, with $m+d=n$. Then there
exists a system of $d$-vectorial defining functions for $M$:
$\rho=(\rho_1,\ldots,\rho_d)$, which are $\Cal C^\omega$ in a
neighborhood $U$ of $0$ in $\C^n$ such that $M\cap U$ coincides with
the zero-set $\{t\in U \: \rho(t,\bar{t})=0\}$, where $\rho(t,\bar{t})
\in \R^d$ satisfies $\rho(0)=0$ and $\partial\rho_1\wedge \cdots
\wedge \partial\rho_m\neq 0$ over $M\cap U$, by genericity of $M$. 
Of course, as a germ at $0$ of a real analytic subset, $(M,0)$ is
independent of the choice of such defining equations. We can assume
that the linear coordinates $t\in\C^n$ are chosen in order that
$T_0M=\C^m\times \R^d$. We can also expand
$\rho(t,\bar{t})=\sum_{\mu,\nu\in \N_*^n} \rho_{\mu,\nu} \, t^{\mu} \,
\bar{t}^{\nu}\in\C\{t,\bar t\}^d$, where
$\rho_{\mu,\nu}=\bar{\rho}_{\nu,\mu}\in \C^d$, $\forall \
\mu,\nu\in\N_*^n$ ({\it cf.} Lemma~3.4 below).  Naturally, we shall 
assume that this $d$-vectorial
power series converges {\it normally} in an open polydisc $U$ through
$0$, a polydisc which we can (and will) assume to be ``big'' after an
eventual dilatation of coordinates, say $U:=2^n\Delta^n=:\Delta_2^n$.

\subhead 3.2.~Segre varieties \endsubhead \!
Now, let us denote by $t_p$ the coordinates of a fixed point $p\in
\Delta^n$. Then the {\it Segre variety} $Q_p:=\{t\in \Delta_2^n\:
\rho(t,\bar t_p):=\sum_{\mu,\nu\in\N_*^n}
\rho_{\mu,\nu} \, t^\mu \, \bar t_p^\nu =0\}$, obtained by a {\it 
polarization} of $\rho$ and by a substitution of $\bar t$ by the
constant $\bar t_p$, is a {\it connected smooth complex
$m$-dimensional submanifold of $\Delta_2^n$}, for all $\bar
t_p\in\Delta^n$ (possibly after a supplementary dilatation of
coordinates).  It is well known that the definition of the Segre
variety $Q_p$ does not depend on the choice of local defining
functions of $M$, that Segre varieties are biholomorphically
invariant, {\it i.e.} $h(Q_p)=Q_{h(p)}'$ for every biholomorphism $h$
of $\Delta_2^n$ where $M':=h(M)$ and the $Q_{p'}'$ are associated with
$p'$, that $q\in Q_p$ iff $p\in Q_q$ and that $p\in Q_p$ iff $p\in M$.
Usually, these properties are established thanks to a complexification of
$\rho$ as follows.

\subhead 3.3.~Complexification of the defining equations \endsubhead \!
To begin with, let us define the complexification $\rho(t,\bar{t})^c:=
\rho(t,\tau):=\sum_{\mu,\nu\in \N^n} \rho_{\mu,\nu} \, t^{\mu} \, 
\tau^{\nu}$, after replacing $\bar{t}$ by an independent variable $\tau$ 
in the series defining $\rho$ and
$\bar{\rho}(t,\tau):=\sum_{\mu,\nu\in \N^n}
\bar{\rho}_{\mu,\nu} \, t^{\mu} \, \tau^{\nu}$, so
$\overline{\rho(t,\tau)}=\bar{\rho}(\bar{t}, \bar{\tau})$.  We can
think that $\tau$ is the ``complexified variable of $\bar t$'' and
write $\tau=(\bar{t})^c$ with superscript $(\bullet)^c$ for
``complexified''. There are also complexifications of $\Cal
C^{\omega}$ vector fields, of $\Cal C^{\omega}$ differential forms. We
write\,: $[\chi(t,\bar{t})]^c=\chi(t,\tau)$, if $\chi$ is $\Cal
C^{\omega}$ and $[\sum_{j=1}^{n} a_j(t,\bar{t}) \partial/\partial t_j+
\sum_{j=1}^n b_j(t,\bar{t}) \partial/\partial \bar{t}_j]^c:=
\sum_{j=1}^n a_j(t,\tau)
\partial/\partial t_j+ \sum_{j=1}^n b_j(t,\tau)\partial/\partial\tau_j$,
whence $(Lf)^c=L^cf^c$. Now, we can state the elementary reality
properties from which follow easily the basic properties of Segre
varieties quoted above.

\proclaim{Lemma~3.4}
As $\rho(t,\bar{t})=\sum_{\mu,\nu\in
\N^n} \rho_{\mu,\nu} \, t^{\mu} \, \bar{t}^{\nu} \in \R^d$ is 
$\Cal C^\omega$ and real-valued, one has\text{\rm \, :}
\roster
\item"{\bf (a)}" \ 
$\bar{\rho}_{\mu,\nu}=\rho_{\nu,\mu}$, $\forall \ \mu, \nu
\in\N_*^n$\text{\rm ;}
\item"{\bf (b)}" \ 
$\rho(t,\tau)\equiv \bar{\rho}(\tau,t)$\text{\rm ;}
\item"{\bf (c)}" \ 
$\rho(t,\tau)=0$ if and only if $\rho(\bar{\tau},\bar{t})=0$.
\endroster
\endproclaim

\demo{Proof}
By reality, $\rho(t,\bar t)\equiv \bar\rho(\bar t,t)$, whence 
{\bf (a)}. Clearly, {\bf (a)} $\Rightarrow$ {\bf (b)} 
$\Rightarrow$ {\bf (c)}.
\qed
\enddemo

\subhead 3.5.~The extrinsic complexification of $M$ \endsubhead \!
We now introduce the most important object, thanks to which the
geometrical aspects of the CR geometry of $M$ will be clearly more
appearant.  To the complexification $\rho(t,\tau)$ of
$\rho(t,\bar{t})$ is canonically (and classically) associated an {\it
extrinsic complexification} $M^c$ of $M$ defined by\,:
$$
(M)^c:=\Cal M:=\{(t,\tau)\in
\Delta^{2n} \:
\rho(t,\tau)=0\}\subset \C^n_t\times \C^n_{\tau}, 
\tag 3.6
$$ 
which is a {\it connected complex $2m+d$-dimensional submanifold of
$\Delta^{2n}$}. As we can embed $\C^n$ in $\C^n_t \times \C^n_{\tau}$
to be the totally real plane, often called the {\it antidiagonal},
$\underline{\Lambda}:=\{(t,\tau)\in
\C^{2n} \: \tau=\bar{t}\}$, {\it i.e.} to be the manifold
graph of the map $t\mapsto \bar{t}$, we can also embed $M$ in
$\underline{\Lambda}\subset \C^n_t\times \C^n_{\tau}$, as being
$M=\{(t,\bar t)\in\Delta^n\times \Delta^n\}$, whence $M$ appears to be
a {\it maximally real submanifold of} $\Cal M$ ({\it cf.} [BER1,2]).
Finally, if $p\in M$, {\it i.e.} $t_p\in M$, we denote by $p^c$ the
point $(p,\bar p)$, {\it i.e.}  with coordinates $(t_p,\bar{t}_p)\in
\Cal M$. Then $p^c=\pi_t^{-1}(\{p\})\cap
\underline{\Lambda}$, where $\pi_t : \C^n_t\times \C^n_{\tau}\to
\C^n_t$, is the projection $(t,\tau)\mapsto t$. We also put $\pi_{\tau}:
\C^n_t\times \C^n_{\tau} \to \C^n_{\tau}$, $(t,\tau)\mapsto \tau$, so that
$p^c=\pi_{\tau}^{-1}(\{\bar{p}\})\cap \underline{\Lambda}$.  We will
keep these notations throughout the article.

\subhead 3.7. Antiholomorphic involution \endsubhead \!
Now, as in [MW], [BER], let us define the antiholomorphic self-map
$\sigma(t,\tau):=(\bar \tau,\bar t)$ of $\C^{2n}$. Clearly, $\sigma$
is involutive, {\it i.e.} $\sigma^2=\text{\rm Id}$. Notice also that
$\sigma(t,\bar{t})=(t,\bar{t})$, which means that $\sigma$ fixes
$\underline{\Lambda}$ point by point, so in particular
$\sigma\v_M=I_M$. By Lemma~3.4 {\bf (c)}, we have\,:

\proclaim{Lemma~3.8}
The complex manifold $\Cal M$ is fixed by $\sigma$, i.e. $\sigma(\Cal
M)=\Cal M$.
\endproclaim

\noindent
In fact, $\sigma$ appears to be the fundamental symmetry of $\Cal M$
which corresponds by complexification to the complex conjugation
$t\mapsto \bar t$ in $\C^n$. Furthermore, it is easy to observe that this
invariance property of $\Cal M$ under $\sigma$ is a characterizing
property of the fact that $\Cal M=M^c$ is a complexification.

\subhead 3.9.~Characterization of complexifications of submanifolds
\endsubhead \!
More generally, we have the following lemma, which we will need in \S6.1.
For its proof, which we leave as an 
(easy) exercise, one has to use Lemma~3.4 {\bf (c)}.

\proclaim{Lemma~3.10}
There is a one-to-one correspondence between real analytic subsets
$\Sigma\subset M$ and complex analytic subvarieties $\Sigma_1$ of
$\Cal M$ satisfying $\sigma(\Sigma_1)=\Sigma_1$ given by
$\Sigma\mapsto \Sigma^c=:\Sigma_1$, with inverse $\Sigma_1 \mapsto
\pi_t(\Sigma_1 \cap \underline{\Lambda})=:\Sigma$.
Furthermore, $\Sigma$ is a smooth submanifold if and only if
$\Sigma_1$ is.
\endproclaim

\head \S4. Segre sets and iterated complexifications \endhead

After these preliminaries, we can present the definitions of higher
order Segre varieties.  Following [BER1,2], [Z1], let $Q_p^0:=\{p\}$,
where $p\in\Delta^n$, let $Q_p^1:=\bigcup_{q\in Q_p^0} Q_q\cap
\Delta^n$, $Q_p^2:=\bigcup_{q\in Q_p^1} Q_q\cap
\Delta^n$ and $Q_p^{k+1} :=\bigcup_{q\in Q_p^k} Q_q \cap
\Delta^n$. Then $Q_p^k$ is called the {\it $k$-th Segre set}.
We would like to mention that in these definitions, some accuracy
about sets upon which unions are taken is necessary, {\it cf.} [BER2],
but we shall not enter into the details before \S7 where we introduce
the sightly different notion of {\it Segre chains}, from which we
shall recover these Segre sets. An analytical definition using the
defining equations $\rho_j(t,\bar t)=0$, $1\leq j\leq d$, of $M$, is
presented as follows by these authors. Consider the {\it iterated
complexifications} $\Cal M_k$, $k\in\N$, of $M$ (introduced
by Zaitsev) to be the sets $\Cal
M_{2j}(p):=\{(p,\tau_1,t_1,\tau_2,t_2,\ldots,t_{j-1},\tau_j)\in
\Delta^{2jn}\:
\rho(p,\zeta_1)=\cdots=\rho(t_j,\tau_j)=0, \
\rho(t_1,\tau_2)=\cdots=\rho(t_{j-1},\tau_j)=0\}$
for $k=2j$ and a similar definition for $\Cal M_{2j+1}$, which are
both {\it complex manifolds}.  Using Lemma~3.4 {\it passim}, it is
easy to see that the Segre sets are projections of the $\Cal M_k$'s
over their last $\C^n$ factor. Up to now, the Segre sets are just sets
and they have few structure, but using the above projection, one sees
that the Segre sets are the images of certain holomorphic maps $v^k$
with source $\C^{mk}$. These maps can be described in coordinates as
follows. Assume that the $d$ scalar equations of $M$ are given by
$z_j=\bar Q_j(w,\zeta,\xi)$ in the notations of \S5.6 below, where
$\bar Q_j(w,\zeta,\xi):=\xi+i\bar\Theta(w,\zeta,\xi)$ and let
$Q(w,\zeta,\xi)$ denote the $\C^d$-valued power series whose
components are the $Q_j(w,\zeta,\xi)$. Notice that we do not assume
that the coordinates are regular, {\it i.e.} that $\bar
Q(w,0,\xi)\equiv 0$. Here, we follow the terminology
of Ebenfelt and use the words ``regular coordinates'', 
instead of ``normal coordinates''. Following the presentation of
[BER2,3], we can write for each integer $k\geq
1$ these holomorphic mappings $v^k\: (\C^{mk},0)\to (\C^n,0)$ as
follows. For $k=0,1$, we define $v^0:= (0,0)\in\C^m\times\C^d$ and
$v^1(w_1):=(w_1,\bar Q(w_1,0,0))$. Then 
$v^2(w_1,w_2):=(w_1,\bar Q(w_1,w_2,Q(w_2,0,0)))$
and $v^3(w_1,w_2,w_3):=(w_1,\bar Q(w_1,w_2,Q(w_2,w_3,
\bar Q(w_3,0,0))))$. Furthermore, for $k=2j, 2j+1$, where $j\geq 1$,
define\,:
$$
\left\{
\aligned
&
v^{2j}(w_1,\ldots,w_{2j}):=
\left(w_1,\bar Q(w_1,w_2,Q(w_2,w_3,\bar Q(w_3,w_4,\ldots \right.\\
& 
\ \ \ \ \ \ \ \ \ \ \ \ \ \ \ \  \left.
\ldots,Q(w_{2j-2},w_{2j-1},\bar Q(w_{2j-1},w_{2j},Q(w_{2j},0,0)))
\ldots )))\right),\\
&
v^{2j+1}(w_1,\ldots,w_{2j+1}):=
\left(w_1,\bar Q(w_1,w_2,Q(w_2,w_3,\bar Q(w_3,w_4,\ldots \right.\\
&
\ \ \ \ \ \ \ \ \ \ \ \ \ \ \ \ \ \ \ \ \left.
\ldots, \bar Q(w_{2j-1},w_{2j},Q(w_{2j},w_{2j+1},\bar Q(w_{2j+1},0,0)))
\ldots )))\right),\\
\endaligned\right.
\tag 4.1
$$
In this definition, the principle of iteration comes from the
set-theoretic definition of the $Q_k^p$'s. Whereas the geometric and
set-theoretic definition of the $Q_k^p$'s needs convergent defining
equations $\rho_j(t,\bar t)$, the maps $v^k$ in eqs.~\thetag{4.2} have
sense both in the category of convergent and of formal mappings.  In
this article, we shall introduce a new geometric point of view on
Segre varieties and as a consequence, we shall modify and geometrize
the definitions of new maps $\psi^{2j}:=
\pi_\tau\circ \Gamma_{2j}$ and $\psi^{2j+1}:=\pi_t\circ
\Gamma_{2j+1}$ ({\it see} eq.~\thetag{6.4} and
eq.~\thetag{7.9}), which are equal to the maps
$v^k$ up to a linear reparametrization with integer coefficients and
up to complex conjugation.

\head \S5. Segre varieties and conjugate Segre varieties. \endhead

\subhead 5.1.~Definitions and basic properties
\endsubhead \! 
Before entering into the geometry of the complexification, we can
define the fundamental pair of {\it Segre varieties} and {\it
conjugate Segre varieties}.  As explained in \S2.4, it is necessary~to
study {\it simultaneously} Segre varieties together with conjugate
Segre varieties.  Of course, there will be no real novelty with the
properties of the second family. In fact, the main novelty of the
point of view only lies in the simultaneous consideration of {\it
both}.  Thus, the Segre variety, usually denoted by $Q_p=\{t\in
\Delta^n\:
\rho(t,\bar t_p)=0\}$, will be denoted {\it in the remainder 
of this article} by 
$$
S_{\bar{t}_p}:=\{t\in \Delta^n \:
\rho(t,\bar{t}_p)=0\}.
\tag 5.2
$$ 
We stress the notation $S_{\bar{t}_p}$ and not
``$S_{t_p}$'', nor ``$S_p$'', with the {\it bar of complex conjugation
over $t_p$ in the index}, as in the expression $\rho(t,\bar t_p)$. We
could also have used the notation $S_{\bar p}$. In fact, this notation
will be consistent with the following (new) definition of {\it
conjugate Segre varieties}\,: 
$$
\overline{S}_{t_p}:= \{\bar{t}\in
\Delta^n \:
\rho(t_p,\bar{t})=0\}.
\tag 5.3
$$ 
Indeed, thanks to Lemma~3.4 {\bf (c)}, we
clearly have $\overline{S}_{t_p}=\overline{S_{\bar{t}_p}}=
\overline{S}_{{\bar{\bar t}}_p}=\overline{S}_{t_p}$, which shows the
consistency of our notation, if, as usual, the set $\overline{E}$
denotes $\{\bar{t}\in \Delta^n\: t\in E\}$ for an arbitrary set
$E\subset \Delta^n$. Also, we have $\bar{t}\in
\overline{S}_{t_p}$ if and only if $t\in S_{\bar{t}_p}$, {\it etc.}
It is easy to check that the Segre and conjugate Segre varieties are
complex analytic submanifolds closed and connected in $U=\Delta_2^n$,
if $t_p\in\Delta^n$, and that they are independent of the choice of a
$\Cal C^{\omega}$ defining function for $M$. The $S_{\bar t_p}$ and
the $\overline{S}_{t_p}$ are, moreover, biholomorphically
invariant. We have indeed $h(S_{\bar t_p})=S_{\bar h(\bar t_p)}$ and
$\bar h(\overline{S}_{t_p})=\overline{S}_{h(t_p)}'$, where $h$ is
as in \S3.2. Now,
as we suggested to duplicate Segre varieties and sets for
completeness, let us employ the notations $S_{\bar t_p}$ instead of
$Q_p$. To mimic and duplicate the set-theoretic definitions
given in \S4 above, we can define\,: $S_{\bar{t}_p}^0:=\{\bar{t}_p\}$, and
$S_{\bar{t}_p}^1:=S_{\bar{t}_p}=
\bigcup_{\bar{t}\in S_{\bar{t}_p}^0} S_{\bar{t}}$,
$S_{\bar{t}_p}^2=
\bigcup_{t\in S_{\bar{t}_p}^1}
\overline{S}_t$,
and then inductively, for $j\in \N_*$, $S_{\bar{t}_p}^{2j} =\bigcup_{t\in
S_{\bar{t}_p}^{2j-1}} \overline{S}_t$ and
$S_{\bar{t}_p}^{2j+1}=\bigcup_{\bar{t}\in S_{\bar{t}_p}^{2j}}
S_{\bar{t}}$.  On the other hand, we can define
$\overline{S}_{t_p}^0:=\{t_p\}$, and
$\overline{S}_{t_p}^1:=\overline{S}_t=\bigcup_{t\in
\overline{S}_{t_p}^0} \overline{S}_t$, 
$\overline{S}_{t_p}^2:=\bigcup_{\bar{t}\in \overline{S}_{t_p}^1}
S_{\bar{t}}$, and inductively, for $j\in \N_*$,
$\overline{S}_{t_p}^{2j}:=\bigcup_{\bar{t}\in
\overline{S}_{t_p}^{2j-1}} S_{\bar{t}}$, and
$\overline{S}_{t_p}^{2j+1} =
\bigcup_{t\in \overline{S}_{t_p}^{2j}}
\overline{S}_t$. As we shall introduce a 
more geometrical point of view in the complexification, 
these definitions are not very important for us. We just mention that
this definition of Segre sets $S_{\bar t_p}^k$'s clearly
coincides with the definition of the $Q_p^k$'s and that we have the 
following elementary properties\,:

\roster
\item"{\bf (1)}"
$\overline{S_{\bar{t}_p}^k}=\overline{S}_{t_p}^k$ and
$S_{\bar{t}_p}^k=\overline{\overline{S}_{t_p}^k}$, $k\in \N$.
\item"{\bf (2)}"
$\overline{h}(S_{\bar{t}_p}^{2j})=
S_{\overline{h}(\bar{t}_p)}^{'2j}$,
$h(S_{\bar{t}_p}^{2j+1})=S_{\overline{h}(\bar{t}_p)}^{'2j+1}$,
$h(\overline{S}_{t_p}^{2j})=\overline{S}_{h(t_p)}^{'2j}$,
$\overline{h}(\overline{S}_{t_p}^{2j+1})=
\overline{S}_{h(t_p)}^{'2j+1}$.
\endroster

\noindent
Here, $h\: M\to M'$ is a local biholomorphism.  In principle, the
above unions should be taken over $t\in S_p^k \cap (\delta\Delta)^n$,
where $\delta>0$ is sufficiently small and $k\leq 3n$, to be
equivalent with the definition that we will give in \S7.8 below after
having analyzed the geometry of the complexification.

\subhead 5.4.~Complexications of Segre varieties \endsubhead \!
Here begins the main topic of the article\,: the appearance of
{\it geometry in the extrinsic complexification}. 
We shall first establish our
``observational'' Theorem~2.10. So, in the sequel, we shall work {\it
exclusively} in the complexification $\Cal M\subset \C_t^n\times
\C_\tau^n$. As in \S2.9, let us define here two very important
objects\,:
$$
\left\{
\aligned
&
\Cal S_{\tau_p}:= (S_{\bar{t}_p})^c:=\{(t,\tau)\in
\Delta^n \times \Delta^n \:
\tau=\tau_p, \ \rho(t,\tau_p)=0\},\\
&
\underline{\Cal S}_{t_p}:=(
\overline{S}_{t_p})^c:=
\{(t,\tau) \in\Delta^n \times \Delta^n \:
t=t_p, \ \rho(t_p,\tau)=0\}.
\endaligned\right.
\tag 5.5
$$ 
Here, $\tau_p\in\Delta^n$ and $t_p\in\Delta^n$ are both fixed points.
These two $m$-dimensional complex submanifolds will be called {\it
complexified Segre varieties} and {\it conjugate complexified Segre
varieties}. Such a terminology is consistent with our previous
definitions in \S5.1, as shows an examination of
eqs.~\thetag{5.2}~\thetag{5.3}~\thetag{5.5}.  Both these two families
$\Cal S_{\tau_p}$ and $\underline{\Cal S}_{t_p}$ are in fact contained
in $\Cal M$. Using Lemma~3.4 {\bf (c)}, it is easy to check that
$\sigma(\Cal S_{\tau_p})=\underline{\Cal S}_{\bar{\tau}_p}$ and
$\sigma(\underline{\Cal S}_{t_p})=\Cal S_{\bar{t}_p}$, which shows a
fundamental symmetry property, simply coming from
the complexification of the
relation $\overline{S_{\bar t_p}}=\overline{S}_{t_p}$.  A particular
case in eqs.~\thetag{5.5} is when the fixed points $\tau_p$ (or $t_p$)
belongs to $M$.  Now, to check Theorem~2.10, we need the following
precisions.  We need explicit defining equations.

\subhead 5.6.~Coordinates\endsubhead \!
As usual, there exist holomorphic coordinates vanishing at $p\in M$,
$(w,z)$, where $w=(w_1,\ldots,w_m)$ and $z=(z_1,\ldots,z_d)$,
$z=x+iy$, such that $T_0^cM=\C^m_w\times \{0\}$, $T_0M=\C_w^m \times
\R_x^d$ and some $d$ {\it real} equations of $M$ are given 
by $y=h(w,\bar{w},x)$, where
$h=\sum_{\alpha,\beta\in\N^m,k\in\N^d} h_{k,\beta,\alpha} \, x^k
w^{\beta}\bar{w}^{\alpha}$, $h(0)=0$, $dh(0)=0$,
$h_{k,\beta,\alpha}\in\C^d$, $\overline{ h_{k,\beta,\alpha}}=
h_{k,\alpha,\beta}$, $\n h_{k,\beta,\alpha}\n\leq c \, 2^{-(\v\alpha\v
+\v\beta\v+\v k\v)}$, $c>0$ (after dilatation). If we replace
$y=(z-\bar{z})/2i$, $x=(z+\bar{z})/2$ and solve in $z$ or in $\bar{z}$
using the analytic implicit function theorem, we obtain as new
equivalent equations for $M$ and for its extrinsic complexification\,:
$$
\left\{
\aligned
&
M: \ \ \ 
z=\bar{z}+i\bar{\Theta}(w,\bar{w},\bar{z}) \ \ \
\ \ \text{\rm or} \ \ \ \ \ \bar{z}=z-i\Theta(\bar{w},w,z),\\
&
\Cal M: \ \ \
z=\xi+i\bar{\Theta}(w,\zeta,\xi) \ \ \ \ \ \text{\rm or} \ \ \ \ \
\xi=z-i\Theta(\zeta,w,z),
\endaligned\right.
\tag 5.7
$$
where $\Theta(\zeta,t)=\sum_{\beta\in\N^m,\mu\in\N^n}
\theta_{\beta,\mu} \, \zeta^\beta t^\mu$, with
$\n \theta_{\beta,\mu}\n\leq c \, 2^{-(\v \beta\v+\v \mu\v)}$, $c>0$.
It is this explicit representation which will help us to define the
invariant objects of Theorem~2.10.  Sometimes we shall use
eqs.~\thetag{5.7} and sometimes we shall set $\bar Q(w,\bar w,\bar
z):=\bar z+i\bar\Theta(w,\bar w,\bar z)$, whence the equations of
$\Cal M$ should be written in a shorter way by $z=\bar Q(w,\tau)$ or
equivalently by $\xi=Q(\zeta,t)$. Finally, we have the following
relation
$$
\Theta(\zeta,w,\xi+i\bar \Theta(w,\zeta,\xi))\equiv 
\bar \Theta (w,\zeta,\xi),
\tag 5.8
$$
which follows from the comparison of the two equivalent equations of
$\Cal M$.  It is easy to see that conversely, to a complex analytic
$d$-vectorial series $\Theta$ satisfying eq.~\thetag{5.8} is
associated a unique {\it real} analytic CR-generic manifold $M$, 
and indeed a series $h(w,\bar w,x)$, such that $M^c$ is given 
by eq.~\thetag{5.7} (this point is proved in the book [BER2], \S4.2, 
but we shall not need it). Very important geometric
objects are the CR vector fields.

\subhead 5.9.~CR (1,0) and anti-CR (0,1) vector fields \endsubhead \!
These {\it pairs} of $m$-vector fields and their complexifications
with holomorphic coefficients are given in vectorial notation by\,:
$$
\left\{
\aligned
&
M: \ \ \
L=\frac{\partial }{\partial
w}+i\bar{\Theta}_w(w,\bar w,\bar z)\frac{\partial }{\partial z} \ \
\ \ \ \hbox{and} \ \ \ \ \ \bar{L}=\frac{\partial }{\partial
\bar w}-i\Theta_{\bar{w}}(\bar w, w,z) \frac{\partial }{\partial
\bar z}\\
&
\Cal M: \ \ \
\Cal L =\frac{\partial }{\partial
w}+i\bar{\Theta}_w (w,\zeta,\xi)
\ \frac{\partial }{\partial z} \ \ \ \
\ \hbox{and} \ \ \ \ \ \underline{\Cal L}=\frac{\partial }{\partial
\zeta}-i\Theta_{\zeta}(\zeta,w,z) \frac{\partial }{\partial \xi}.
\endaligned\right.
\tag 5.10
$$  
In fact, these notations stand for $\Cal L=(\Cal L^1,\ldots,\Cal L^m)$
and $\underline{\Cal L}=(\underline{\Cal L}^1,\ldots,
\underline{\Cal L}^m)$. Clearly, the second row is the complexification 
of the first. In other words, $L^c= \Cal L$ and $\bar{L}^c=
\underline{\Cal L}$. Here, we drop
the indices, for the reason that the behaviour of these (obviously)
commuting $\Cal L^j$ and commuting (too) $\underline{\Cal L}^j$ is
formally analogous to the behaviour of a single vector field. In fact,
Frobenius integrable $m$-dimensional distributions behave like simple
vector fields. We shall call them {\it $m$-vector fields}. We denote
by $\Cal F_{\Cal L}$ and $\Cal F_{\underline{\Cal L}}$ the two flow
foliations that they induce. Of course, the distribution defined by 
the pair $\{\Cal L,\underline{\Cal L}\}$ is {\it not} in general
Frobenius integrable, unless $M$ is Levi-flat. This is why Lie
brackets and orbits come in consideration. At this subject,
we shall need the following Lemma~5.10 below in \S6.5. Granted
the equivalences {\bf (b)} $\Longleftrightarrow$ {\bf (c)}
and {\bf (d)} $\Longleftrightarrow$ {\bf (e)} which are explained in
[Sus] (using Chow's theorem and Nagano's theorem) and that
we shall admit, the following lemma is 
immediate\,:

\proclaim{Lemma~5.11}
Let $p\in M$. Then
the following properties are equivalent\text{\rm \, :}
\roster
\item"{\bf (a)}" \
$\text{\rm Lie}_p \ (T^{1,0}M, T^{0,1} M)=\C\otimes T_pM$.
\item"{\bf (b)}" \
$\text{\rm Lie}_p \ (T^cM)=T_pM$.
\item"{\bf (c)}" \
$(M,p)$ is $T^cM$-orbit-minimal.
\item"{\bf (d)}" \
$\text{\rm Lie}_{p^c} \ (\Cal L, \underline{\Cal L})=T_{p^c}\Cal M$.
\item"{\bf (e)}" \
$(\Cal M,p^c)$ is $\{\Cal L, \underline{\Cal L}\}$-orbit-minimal. \qed
\endroster
\endproclaim

\noindent
Indeed, the only remaining equivalence to be checked\,:
{\bf (b)} $\Longleftrightarrow$ {\bf (d)}, is trivial, because the
complexification of an $\{L,\bar L\}$-Lie bracket of arbitrary length 
is the corresponding $\{\Cal L,\underline{\Cal L}\}$-Lie bracket 
of the complexified vector fields.

\subhead 5.12.~Complexified Segre varieties \endsubhead \!
Now, we achieve the proof of Theorem~2.10.
At first, we observe that the two Segre varieties and
conjugate Segre varieties, which can be rewritten as
$$
S_{\bar{t}_p}: \ z=\bar{z_p}+i\bar{\Theta}(w,\bar{w}_p,\bar{z}_p) \ \
\ \ \ \hbox{and}
\ \ \ \
\overline{S}_{t_p}: \
\bar{z}=z_p-i\Theta(\bar w_p,w_p,z_p)
\tag 5.13
$$ 
admit two different complexifications in $\Cal M$, $\Cal S_{\tau_p}$ and
$\Cal S_{t_p}$, which now can be rewritten in terms of our
choice of coordinates as
$$
\left\{
\aligned
&
\Cal S_{\tau_p}=
\Cal S_{\zeta_p,\xi_p}: \ \ \ \zeta=\zeta_p, \ \xi
=\xi_p, \ z=\xi_p+i\bar{\Theta}(w,\zeta_p,\xi_p) \ \ \ \ \ \ \ \ \ \ 
\hbox{and} \\
&
\underline{\Cal S}_{t_p}=
\underline{\Cal S}_{w_p,z_p}: \ \ \ w=w_p, \ z=z_p, \
\xi=z_p-i\Theta(\zeta,w_p,z_p).
\endaligned\right.
\tag 5.14
$$
It is well known that the first vector fields $L$ in~\thetag{5.10} are
tangent to the (classical) Segre varieties $S_{\bar t_p}$, whence
$\Cal L$ also is tangent to the complexification $\Cal S_{\tau_p}$
({\it cf.}~[DW]). Indeed, one verifies immediately that $\Cal L \Cal
S_{\tau_p}\equiv 0$.  A similar relation holds between $\bar L$ and
$\overline{S}_{t_p}$ and between $\underline{\Cal L}$ and
$\underline{\Cal S}_{t_p}$, that is to say $\underline{\Cal L}
\underline{\Cal S}_{t_p}\equiv 0$. More precisely, thanks to this
complexification and for dimensional reasons, an inspection of
eqs.~\thetag{5.7}, \thetag{5.10},
\thetag{5.13} and \thetag{5.14} readily shows the following\,:

\smallskip

\roster
\item"{\bf 1.}"
The $\Cal S_{\tau_p}$ and the $\underline{\Cal S}_{t_p}$ form {\it
families of integral complex analytic manifolds for the 
$m$-vector fields $\Cal L$ and
$\underline{\Cal L}$} respectively.
\item"{\bf 2.}"
The $\Cal S_{\tau_p}$ are the leaves of the flow foliation $\Cal F_{\Cal
L}$ of $\Cal M$ by $\Cal L$ and the $\underline{\Cal S}_{t_p}$ are the
leaves of the flow foliation $\Cal F_{\underline{\Cal L}}$ of $\Cal M$
by $\underline{\Cal L}$.
\item"{\bf 3.}"
Finally, we clearly have $\Cal S_{\tau_p}=\Cal M
\cap \pi_\tau^{-1}(\tau_p)$ and $\underline{\Cal S}_{t_p}=\Cal M\cap 
\pi_t^{-1}(t_p)$. Equivalently, the complexified Segre varieties are
simply the intersection of $\Cal M$ with the coordinate planes. 
\endroster

\smallskip
\noindent
In conclusion, Theorem~2.10 is established. There remain only elementary 
details to be checked and for which enough precisions are already
given in the text. \qed

\remark{Remark} 
Local foliations by Segre varieties in ambient space have been
introduced by Sharipov and Sukhov for interesting application to the
algebraic regularity mapping problem as follows ({\it see} [SS],
[CMS,\S4]). If $H$ is a $d$-dimensional complex manifold through the
origin with $T_0H \oplus T_0^cM=T_0\C^n$, it is easy to see that
$\bigcup_{q\in H} Q_q$ makes a holomorphic foliation $\Cal F_H$ by the
$m$-dimensional complex leaves $Q_q$ of a neighborhood of $0$ in
$\C^n$, but this foliation depends on the choice of $H$. On the
other hand, the double foliation $\{\Cal F_{\Cal L},\Cal
F_{\underline{\Cal L}}\}$ is intrinsic and biholomorphically
invariant. Finally, we mention that in codimension $d\leq 2$, 
minimality of $M$ at $p$ is equivalent to Segre transversality of
$M$ at $p$, {\it see} [M2].
\endremark

\head \S6. Complexification of orbits of CR vector fields on $M$\endhead

\subhead 6.1.~Flows \endsubhead \!
In this paragraph, we now begin the (slightly longer) proof of Theorems~2.3
and~2.6.
Let $w=(w^1,\ldots,w^m)\in \C^m$ and $p^c\in\Cal M$. We will denote by
$\Cal L_w(p^c)$ the composition of flows 
$\Cal L_{w^m}^m\circ \cdots \circ \Cal
L_{w^1}^1(p^c)$ alluded to in \S2.2. In fact, for every permutation $\varpi:
[\![1,m]\!]\to [\![1,m]\!]$, we have $\Cal
L_{w_{\varpi(m)}}^{\varpi(m)}\circ \cdots
\circ\Cal L_{w_{\varpi(1)}}^{\varpi(1)}(p^c)= 
\Cal L_{w_m}^m\circ \cdots \circ 
\Cal L_{w_1}^1(p^c)$, since the $\Cal L^j$'s commute.
The vectorial $m$-flows of $\Cal L$ and of $\underline{\Cal L}$ on
$\Cal M$ are simply given in coordinates by the following relations\,:
$$
\left\{
\aligned
&
\Cal L_w(w_p, \, \bar z_p+i\bar{\Theta}(w_p, \bar w_p,\bar z_p), \,
\bar w_p, \, \bar z_p)=\\
&
\ \ \ \ \ \ \ \ \ \ \ =\left(w_p+w, \, 
\bar z_p+i\bar{\Theta}(w_p+w, \bar w_p, 
\bar z_p), \, \bar w_p, \, \bar z_p\right),\\
&
\underline{\Cal L}_{\bar w}(w_p, \, z_p, \, \bar w_p, \,
z_p-i\Theta(\bar w_p, w_p,z_p))=\\
&
\ \ \ \ \ \ \ \ \ \ \ =\left(w_p, \, z_p, \, \bar w+\bar w_p, \,
z_p-i\Theta(\bar w+\bar w_p, \, w_p, \, z_p)\right).
\endaligned\right.
\tag 6.2
$$
Applying $\sigma$ to the eqs.~\thetag{6.2}, we get $\sigma(\Cal
L_{w}(p^c))=\underline{\Cal L}_{\bar{w}}(p^c)$ and
$\sigma(\underline{\Cal L}_{\zeta}(p^c))=\Cal L_{\bar{\zeta}}(p^c)$.
We shall study closely the concatenated flow maps $\Gamma_k(w_{(k)})$
defined in \S2.2. To understand them formally, let us write down these
first three concatenated flows maps at $p^c=0$, after choosing the
(shorter) notation $\xi=Q(\zeta,t)$ instead of
$\xi=z-i\Theta(\zeta,t)$ for the equations of $\Cal M$.  We express
the right-hand sides in $(w,z,\zeta,\tau)$ coordinates. Using
~\thetag{6.2}, we can write these three
first concatenated flow\,:
$$
\left\{
\aligned
&
\Cal L_{w_1}(0)=\left(w_1, \, \bar Q(w_1,0,0), \, 0, \, 0\right).\\
&
\underline{\Cal L}_{w_2}(
\Cal L_{w_1}(0))=\left(w_1, \, \bar Q(w_1,0,0), \, 
w_2, \, Q(w_2,w_1,\bar Q(w_1,0,0))\right).\\
&
\Cal L_{w_3}(
\underline{\Cal L}_{w_2}(
\Cal L_{w_1}(0))))=\left(w_1+w_3, \,
\bar Q(w_1+w_3,w_2,Q(w_2,w_1,\bar Q(w_1,0,0))),\right.\\
&
\left. \ \ \ \ \ \ \ \ \ \ \ \ \ \ \ \ \ \ \ \ \ \ \ \ \ \ \ \ 
\ \ \ \ \ \ \ \ \ \ \ 
,w_2,\, Q(w_2,w_1,\bar Q(w_1,0,0) )\right).\\
\endaligned\right.
\tag 6.3
$$
It is easy to check that the following relations hold  between the 
$\Gamma_k$ and the $v^k$:
$$
\left\{
\aligned
&
v^1(w_1)=\pi_t(\Gamma_1(w_1)),\\
&
\bar v^2(w_2,w_1)=\pi_\tau(\Gamma_2(w_1,w_2)),\\
&
v^3(w_3+w_1, \, w_2,w_1)=\pi_t(\Gamma_3(w_1,w_2,w_3))\\
&
\bar v^4(w_4+w_2, \, w_3+w_1, \, w_2,w_1)=\pi_\tau(
\Gamma_4(w_1,w_2,w_3,w_4)),\\
&
v^5(w_5+w_3+w_1, \, w_4+w_2, \, w_3+w_1, \, w_2,w_1)=\\
&
\ \ \ \ \ \ \ \ \ \ \ \ \ \ \ \ \ \ \ \ \ \ \ \ \ \ \ \ \ \ 
\ \ \ \ \ \ \ \ \
=\pi_t(\Gamma_5
(w_1,w_2,w_3,w_4,w_5)),
\endaligned\right.
\tag 6.4
$$
{\it etc.}, which shows that the maps $v^k$ (formal definition) and
$\Gamma_k$ (set-theoretic definition) coincide, up to a
reparametrization with integer coefficients and up to complex
conjugation.  Let us now establish Theorem~2.3~{\bf (c)}.

\proclaim{Proposition~6.5}
There is a one-to-one correspondence
between the CR orbits and their
extrinsic complexifications\text{\rm \,:}
\roster
\item"{\bf (1)}" 
$\Cal O_{CR}(M,p)^c=\Cal O_{\Cal L, \underline{\Cal
L}}(\Cal M, p^c)$, and 
\item"{\bf (2)}" 
$\Cal O_{CR}(M,p)=\pi_t(\underline{\Lambda}\cap
\Cal O_{\Cal L, \underline{\Cal L}} (\Cal M, p^c))$.
\endroster
\endproclaim

\demo{Proof}
By the Orbit Theorem of \S1.1, $\Cal O_{CR}(M,p)$ is a $\Cal C^\omega$
closed submanifold of $M$ through $p$.  Thus, let $\Cal O$ be a small
open connected manifold-piece of $\Cal O_{CR}(M,p)$ through $p$, and
let $\Cal O^c$ be its extrinsic complexification. Because $L|_{\Cal O}$ and
$\bar{L}|_{\Cal O}$ are tangent to $\Cal O$, the {\it generic
uniqueness principle} ({\it via} $\Cal O\subset
\underline{\Lambda}$, where $\underline{\Lambda}$ is
maximally real) entails that $\Cal L|_{\Cal O^c}$ and
$\underline{\Cal L}|_{\Cal O^c}$ are tangent to $\Cal O^c$. Therefore
$\Cal O^c$ is an integral manifold for $\{\Cal L, \underline{\Cal
L}\}$ through $p^c$, whence $\Cal O^c \supset \Cal O_{\Cal L,
\underline{\Cal L}}(\Cal M, p^c)$, since a characterizing property of the
orbit $\Cal O_{\Cal L,
\underline{\Cal L}}(\Cal M, p^c)$ is to say that it is the {\it smallest
integral manifold-piece} for $\{\Cal L, \underline{\Cal
L}\}$ through $p^c$.

Conversely, Let $\Cal N$ be a manifold-piece of $\Cal O_{\Cal L,
\underline{\Cal L}}(\Cal M, p^c)$ through $p^c$. We have
just shown that $\Cal N\subset \Cal O^c$, hence to finish the proof,
we want to show that $\Cal N\supset \Cal O^c$.  We claim that we have
$\sigma(\Cal N)=\Cal N$ as germs at $p^c$. Indeed, By definition, the
orbit is the following set of endpoints of concatenations of flows of
$\Cal L$ and of flows of $\underline{\Cal L}$ (notice that because
$\Cal L_{w_2} \circ \Cal L_{w_1}= \Cal L_{w_1+w_2}$ and
$\underline{\Cal L}_{\zeta_2} \circ \underline{\Cal L}_{\zeta_1}=
\underline{\Cal L}_{\zeta_1+\zeta_2}$ but $\Cal L$ and
$\underline{\Cal L}$ do not commute, there can be only two different
kinds of concatenated flow maps\,; we do not use the abbreviated
notation $\Gamma_k$ here)\,:
$$
\left\{
\aligned
&
\Cal O_{\Cal L,
\underline{\Cal L}}(\Cal M, p^c)= \{ \Cal L_{w_k} \circ \cdots \circ \Cal L_{w_2}
\circ
\underline{\Cal L}_{\zeta_1} \circ \Cal L_{w_1} (p^c) \: \\
&
\ \ \ \ \ \ \ \ \ \ \ \ \ \ \ \ \ \ \ \ \ \ \ \ \ \ \ \ \ \ 
\: \left. w_1, \zeta_1, w_2,\ldots,w_k\in \C \ \text{\rm small}, 
k\in \Cal N_*\right\} \
\bigcup\\
&
\ \ \ \ \ \ \ \ \ \ \ \ \ \ \ \
\bigcup \ \{\underline{\Cal L}_{\zeta_k}
\circ \cdots \circ\underline{\Cal L}_{\zeta_2} \circ \Cal L_{w_1}
\circ \underline{\Cal L}_{\zeta_1} (p^c) \:\\
&
\ \ \ \ \ \ \ \ \ \ \ \ \ \ \ \ \ \ \ \ \ \ \ \ \ \ \ \ \ \ \
\: \zeta_1,w_1,\zeta_2,\ldots,\zeta_k\in \C \ \text{\rm small}, k\in \N_*\}:=
E\cup F.
\endaligned\right.
\tag 6.6
$$ 
Now an examination of eqs.~\thetag{6.2} shows that we have $\sigma(\Cal
L_w(q))=\underline{\Cal L}_{\bar{w}}(\sigma(q))$ and
$\sigma(\underline{\Cal L}_{\zeta}(q))=
\Cal L_{\bar{\zeta}}(\sigma(q))$, for each $q\in \Cal M$. Consequently\,:
$$
\left\{
\aligned
&
\sigma(\Cal L_{w_k}\circ \underline{\Cal L}_{\zeta_{k-1}} \circ
\cdots \circ \Cal L_{w_1}(p^c))= \underline{
\Cal L}_{\bar{w}_k}(\sigma( {\underline{\Cal L}_{\zeta_{k-1}}}\circ
\cdots \circ \Cal L_{w_1} (p^c)))=\\
&
=\underline{\Cal L}_{\bar{w}_k} \circ
\Cal L_{\bar{\zeta}_{k-1}}\circ \cdots \circ\sigma(\Cal L_{w_1}(p^c))=
\underline{\Cal L}_{\bar{w}_k} \circ 
\Cal L_{\bar{\zeta}_{k-1}}\circ \cdots \circ \underline{
\Cal L}_{\bar{w}_1}(p^c),
\endaligned\right.
\tag 6.7
$$ 
since $\sigma(p^c)=p^c$. This proves $F=\sigma(E)$, hence $\sigma(\Cal
O_{\Cal L, \underline{\Cal L}}(\Cal M, p^c))= \Cal O_{\Cal L, 
\underline{\Cal L}}(\Cal M, p^c)$, so we have $\sigma(\Cal N)=\Cal N$ 
as announced.

By the theorem of Nagano, $\Cal N$ is smooth at $p^c$ and satisfies
$\sigma(
\Cal N)=\Cal N$. By Lemma~3.10, there exists $N:= \pi_t(\Cal N \cap
\underline{\Lambda})$ a unique germ at $p$ of a $\Cal C^{\omega}$
submanifold $N\subset M$ such that $N^c=\Cal N$.

et us denote $\Cal N=\{\rho(t,\tau)=0, \chi(t,\tau)=0\}$, so that
$N=\{\rho(t,\bar{t})=0, \chi(t,\bar{t})=0\}$. Then $\Cal L\rho=0$,
$\underline{\Cal L}\rho=0$, $\Cal L \chi=0$, $\underline{\Cal L}
\chi =0$ on $\{\rho=\chi=0\}$, since $\Cal N$ is an $\{\Cal L,
\underline{\Cal L}\}$-integral manifold. Therefore, after restriction
to $\{\tau=\bar{t}\}=\underline{\Lambda}$, we have $L\rho=0$,
$\bar{L}\rho=0$, $L\chi=0$ and $\bar{L}\chi=0$ on
$\{\rho(t,\bar{t})=0,
\chi(t,\bar{t})=0\}=N$, so that $N$ is an $\{L,\bar{L}\}$-integral
manifold. Thus by the minimality property of CR-orbits, we have
$N\supset
\Cal O$ as germs at $p$. By complexifying, we get $\Cal N\supset 
\Cal O^c$, as desired.
\qed
\enddemo

\smallskip
\noindent
In conclusion, Theorem~2.3 is established. It remains only to check
{\bf (a)} and {\bf (b)}, for which enough precisions 
are already given in the text. \qed

\subhead 6.8.~Summary\endsubhead \!
Of course, it follows at once from Proposition~6.5 that $M$ is
$\{L,\bar{L}\}$-minimal at $p$ {\it if and only if} \, $\Cal M$ is
$\{\Cal L,\underline{\Cal L}\}$-minimal at $p^c$. This gives the
equivalence {\bf (d)} $\Longleftrightarrow$ {\bf (e)} of
Theorem~2.6. It remains to establish the remaining equivalences.
Motivated by the construction of Segre sets by 
Baouendi-Ebenfelt-Rothschild, we shall introduce some special 
notation for the sets
appearing in eq.~\thetag{6.6}, and call them {\it Segre
$k$-chains}. A preliminary presentation of them was given in
\S2.2. Crucially, in their definition, {\it the Segre $k$-th
chains will be endowed with holomorphic maps $\Gamma_k$ coming from
their geometric definition}.

\head \S7. Complexified Segre k-chains and Segre sets \endhead

\subhead 7.1.~Segre chains as $k$-th orbits chains of vector fields
\endsubhead \!
At first, we come back to the concatenated
flow maps in eqs.~\thetag{6.6} to analyze them.
Looking at eqs.~\thetag{5.14} and~\thetag{6.2},
we see that the complexified Segre varieties of a point $p^c\in\Cal M\cap
\underline{\Lambda}$, can be defined by $\Cal S_{\tau_p}:=
\{\Cal L_w(p^c)\in \Delta^{2n} \: w\in
\Delta^m \}$ and $\underline{\Cal S}_{t_p}
:= \{\underline{\Cal L}_{\zeta}(p^c)\in\Delta^{2n} \:
\zeta\in \Delta^m\}$. At order $k=2$, we can define\,:
$$
\left\{
\aligned
&
\Cal S_{\tau_p}^2=\{\underline{\Cal L}_{\zeta_1}(
\Cal L_{w_1} (p^c))\in \Delta^{2n} \:
w_1, \zeta_1\in
\delta \Delta^m\},\\
&
\underline{\Cal S}_{t_p}^2=\{\Cal L_{w_1}
(\underline{\Cal L}_{\zeta_1}(p^c))\in \Delta^{2n}
\: \zeta_1,w_1
\in \delta \Delta^m\}.
\endaligned\right.
\tag 7.2
$$ 
Here, $\delta>0$ will be chosen in a while.  More generally, let us
define the $\Cal S_{\tau_p}^k$ and $\underline{\Cal S}_{t_p}^k$. We
recall at first that, because only two
``starting actions'' $\Cal L_{w_1}(p^c)$ and $\underline{\Cal
L}_{\zeta_1}(p^c)$ can make a difference in a concatenation of flows
of $\Cal L$ and of $\underline{\Cal L}$, {\it there can exist only two
different families of Segre $k$-chains}. We recall also that we
have abbreviated in \S2.2 by $w_{(k)}\mapsto
\Gamma_k(w_{(k)})$ (resp. by $w_{(k)}\mapsto
\underline{\Gamma}_k(w_{(k)})$) the concatenated flow maps starting
by $\Cal L_{w_1}$ (resp. starting by $\underline{\Cal L}_{w_1}$).
Here, $\Gamma_k$ and $\underline{\Gamma}_k$ depend on $p^c$, but as we
shall essentially fix $p^c$ in a while, we shall not introduce a
specific notational index for $p^c$ in $\Gamma_k$ and
$\underline{\Gamma}_k$. Further, we have observed the relation
$\sigma(\Gamma_k(w_{(k)}))=\underline{\Gamma}_k
(\overline{w_{(k)}})$. Now, since we are interested in local
considerations, let us choose $\delta>0$ such that, for all $k\leq
3n$, all $w_{(k)}\in (\delta\Delta^m)^k$ and all $p\in
((\frac{1}{2}\Delta)^{n}\times (\frac{1}{2}\Delta)^{n})) \cap
\Cal M$, then $\Gamma_k(w_{(k)})
\in\Cal M$ and $\underline{\Gamma}_k(w_{(k)})\in \Cal M$ too. Thus, to be
explicit, the {\it Segre $k$-chains} of $p^c\in \Cal M$ are the
sets\,:
$$
\left\{
\aligned
&
\Cal S_{\tau_p}^{2j}:=\{\underline{\Cal L}_{\zeta_j}\circ
\Cal L_{w_j}\circ\cdots\circ\underline{\Cal L}_{\zeta_1}
\circ\Cal L_{w_1}(p^c)\: \, w_1,\zeta_1,\ldots,w_j,\zeta_j\in 
\delta\Delta^m\}\\
&
\Cal S_{\tau_p}^{2j+1}:=\{\Cal L_{w_{j+1}}\circ
\underline{\Cal L}_{\zeta_j}
\circ\cdots\circ\underline{\Cal L}_{\zeta_1}\circ 
\Cal L_{w_1}(p^c)\: \, w_1,\zeta_1,\ldots,w_{j+1}\in
\delta\Delta^m\}\\
&
\underline{\Cal S}_{t_p}^{2j}:=\{\Cal L_{w_j}\circ
\underline{\Cal L}_{\zeta_j}\circ\cdots\circ \Cal L_{w_1}\circ 
\underline{\Cal L}_{\zeta_1}(p^c) \: \,
\zeta_1,w_1,\ldots,\zeta_j,w_j\in \delta\Delta^m\} \\
&
\underline{\Cal S}_{t_p}^{2j+1}:=\{\underline{\Cal L}_{\zeta_{j+1}}\circ
\Cal L_{w_j}\circ\cdots\circ
\Cal L_{w_1}\circ \underline{\Cal L}_{\zeta_1}(p^c) \: \,
\zeta_1,w_1,\ldots,\zeta_{j+1}\in \delta\Delta^m\}.
\endaligned\right.
\tag 7.3
$$
for $k=2j$ or $k=2j+1$, where $j\in \N$ and $k\leq 3n$. Clearly, we
have $\Cal S_{\tau_p}^k\subset \Cal M$ and $\underline{\Cal
S}_{t_p}^k\subset
\Cal M$. As $\sigma(\Cal L_w(q))=\underline{\Cal L}_{\bar w}(\sigma(q))$,
we have $\sigma({\Cal S}_{\tau_p}^k)=\underline{\Cal S}_{\bar t_p}^k$.
We shall now endeavour a detailed study of these concatenated
flow maps which will enable us to offer a slight refinement of the 
statements in Theorem~2.6.

\subhead 7.4.~ Segre type and multitype\endsubhead \! 
To begin with, we shall denote in the sequel by $\genrk(\varphi)$ the
{\it generic} rank of a holomorphic map $\varphi\: X\to Y$ of
connected complex manifolds. Here of course,
$\genrk(\Gamma_1)=\genrk(\underline{\Gamma}_1)=m$ and
$\genrk(\Gamma_2)=\genrk(\underline{\Gamma}_2)=2m$, which is also
evident in eqs.~\thetag{6.3}. We set $e_1:=
\genrk(\Gamma_3)-2m$ and, by induction $e_{k+1}:=\genrk(\Gamma_{k+1})
-e_k-2m$, whence $\genrk(\Gamma_k)=m+m+e_1+\cdots+e_k$ if $k\geq 3$,
and similarly, we can define the sequence $\underline{e}_k$ for
$\underline{\Gamma}_k$. We notice at once that we have
$\underline{e}_k=e_k$, since
$\sigma(\Gamma_k(w_{(k)}))=\underline{\Gamma}_k(\overline{w_{(k)}})$.
We claim that $e_l=0$ for all $l\geq k+1$ if $e_{k+1}=0$ and $e_k\neq
0$.  Indeed, we first choose a point $w_{(k)}^*$ arbitrarily close to
the origin in $\C^{mk}$ such that $\Gamma_k$ has (necessarily locally
constant) rank equal to $2m+e_1+\cdots+e_k$ at $w_{(k)}^*$ and we set
$q:=\Gamma_k(w_{(k)}^*)\in \Cal M$. Then by the rank theorem, the
image $H$ of a neighborhood $\Cal W^*$ of $w_{(k)}^*$ is a complex
manifold.  We claim that $\Cal L$ and $\underline{\Cal L}$ are both
tangent to $H$. For instance, to fix ideas, we assume that $k$ is even
(the odd case will be similar). Thus we can write $\Gamma_k(w_{(k)})=
\underline{\Cal L}_{w_k} (\cdots \Cal L_{w_1}(0))$, {\it i.e.} the chain 
$\Gamma_k$ ends-up with a $\underline{\Cal L}$. This shows that $H$ is
fibered by the leaves of $\Cal F_{\underline{\Cal L}}$, so
$\underline{\Cal L}$ is already tangent to $H$ at every point.  On the
other hand, if $\Cal L$ were not tangent to $H$ at every point, 
the chain $\Gamma_{k+1}=\Cal L_{w_{k+1}} (\Gamma_k(w_{(k)}))$
would escape from $H$ and we would have $r_{k+1}>r_k$, a contradiction.
Finally, as $\Cal L$ {\it and} $\underline{\Cal L}$ are
both tangent to $H$, it follows that their local flow at $q$ is contained in 
$H$, whence the range of the subsequent $\Gamma_l$, $l\geq k+1$,
is contained in $H$. Because they are holomorphic, this shows that 
their generic rank does not go beyond $r_k$, {\it q.e.d.}

Consequently, there exists a well-defined integer $\kappa_p\geq
0$ with $\kappa_p\leq d$ such that $e_1>0,\cdots,e_{\kappa_p}>0$ and
$e_l=0$, for all $l\geq
\kappa_p+1$. We call the integer $\mu_p:=2+\kappa_p$ the {\it Segre type}
of $\Cal M$ at $p^c$ and we call the $\mu_p$-tuple
$(m,m,e_1,\ldots,e_{\kappa_p})$ the {\it Segre multitype} of $\Cal M$
at $p^c$. This {\it Segre multitype} simply recollects all
the jumps of generic rank of the 
$\Gamma_k$'s. It is clear that Segre type and multitype are biholomorphic
invariants, because the Segre foliations defined by $\Cal L$ and
$\underline{\Cal L}$ are so. To summarize, we have\,:
\roster
\item"{\bf 1)}"
$\text{\rm gen-rk}_{\C} (\Gamma_{k+2})=
2m+e_1+\cdots+e_k= \text{\rm gen-rk}_{\C} 
(\underline{\Gamma}_{k+2})$, \ $0\leq k \leq \kappa_p$,
\item"{\bf 2)}"
$\text{\rm gen-rk}_{\C} (\Gamma_{k+2})=2m+e_1+\cdots+e_{\{\kappa_p\}}=
\text{\rm gen-rk}_{\C} (\underline{\Gamma}_{k+2})$, \
$\kappa_p\leq k \leq 3n-2$.
\endroster
The main advantage of dealing with $\Cal M$, $\Cal L$,
$\underline{\Cal L}$, $\Gamma_k$ lies in the fact that all these
objects are coordinate-free.  Even the two projections $\pi_t$ and
$\pi_{\tau}$ could be defined abstractly, because their fibers are the
leaves of the Segre foliations. Correspondingly, the following
geometric statement, which will be re-interpreted
thanks to a more general construction in \S9 below, should be
understood in a coordinate-free style. This statement finishes
the proof of Theorem~2.6, except of part {\bf (h)}.

\proclaim{Theorem~7.5}
There exist some $w_{(\mu_p)}^*\in\C^{m\mu_p}$ arbitrarily close to
the origin of the form $w_{(\mu_p)}^*=(w_1^*,\ldots,w_{\mu_p-1}^*,0)$
and small neighborhoods $\Cal W^*$ of $w_{(\mu_p)}^*$ in
$(\delta\Delta^m)^{\mu_p}$ such that, if we denote
$\omega_{(\mu_p-1)}^*:=(-w_{\mu_p-1}^*,\ldots, -w_1^*)$, then we 
have\text{\rm :}
\roster
\item"{\bf 3)}"
The map $\Gamma_{\mu_p}$ is of rank $2m+e_1+\cdots+e_{\kappa_p}$ at
$w_{(\mu_p)}^*$.
\item"{\bf 4)}"
$\Gamma_{2\mu_p-1}(w_{(\mu_p)}^*,\omega_{(\mu_p-1)}^*)=p^c$.
\item"{\bf 5)}"
The restricted map $\Gamma_{2\mu_p-1}\: 
\Cal W^*\times \omega_{(\mu_p-1)}^* \to (\Cal M,p^c)$
is of constant rank equal to $2m+e_1+\cdots+e_{\kappa_p}$.
\item"{\bf 6)}"
The image
$\Gamma_{2\mu_p-1} (\Cal W^*\times
\omega_{(\mu_p-1)}^*)$ is a manifold-piece $\Cal O^c$ of the
complexified CR orbit of $M$ through $p^c$ 
$($or equivalently of $\Cal O_{\Cal L,
\underline{\Cal L}}(\Cal M, p^c))$.
\item"{\bf 7)}"
$2m+e_1+\cdots+e_{\kappa_p}=\text{\rm
dim}_{\C} \, \Cal O^c=\text{\rm dim}_{\C} \,
\Cal O_{\Cal L, \underline{\Cal L}}(\Cal M, p^c)=
\text{\rm dim}_{\R} \, \Cal O$.
\endroster
\endproclaim

\noindent
As this statement is quite technical, it is perharps worth to see what 
happens in the (considerably more simple) hypersurface case and to give
precise elementary examples before going into the details of the proof.

\proclaim{Corollary~7.6} 
In particular, let $M$ be a real $\Cal C^\omega$ hypersurface, 
\text{\rm i.e.} let $d=1$. Then
\roster
\item"{\bf 8)}" 
$M$ is minimal at $p$
\ $\Longleftrightarrow$ \ $\mu_p=3$ \ $\Longleftrightarrow$ \
$\kappa_p=1$.
\item"{\bf 9)}"
$M$ is nonminimal at $p$ \ $\Longleftrightarrow$ \  $\mu_p=2$ 
\ $\Longleftrightarrow$ \
$\kappa_p=0$.
\endroster  
\endproclaim

\demo{Proof} 
Firstly, we give the explanation in terms of the above definitions.
Indeed, the $2m$-dimensional complex submanifold $\{\underline{\Cal
L}_{\zeta_1} (\Cal L_{w_1} (p^c)) : \ w_1 , \zeta_1 \in
\delta\Delta^m\}$ is already of codimension only one in $\Cal M$, {\it
i.e.} of dimension $2m$. Then either $\text{\rm gen-rk}_{\C}
(\Gamma_3)= 2m+1$ or $\text{\rm gen-rk}_{\C} (\Gamma_3)=2m$, 
and the corollary follows. Secondly, we explain
the same statement in a pedestrian way, using coordinates and relate
it for completeness to this first proof of Corollary~7.6. We assume that
$\Cal M$ is given by eqs.~\thetag{5.7}, where $w\in \C^{n-1}$ and
$z\in \C$, and that the Segre variety
associated with the origin $S_{\bar 0}$ is straightened to be the
complex hyperplane $\{(w,0)\}$. Of course, this can be achieved
after an eventual biholomorphism. Since $M$ is
minimal if and only if $S_{\bar 0}$ is not 
contained in $M$, then we have following dichotomy
expressed simply in terms of the single ($d=1$) defining equation\,:
\roster
\item"{\bf 8')}" 
$M$ is minimal at $p$
\ $\Longleftrightarrow$ \ $\Theta(\zeta,w,0)\not\equiv 0$.
\item"{\bf 9')}"
$M$ is nonminimal at $p$ \ $\Longleftrightarrow$ \  
$\Theta(\zeta,w,0)\equiv 0$.
\endroster
Now, looking at eqs.~\thetag{6.3}, we can compute the flow maps 
explicitely\,:
$$
\left\{
\aligned
&
\Cal L_{w_1}(0)=(w_1, \, 0, \, 0, \, 0),\\
&
\underline{\Cal L}_{w_2}(
\Cal L_{w_1}(0))=(w_1, \, 0, \, 
w_2, \, -i\Theta(w_2,w_1,0)),\\
&
\Cal L_{w_3}(
\underline{\Cal L}_{w_2}(
\Cal L_{w_1}(0))))=(w_1+w_3, \, -i\Theta(w_2,w_1,0)+i\bar\Theta(
w_1+w_3,w_2,\\
&
\ \ \ \ \ \ \ \ \ \ \ \ \ \ \ \ \ \ \ \ \ \ \ \ \ \ \ \
-i\Theta(w_2,w_1,0)),
w_2,\, -i\Theta(w_2,w_1,0)).
\endaligned\right.
\tag 7.7
$$
Clearly, if $\Theta(\zeta,w,0)\equiv 0$, then $\Gamma_3$ is of generic
rank equal to $2n-2$, and so on for all the subsequent
$\Gamma_k$'s. We thus recover {\bf 9)} above. On the other hand, if
$\Theta(\zeta,w,0)
\not\equiv 0$, it is also clear with this explicit representation after
performing elementary verifications that $\Gamma_3$ is of 
(maximal possible) generic rank
equal to $2n-1$. We thus recover {\bf 8)} above, which completes a
second proof of Corollary~7.6.
\qed 
\enddemo

\example{Examples~7.8}
We now give a simple example in the hypersurface case which illustrates
statements {\bf 5)} and {\bf 6)} of Theorem~7.5 in a very concrete
way. We let $M$ be the hypersurface of $\C^2$ of equation 
$z=\bar z+iw^2\bar w^2$. We choose $p=0$ and here $2\mu_0-1=5$. 
We compute\,:
$$
\left\{
\aligned
&
\Gamma_1(w_1)=
(w_1, \, 0, \, 0, \, 0)\\
&
\Gamma_2(w_1,w_2)=
(w_1, \, 0, \, w_2, \, -i w_1^2 w_2^2)\\\
&
\Gamma_3(w_1,w_2,w_3)=(w_1+w_3, \,iw_2^2[w_3^2+2w_1w_3], \,
w_2, \, -i w_1^2 w_2^2)\\
&
\Gamma_4(w_1,w_2,w_3,w_4)=(w_1+w_3, \, 
iw_2^2[w_3^2+2w_1w_3], \, w_2+w_4, \, \\
& 
\ \ \ \ \ \ \ \ \ \ \ \ \ \ \ \ \ \ \ \ \ \ \ \ \ \ \ \ \ \
iw_2^2[w_3^2+2w_1w_3]-i
[(w_2+w_4)(w_1+w_3)]^2)\\
&
\Gamma_5(w_1,w_2,w_3,w_4,w_5)=(w_1+w_3+w_5, \,
iw_2^2[w_3^2+2w_1w_3]-\\
&
\ \ \ \ \ \ \ \ \ \ \ \ \ \ \ \
-i[(w_2+w_4)(w_1+w_3)]^2)+i[(w_1+w_3+w_5)(w_2+w_4)]^2, \,  \\
&
\ \ \ \ \ \ \ \ \ \ \ \ \ \ \ \
w_2+w_4, \, iw_2^2[w_3^2+2w_1w_3]-i
[(w_2+w_4)(w_1+w_3)]^2).
\endaligned\right.
\tag 7.9
$$
The maps $\Gamma_k$ have range in $\Cal M$, on which either the
coordinates $(w,z,\zeta)$ or $(w,\zeta,\xi)$ can be chosen.  We do the
first choice for $k$ even and the second choice for $k$ odd. Thus, we
view $\Gamma_5$ as a map $\C^5\to \C_{(w,\zeta,\xi)}^3$, {\it i.e.} we
forget the second $z$-coordinate in the above expression of $\Gamma_5$.
Now, computing the $3\times 5$ Jacobian matrix of at the point
$(w_{(3)}^*,\omega_{(2)}^*)$ as in Theorem~7.5 which is necessarily of
the form $(w_1^*,w_2^*,0,-w_2^*,-w_1^*)$, and for which we clearly
have $\Gamma_5(w_1^*,w_2^*,0,-w_2^*,-w_1^*)=0$, we see that the
determinant of the first $3\times 3$ submatrix is equal to
$2iw_1^*(w_2^*)^2$. Thus, it is nonzero for an arbitrary choice of
$w_1^*\neq 0$ and $w_2^*\neq 0$. By the way, the question arises
whether the integer $(2\mu_p-1)$ in Theorem~7.5 is optimal. Incidentally, 
this example shows that it is optimal. Indeed, if we ask whether there
exists $w_{(4)}^*=(w_1^*,w_2^*,w_3^*,w_4^*)$ such that
$\Gamma_4(w_{(4)}^*)=0$ and the rank at $w_{(4)}^*$ of the
differential of $\Gamma_4$ equals $3$ (the dimension of $\Cal M$), then
looking at eqs~\thetag{7.9}, we get first $w_1^*+w_3^*=0$,
$(w_2^*)^2w_3^*[w_3^*+2w_1^*]=0$ and $w_2^*+w_4^*=0$, thus $w_{(4)}^*$
is necessarily of the form $(0,w_2^*,0,-w_2^*)$ or
$(w_1^*,0,-w_1^*,0)$. Viewing now $\Gamma_4$ as a map $\C^4\to
\C_{(w,z,\zeta)}^3$, and computing its $3\times 4$ Jacobian matrix at
such points, one sees that it is of rank $2$, which proves
the claim.

Let us now give an example in codimension $d=2$. To illustrate
part {\bf 3)} of Theorem~7.5, we consider the generic manifold
$M\subset \C^3$ given by $z_1=\bar z_1+iw\bar w$ and
$z_2=\bar z_2+iw\bar w(w+\bar w)$. The reader can check that 
the map $(w_1,w_2,w_3,w_4)\mapsto 
(w_1+w_3, \, iw_2w_3, \, i(w_1+w_3)(w_1+w_2+w_3)-iw_1w_2(w_1+w_2), \,
w_2+w_4)$ is indeed of generic rank equal to $4$. He can also compute
$\Gamma_7$ to test the other statements.
\endexample

\demo{Proof of Theorem~7.5}
Now, we can proceed to the general arguments, which just
use elementary properties of flows of vector fields, as in 
[Sus]. According to {\bf 2)},
$\Gamma_{\mu_p}$ is of generic rank
$2m+e_1+\cdots+e_{\kappa_p}$. Consequently, for every point
$w_{(\mu_p)}^*\in(\delta\Delta^m)^{\mu_p}$ outside of some proper
complex subvariety, then the map $\Gamma_{\mu_p}$ is of rank
$2m+e_1+\cdots+e_{\kappa_p}$ {\it at} \, $w_{(\mu_p)}^*$. In fact, we
claim that we can even choose such a $w_{(\mu_p)}^*$ of the form
$(w_1^*,\ldots,w_{\mu_p-1}^*,0)$, {\it i.e.}  with
$w_{\mu_p}^*=0$. Indeed, as $\Gamma_{(k)}(w_{(k)}) =[\Cal L
\,\text{\rm or} \, \underline{\Cal
L}]_{w_k}(\Gamma_{(k-1)}(w_{(k-1)}))$, the following easy consequence
of the fact that the flow maps are local biholomorphisms 
yields the existence of such a
$w_{(\mu_p)}^*$ of this form. Let us state this property independently.

\proclaim{Lemma~7.10}
Let $w\in\C^m$, $w'\in\C^{m'}$, let $\Gamma(w')\in\C^\nu$ be
holomorphic in $w'$ with $\Gamma(0)=0$ and let $\varphi(w,w'):= {\Cal
L}_w(\Gamma(w'))$ or $\varphi(w,w'):=
\underline{\Cal L}_w(\Gamma(w'))$. Then $\varphi$
attains its maximal rank at some points of the form $(0,{w'}^*)$. \qed
\endproclaim

\noindent
Now, we fix such a $w_{(\mu_p)}^*$ of the form
$(w_1^*,\ldots,w_{\mu_p-1}^*,0)$, which satisfies {\bf 3)}
and we check that it satisfies the other claims.
Let $\omega_{(\mu_p-1)}^*:=(-w_{\mu_p-1}^*,\ldots,-w_1^*)$.
First, {\bf 4)}
is easy\,: suppose for instance $\mu_p$ is even, then we have
$\Gamma_{2\mu_p-1}(w_{(\mu_p)}^*,\omega_{(\mu_p-1)}^*)=
\Cal L_{-w_1^*}\circ\cdots\circ \Cal L_{-w_{\mu_p-1}^*}\circ
\underline{\Cal L}_0\circ \Cal L_{w_{\mu_p-1}^*}\circ \cdots \circ
\Cal L_{w_1^*}(p^c)=p^c$, because $\underline{\Cal L}_0(q)= q$, 
$\Cal L_{-w}\circ \Cal L_w\equiv
\hbox{Id}$ and $\underline{\Cal L}_{-\zeta}\circ
\underline{\Cal L}_{\zeta}=\hbox{Id}$. The odd case is similar. 
Now, we proceed to {\bf 5)}. The restricted map
$w_{(\mu_p)}\mapsto\Cal L_{-w_1^*}\circ\cdots\circ \Cal
L_{-w_{\mu_p-1}^*}\circ
\Gamma_{\mu_p}(w_{(\mu_p)})$ (again written in case $\mu_p$ is even), 
is clearly of rank $2m+e_1+\cdots+e_{\kappa_p}$ at the point 
$w_{(\mu_p)}^*$,
because the maps $q\mapsto \Cal L_{-w_1^*}\circ\cdots\circ \Cal
L_{-w_{\mu_p-1}^*}(q)$ is a local biholomorphism, by definition of
flows. Notice that $\Gamma_{2\mu_p-1}$ is then of {\it constant rank}
equal to $2m+e_1+\cdots+e_{\kappa_p}$ in a neighborhood of
$(w_{(\mu_p)}^*,\omega_{(\mu_p-1)}^*)$ in $\Cal W^*\times
\omega_{(\mu_p-1)}^*$, since, by {\bf 2)},
$2m+e_1+\cdots+e_{\kappa_p}$ is already the maximum value of all the
generic ranks of the $\Gamma_k$'s. This proves {\bf 5)}.  By
definition, the orbit $\Cal O_{\Cal L,\underline{\Cal L}}(\Cal M,p^c)$
is the union of the ranges of the maps $\Gamma_k$ and of the
$\underline{\Gamma}_k$'s. It is easy to check that this double union
coincides in fact with the union of only the $\Gamma_k$'s (or of only the
$\underline{\Gamma}_k$'s), simply because, setting $w_1=0$, we have
$\Gamma_k(0,w_2,\ldots,w_k)\equiv
\underline{\Gamma}_{k-1}(w_2,\ldots,w_k)$.
Thanks to the constant rank property {\bf 5)}, we already know that
this orbit contains the $(2m+e_1+\cdots+e_{\kappa_p})$-dimensional
manifold-piece passing through $p^c$\,: 
$\Cal N:=\Gamma_{2\mu_p-1} (\Cal
W^*\times\omega_{(\mu_p-1)}^*)$. Because by {\bf 2)} the next generic
ranks for $k\geq 2\mu_p-1$ do not increase and because of the
principle of analytic continuation, we then deduce that all the ranges
of the subsequent $\Gamma_k$'s are contained in this manifold
piece $\Cal N$ and it follows that
$\Cal L$ and $\underline{\Cal L}$ are
tangent to this manifold-piece. 
In conclusion, we get {\bf 6)} and {\bf 7)}, which completes
the proof of Theorem~7.5.
\qed
\enddemo

\remark{Remarks}
{\bf 1.}~Of course, we get the same statements in Theorem~7.5 as well
with $\underline{\Gamma}_{2\mu_p-1}$ instead of $\Gamma_{2\mu_p-1}$.

{\bf 2.}~The above proof is a special application of Sussmann's
construction ([Sus] studies mainly the 
$\Cal C^\infty$ category) to the particular case where the distributions are
holomorphic, which simplifies greatly the reasonings. A generalization
of it is provided in
\S9 below.
\endremark

\smallskip

We have argued of the simplicity of the hypersurface case.
Another particular case is $m=\text{\rm dim}_{CR} M = 1$. 
We clearly have the following.

\proclaim{Corollary~7.11}
If $\dim_{CR}M=1$, then $e_1=\ldots= e_{\kappa_p}=1$. Thus $\mu_p=d+2$
and $\kappa_p=d$ if and only if $M$ is minimal at $p$.
\endproclaim

\remark{Remark}
The question arises whether one can produce examples of $M$'s where
the Segre multitype assumes arbitrary possible values
$(m,m,e_1,\ldots,e_{\kappa_p})$ with $1\leq e_j\leq m$, provided
$e_1+\cdots+e_{\kappa_p}=d$. We cannot answer this question.
\endremark

\subhead 7.12.~Segre sets in ambient space \endsubhead \!
It is interesting to relate our construction with the construction of
Segre sets given in [BER1,2] which was
presented informally in \S4. To this aim, it is convenient to define
Segre sets in ambient space as certain projections of Segre chains.  
The germs $S_{\bar{t}_p}^{2j+1}:=\pi_t(\Cal S_{\bar{t}_p}^{2j+1})\subset
\Delta^n$, $\overline{S}_{t_p}^{2j+1}:=\pi_{\tau}(\underline{\Cal
S}_{t_p}^{2j+1})\subset
\Delta^n$, $S_{\bar{t}_p}^{2j}:=\pi_{\tau}(\Cal S_{\bar{t}_p}^{2j})\subset
\Delta^n$ and $\overline{S}_{t_p}^{2j}:=
\pi_t(\underline{\Cal S}_{t_p}^{2j})\subset \Delta^n$ will be called the {\it
Segre $k$-sets and conjugate Segre $k$-sets $(k=2j$ or $k=2j+1)$.}
Notice that because of eqs.~\thetag{6.2} (see also
eqs.~\thetag{6.3},~\thetag{7.7}), the action of the flow of $\Cal L$
leaves unchanged the $(\zeta,\xi)$-coordinates, and vice versa, the
action of the flow $\underline{\Cal L}$ leaves unchanged the
$(w,z)$-coordinates. This is why in the definition of Segre $k$-sets,
we alternately project in the $\C_t^n$-space and in the $\C_\tau^n$
space.  Using eqs.~\thetag{6.2-3-4}, the reader can easily observe that
this definition coincides with the definition of [BER2] 
in the form given in \S4, {\it modulo} a correct
correspondence of the ranges of the parameters $w_{(k)}$. Further, with
this new definition, we recover the definition of {\it conjugate
pairs} of Segre sets given in \S5.2.  Now, let
us define the maps $\psi^1(w_{1}):=\pi_t(
\Gamma_1 (w_{1}))$, $\psi^2(w_1,w_2):=\pi_{\tau}(\Gamma_2 (w_1,w_2))$ 
and more generally\,:
$$
\psi^{2j}(w_{(2j)}):=
\pi_{\tau}(\Gamma_{2j}(w_{(2j)}))
 \ \ \ \text{\rm and} \ \ \psi^{2j+1}(w_{(2j+1)}):=\pi_t(
\Gamma_{2j+1}(w_{(2j+1)})).
\tag 7.13
$$ 
Similarly also, we can define the maps $\underline{\psi}^k$ by
$\underline{\psi}^{2j}(w_{(2j)}):=\pi_t(\underline{\Gamma}_{2j}(w_{(2j)}))$
and $\underline{\psi}^{2j+1}(w_{(2j+1)}):=
\pi_{\tau}(\underline{\Gamma}_{2j+1}(w_{(2j+1)}))$.
Again, by an inspection of eqs.~\thetag{6.2-3-4}, the reader can
observe that up to a reparametrization and up to complex conjugation,
these maps are ``the same'' as the maps $v^k$, $\bar v^k$ of
\S4. We need the following.

\proclaim{Lemma~7.14}
For $0\leq k \leq \kappa_p$, we have\text{\rm \,:}
$$
m+ \text{\rm gen-rk}_{\C} (\psi^{k+1})= \text{\rm
gen-rk}_{\C} (\Gamma_{k+2})=2m+e_1+\cdots+e_k,
\tag 7.15
$$ 
and $\text{\rm gen-rk}_{\C} (\psi^{k+1})=
m+e_1+\cdots+e_{\kappa_p}$ for $\kappa_p \leq k \leq 3n-2$.
\endproclaim
 
\remark{Remark}
Of course, this same statement also holds with $\underline\Gamma_{k+2}$ and
$\underline{\psi}^{k+1}$ instead.
\endremark

\demo{Proof}
For $k=0$, we have $\psi^1(w_1)=(w_1, i \bar \Theta(w_1,0,0))$,
whence $\text{\rm gen-rk}_{\C} (\psi^1)=m$
obviously. Recall that, by eq.~\thetag{6.3}, we have
$$
\Gamma_2(w_1,w_2)=
\left(w_1, \, i\bar\Theta(w_1,0,0), \, 
w_2, \, i\bar\Theta(w_1,0,0) -i\Theta(w_2,w_1,\bar\Theta(w_1,0,0))\right),
\tag 7.16
$$ 
so $m+\text{\rm gen-rk}_{\C} (\psi^1)= \text{\rm gen-rk}_{\C}
(\Gamma_2)=2m$. More generally, for $k=2j$, we have:
$$
\left\{
\aligned
&
\Cal L_{w_{2j+1}}(\Gamma_{2j}(w_{(2j)}))=
\Cal L_{w_{2j+1}}(
w(w_{(2j)}), z(w_{(2j)}),\zeta(w_{(2j)}),\xi(w_{(2j)}))=\\
&
\ \ \ \ \ \ \ \ \ \ \ \ \ \ \ \ \ \ \ \ \ \ \ \ 
=\left(w_{2j+1}+w(w_{(2j)}), \,
\xi(w_{(2j)})+i\bar{\Theta}(w_{2j+1}+\right.\\
&
\ \ \ \ \ \ \ \ \ \ \ \ \ \ \ \ \ \ \ \ \ \ \ \ \
\left. \ \ +w(w_{(2j)}), \zeta(w_{(2j)}),
\xi(w_{(2j)})), \, \zeta(w_{(2j)}), \, \xi(w_{(2j)})\right).
\endaligned\right.
\tag 7.17
$$
As in Example~7.8, we choose the coordinates $(w,\zeta,\xi)$ on $\Cal
M$, whence we consider the map $\Gamma_{2j+1}$ 
in eq.~\thetag{7.17} to have range in
$\C_{(w,\zeta,\xi)}^{2m+d}$. It is then the map
$(w_{(2j)},w_{2j+1})\mapsto (w_{2j+1}+w(w_{(2j)}),
\zeta(w_{(2j)}), \, \xi(w_{(2j)})$. Then it follows
immediately that
$$
\text{\rm gen-rk}_{\C} (\Gamma_{2j+1})=
m+\text{\rm gen-rk}_{\C} [w_{(2j)}\mapsto (\zeta(w_{(2j)}),
\xi(w_{(2j)}))]= m+\text{\rm gen-rk}_{\C} \psi^{2j}. 
\tag 7.18
$$
This completes the proof of Lemma~7.14.
\qed
\enddemo

We now define the {\it Segre type of $M$ at $p\in M$} (not to be
confused with $\mu_p$) to be the smallest integer $\nu_p$ satisfying
$\text{\rm gen-rk}_{\C} (\psi^{\nu_p})=\text{\rm gen-rk}_{\C}
(\psi^{\nu_p+1})$.  By eq.~\thetag{7.15}, we readily observe that in
fact, we have $\nu_p=\kappa_p+1$ and $\nu_p=\mu_p-1$. The Segre type
{\it of $M$} can be related to its CR orbits as follows in the next
subparagraph. These last results will close up our presentation of the
general theory of Segre chains.

\subhead 7.19.~Intrinsic complexification of CR-orbits \endsubhead \!
By the {\it intrinsic complexification $N^{i_c}$ of a real CR manifold
$N$}, we understand the smallest complex analytic manifold containing
$N$ in $\C^n$, which satisfies $\dim_\C N^{i_c}=\dim_{CR} N+\codim_{CR}
N$.  Let $\Cal O_p$ denote a manifold-piece of $\Cal O_{CR}(M,p)$
through $p$ and let $\Cal O_p^{i_c}$ be its intrinsic
complexification.  By construction, the ranges of the $\psi^{2j}$'s
are contained in $\C_\tau^n$, but we will prefer to work in $\C_t^n$
(although it is equivalent in principle to work in $\C_\tau^n$), hence
we shall consider the $\underline{\psi}^{2j}$'s. We can now re-prove
that $\text{\rm gen-rk}_{\C} (\underline{\psi}^{\nu_p})=\text{\rm
dim}_{\C} \Cal O_p^{i_c}$ (a theorem due to
Baouendi-Ebenfelt-Rothschild, {\it see} [BER1,2]) and that the range
of $\underline{\psi}^{2\nu_p}$ contains a manifold-piece of $\Cal
O^{i_c}$ through $p$, as announced in \S2.5.

\proclaim{Theorem~7.20}
There exist some points 
$\underline{v}_{(2\nu_p)}^*\in\C^{m2\nu_p}$ arbitrarily close to
the origin and small neighborhoods $\Cal V^*$ of 
$\underline{v}_{(2\nu_p)}^*$ in
$(\delta\Delta^m)^{2\nu_p}$ such that we have\text{\rm \,:}
\roster
\item"{\bf 10)}"
$\underline{\psi}^{2\nu_p}(\underline{v}_{(2\nu_p)}^*)=p$.
\item"{\bf 11)}"
The map $\underline{\psi}^{2\nu_p}$ is of constant rank 
$m+e_1+\cdots+e_{\kappa_p}$
near $\underline{v}_{(2\nu_p)}^*$ in $\Cal V^*$.
\item"{\bf 12)}"
$\underline{\psi}^{2\nu_p} (\Cal V^*)$ is a manifold-piece $\Cal O^{i_c}$ of the
intrinsically complexified CR orbit of $M$ through $p$.
\item"{\bf 13)}"
$m+e_1+\cdots+e_{\kappa_p}=\text{\rm
dim}_{\C} \, \Cal O^{i_c}=
\text{\rm dim}_{CR} \, \Cal O+\codim_{CR} \, \Cal O$.
\endroster
\endproclaim

\demo{Proof}
Recall that in view of Theorem~7.5, there exists
$\underline{w}_{(2\mu_p-1)}^*\in (\delta\Delta^m)^{2\mu_p-1}$ with
$\underline{\Gamma}_{2\mu_p-1}(\underline{w}_{(2\mu_p-1)}^*)=p^c$,
such that $\Gamma_{2\mu_p-1}$ is of rank $2m+e_1+\cdots+e_{\kappa_p}$
at $\underline{w}_{(2\mu_p-1)}^*$. Looking again at eq.~\thetag{7.17}
(for $k=2j+1$ odd, which we have not written, but the corresponding
equation is similar), and using the chain
rule, we deduce that $\underline{\psi}^{2\mu_p-2}$ is of rank
$m+e_1+\cdots+e_{\kappa_p}$ at the point
$\underline{v}_{2\nu_p}^*:=(\underline{w}_1^*,\ldots,
\underline{w}_{2\mu_p-2}^*)$ and that
$\underline{\psi}^{2\nu_p}(\underline{v}_{(2\nu_p)}^*)=p$ (recall
$\nu_p=\mu_p-1$). This yields {\bf 10)} and {\bf 11)}. For reasons of
dimension, we already know that $\dim_\C \, \Cal O^{i_c}$ must be
equal to $m+e_1+\cdots+e_{\kappa_p}$, since $\dim_{CR} \, \Cal O=m$
and $\dim_\C \, \Cal O^c=m+\dim_\C \, \Cal O^{i_c}$, quite generally.
This is {\bf 13)}. Finally, to deduce {\bf 12)}, we claim that it
can be observed that the range of $\underline{\psi}^{2\nu_p}$ is {\it
a priori} contained in $\Cal O^{i_c}$, and afterwards for dimensional
reasons, the image $\underline{\psi}^{2\nu_p}(\Cal V^*)$ will 
necessarily be a manifold-piece of $\Cal O^{i_c}$ through $p$. 
To complete this observation, we introduce holomorphic coordinates
$(w,z_1,z_2)\in \C^m\times \C^{e_1+\cdots+e_{\kappa_p}}\times 
\C^{n-m-e_1-\cdots-e_{\kappa_p}}$ vanishing at 
$p$ in which the equation of $\Cal O^{i_c}$ 
is $\{z_2=0\}$ and $T_0M=\C_w^m\times \R_x^d$, which is possible.
Using the assumption that $M\cap \{z_2=0\}$ is smooth and of CR
dimension $m$, one shows that the equations of $\Cal M$ can then 
be written 
in the form $z_1=\xi_1+i\bar \Theta_1(w,\zeta,\xi_1,\xi_2)$ 
and $z_2=\xi_2[1+i\bar\Theta_2(w,\zeta,\xi_1,\xi_2)]$.
Then an inspection of the inductive construction of the maps
$\underline{\Gamma}_k$ shows that they have range contained in 
$\{z_2=0, \xi_2=0\}$, whence the maps $\underline{\psi}_{2j}$ have
range in $\{z_2=0\}$, as announced.

The proof of Theorem~7.20 is complete.
\qed
\enddemo

\example{Example~7.21}
Looking at the map $\Gamma_4$ in eq.~\thetag{7.9}, we see that the integer
$2\nu_p=2\mu_p-1$ satisfying the assertions {\bf 10)} and {\bf 11)} of
Theorem~7.20 is in general optimal.
\endexample

This last observation finishes up our presentation of the general
theory of Segre chains.  Our next goal will be to give some nontrivial
illustrations and examples of the theory, by specifying what are the
constraints between the Segre multitype and various other CR
invariant numbers in low codimension. Numerous refined examples of CR
manifolds are spread over the literature and follow for instance from
the study of infinitesimal CR automorphisms (especially in the Russian
school) and from the study of normal forms (Chern-Moser, Loboda,
Ezhov, Schmalz, Ebenfelt and others). In the next very 
concrete paragraph, we shall begin with a comparison
between H\"ormander numbers and Segre multitype.

\head \S8. 
H\"ormander numbers, Minimality, holomorphic degeneracy in low
codimension and various examples
\endhead

\subhead 8.1.~Preliminaries\endsubhead 
As the Minimality Criterion was first established by means of the
H\"ormander numbers of the system $\{L^1,\ldots, L^m,\bar L^1,\ldots
,\bar L^m\}$ of CR vector fields (which are equal to the H\"ormander
numbers of the complexified system $\{\Cal L^1,\ldots, \Cal
L^m,\underline{\Cal L}^1,\ldots,\underline{\Cal L}^m\}$) and
as an important part of the monograph [BER2] is devoted to the
exposition of the proof of this criterion, we find enough motivation
to ask whether or not these numbers are intrinsically related to the
invariants $\kappa_{p}$ and $e_1,\ldots,e_{\kappa_{p}}$ of the pair of
Segre foliations $\{\Cal F_{\Cal L}, \Cal F_{\underline{\Cal
L}}\}$. One could say that the H\"ormander numbers are {\it differential}
in nature wheread the Segre multitype is {\it geometric} in nature.
We need at first some definition and notation.  If $\Cal D=\{\Cal
L^1,\ldots,\Cal L^m,\underline{\Cal L}^1,\ldots,
\underline{\Cal L}^m\}$ is our collection of complexified CR 
{\it $1$-vector fields} (which is of cardinal $2m$), we denote for
$k\geq 1$ by $\Cal D^k$ the set of vector fields including $\Cal D$
and the multiple Lie brackets of length $\leq k$ of the $\Cal L^i$ and
the $\underline{\Cal L}^j$. Let $p\in M$. We assume that $\Cal M$ is
$\Cal D$-minimal at $p^c$. Following [BER2,\S3.4], we define the
integers $h_{p}\in \N_*$ and $\mu_0,\mu_1,\ldots,\mu_{h_{p}} \in \N_*$
called the {\it H\"ormander numbers} of $\Cal D$, together with their
{\it multiplicities} $l_0,l_1,\ldots,l_{h_p}$ as follows.  The number
$\mu_0$ is equal to $1$ and we define another number $l_0:=2m$,
the {\it multiplicity} of $\mu_0$. The number $\mu_1$ is the smallest
integer for which there exists a multiple Lie bracket at $p^c$ of length
$\mu_1$ which is not in the span of $\Cal D(p^c)\subset T_{p^c} \Cal
M$. We define the subspace $\Cal D^{\mu_1}(p^c)\subset T_{p^c}\Cal M$
to be the linear span of $\Cal D(p^c)$ and the values at $p^c$ of all
commutators of vector fields in $\Cal D$ of length $\mu_1$. We define
$l_1$ to be
$$
l_1:= \dim_\C \, \Cal D^{\mu_1}(p^c) - 2m.
\tag 8.2
$$
We define inductively the numbers $\mu_1<\mu_2<\cdots<\mu_{h_{p^c}}$
and special linear subspaces $\Cal D(p^c)=\Cal D^{\mu_0}(p^c)\subset \Cal
D^{\mu_1}(p^c) \subset \cdots \subset \Cal D^{\mu_{h_{p}}}(p^c)\subset
T_{p^c}\Cal M$ as follows. The number $\mu_{k+1}$ is obtained as the
smallest integer for which there exists a multiple Lie bracket at $p^c$ of
vector fields in $\Cal D$ of length $\mu_{k+1}$ which does not belong
to $\Cal D^{\mu_k}(p^c)$. The subspace $\Cal D^{\mu_{k+1}}$ is then defined
as the span of $\Cal D^{\mu_k}(p^c)$ and the values at $p^c$ of all
multiple Lie brackets of vector fields in $\Cal D$ of length
$m_{k+1}$. We define
$$
l_{k+1}:= \dim_\C \, \Cal D^{\mu_{k+1}}(p^c) - \dim_{\C} \, 
\Cal D^{\mu_k}(p^c) =
\dim_\C \, \Cal D^{\mu_{k+1}}(p^c)-2m-\sum_{i=1}^k l_i.
\tag 8.3
$$
It is clear that this process terminates after a finite number of
steps and that $\Cal D^{h_{p}}(p^c)=T_{p^c} \Cal M$ by the minimality
assumption. The number $l_k$ is called the {\it multiplicity} of the
H\"ormander number $\mu_k$. Now, the question arises whether there
exist some peculiar links between the H\"ormander numbers of $\Cal D$
together with their multiplicities and the Segre multitype of $\Cal
D$. Of course, we must have
$$
2m+l_1+\cdots+l_{h_{p}}=2m+e_1+\cdots+e_{\kappa_{p}}=n.
\tag 8.4
$$
Sometimes, in order to insist on the dependence of these numbers on
the reference point $p$, where $p\in M$, we shall write them $e_k(p)$,
$\mu_k(p)$, $l_k(p)$, as we have already written $\kappa_{p}$,
$\mu_{h_{p}}$. In fact, it can be established that for a
Zariski-generic point $p\in M$, the H\"ormander numbers
$\mu_k(p)$ are all equal to $k$ and also, their multiplicities $l_k(p)$ are
constant in a neighborhood of $p$. 
We shall not enter into the details and simply admit here this
property. Of course, the same statement holds for the numbers
$e_k(p)$. Thus there exists a {\it common} Zariski-open subset
$U_{gen}$ of $M$ such that for $p\in U_{gen}$, then all these numbers
are constant equal to $h_{gen}$, $\kappa_{gen}$, $\mu_{k,gen}(=k)$,
$l_{k,gen}$, $1\leq k\leq h_{gen}$ and $e_{k,gen}$, $1\leq k\leq
\kappa_{gen}$.  Finally, one can check that the dependence of all
these numbers is lower semi-continuous with respect to the point $q\in
M$. We shall use these properties in our discussion of some examples
below. As a matter of fact, we shall assume throughout the remainder
of \S8 that all our manifolds $M$ are {\it minimal at every point}. If
$(M,p)$ is minimal, this can always be assumed after picking a sufficiently
small manifold-piece of $M$ around $p$.

\subhead 8.5.~Comparison between H\"ormander numbers and Segre 
multitype \endsubhead Af\-ter these preliminaries are finished, we now
claim that, apart some very peculiar relations which hold in small CR
dimension and small codimension, then eq.~\thetag{8.4} is the only
general relation between the H\"ormander numbers of $\{\Cal L,
\underline{\Cal L}\}$ together with their multiplicities and the Segre
multitype of $\{\Cal L, \underline{\Cal L}\}$, {\it even between the
generic ones}.  This will show that in the study of the CR invariants
(and normal forms) of a real analytic CR manifold {\it at a generic
point}, the generic Segre multitype should add some new independent
information.  Let us argue the assertion of independence with some
examples. First, the minimal hypersurface in $\C^2$ given by $z=\bar
z+i w^{\tilde{\mu}}\bar w^{\tilde{\mu}}$ shows that
$\mu_1=2\tilde{\mu}$ can be arbitrarily large whereas
$e_1(0)=e_{1,gen}=1=l_1(0)=l_{1,gen}$.  As we have already noted, the
hypersurface case ($d=1$) is very particular.  We have observed that
$e_1(p)=1$ for all $p\in M$ (recall that $M$ is everywhere minimal) and of
course also $l_1(p)=1$, $\forall \, p\in M$, whence
$e_{1,gen}=l_{1,gen}=1$. In fact, this relation also follows
immediately from equation~\thetag{8.4}. Now, the $2$-codimensional
generic manifold in $\C^4$ given by $z_1=\bar z_1+iw_1\bar w_1$ and
$z_2=\bar z_2+iw_2^2\bar w_2^2$ has $m=2$, $d=2$, $e_1(0)=2$ but
$l_2(0)=1\neq e_1(0)$. However, $e_{1,gen}=l_{1,gen}=2$ for it. In the
sequel, we shall examine more thoroughly the case $m=d=2$.  After an
inspection of the three different types of CR manifolds of codimension
two in $\C^4$ with nondegenerate Levi-form given by Loboda [Lo] ({\it
see} eqs.~\thetag{8.18-19-20} below), ones sees that $e_{1,gen}=2$ if
$l_{1,gen}=2$. To end-up the comparison between $l_{1,gen}$ and
$e_{1,gen}$, we shall produce a nontrivial example with $l_{1,gen}=1$
but $e_{2,gen}=2$ ({\it see} \S8.7).

On the other hand, in CR dimension $m=1$ and general codimension $d$,
we have seen that the numbers $e_k(p)=e_{k,gen}$ are all equal to $1$
and that $\kappa_{p}=\kappa_{gen}=d$, whereas there is no reason why
all the numbers $l_k$ should be equal to $1$. In CR dimension $m=1$,
after inspecting the multiple Lie brackets of lengths $\leq 3$, one
can easily see that $l_{2,gen}=1$ and $l_{3,gen}=1$ necessarily,
because there are two Lie brackets of length $2$ of the system
$\{L^1,\bar L^1\}$ ($m=1$) which are opposite to each other and there
are four multiple Lie brackets of length $3$ which are either opposite
or conjugate to each other. However there may exist two multiple Lie
brackets of lentgh $4$ which can really differ from each other, thus
making possible $l_{4,gen}=2$. Consequently, the codimension $d\geq 4$
is necessary to construct an example with $l_{3,gen}=2$ (whereas
$e_{3,gen}=1$).  Here is such an example in
$\C_{(w,z_1,z_2,z_3,z_4)}^5$\,:
$$
M:
\left\{
\aligned
&
z_1=\bar z_1+i w\bar w, \ \ \ \ \ \ \ \ \ \ \ \ \ \ \ \ \ \ 
z_2=\bar z_2+i w \bar w(w+\bar w)\\
&
z_3=\bar z_3+i w \bar w(w^2+\bar w^2), \ \ \ \ \
z_4=\bar z_4+iw^2\bar w^2.
\endaligned\right.
\tag 8.6
$$

\subhead 8.7.~Characterization of holomorphic degeneracy by Segre
multitype in codimension $2$ and CR dimension $2$ \endsubhead As announced,
we study in this paragraph
the case $d=2$, $m= 2$, with $M$ minimal at {\it
every point}. If $M$ is connected (which we will also suppose), then
there exists an integer $\kappa_M$ with $0\leq \kappa_M
\leq 2$ such that $M$ is biholomorphic to a product $\underline{M}
\times \Delta^{\kappa_M}$ in a neighborhood of a Zariski-generic point
$p\in M$, {\it see} [BER1,2], [M2]. We say that $M$ is {\it
$\kappa_M$-holomorphically degenerate}. If $\kappa_M=2$, then
$\underline{M}$ is a maximally real submanifold of $\C^2$ and $M$ is
biholomorphic to $\R^2\times \Delta^2$ near such Zariski-generic
points.  In fact, this implies that $M$ is Levi-flat, and equivalent
to $\R^2\times \Delta^2$ in a neighborhood of {\it every} point.  Thus
$e_1(p)=0$, $\forall \ p\in M$, the CR-orbits are the leaves
$x_p\times \Delta^2$ and $M$ cannot be minimal, contradiction. On the
other hand, if $\kappa_M=1$, $M$ can well be minimal at every point,
and it is easy to see that for each Zariski-generic point $p\in M\cong
\underline{M}\times \Delta$, then $e_1(p)=e_2(p)=1$ and $e_3(p)=0$,
since $d=2$.  By the upper semi-continuity of the numbers $e_k(q)$
with respect to $q\in M$, it follows that $e_1(p)=e_2(p)=1$ for {\it
all} points $p\in M$. In fact, we shall prove that the property
$e_1(p)=1$ $\forall \ p\in M$, for $M$ minimal at every point of
codimension $d=2$ and CR dimension $m\geq 2$, {\it characterizes}
$1$-holomorphic degeneracy. This will show an interesting and
unexpected illustration of the Segre multitype and will provide some
new nontrivial examples.  By the way, before proceeding to the
codimension two case, we would like to remind that in codimension
$d=1$ and CR dimension $m=1$, it is already well know that an
everywhere minimal $M\subset \C^2$ is holomorphically nondegenerate
($\kappa_M=0$).  The following proposition generalizes such an
observation to codimension $2$ and CR dimension $2$. Its interest lies
in the analysis of the local form of an arbitrary generic
six-dimensional $M\subset \C^4$ that we shall delineate during the
course of its proof.

\proclaim{Proposition~8.8}
Let $M\subset \C^n$ be a connected CR-generic real analytic manifold of
codimension $d=2$ and CR dimension $m= 2$ which is minimal at every
point. Then $1\leq e_1(p)\leq 2$ for all
$p\in M$ and the following properties are equivalent
\roster
\item"{\bf (a)}"
$e_{1,gen}=1$.
\item"{\bf (b)}"
$e_1(p)=1$ for all $p\in M$.
\item"{\bf (c)}"
$M$ is $1$-holomorphically degenerate.
\endroster
Consequently $M$ 
is holomorphically nondegenerate if and only if
$e_{2,gen}=2$.
\endproclaim

\remark{Remark}
According to [BER], $M$ is holomorphically nondegenerate if
and only if its {\it Levi type} $\ell$ is finite at a generic point 
and then $1\leq \ell_{gen}\leq m$. By definition, the Levi-type at $p\in M$
is the smallest integer $k$ such that
$$
\text{\rm Span} \ 
\{\bar L^\beta \nabla \rho_j(p,\bar p)\: \v \beta \v \leq k, 1\leq j\leq
d\}=\C^n,
\tag 8.9
$$
where $\nabla\rho_j$ denotes the holomorphic gradient $\partial
\rho_j/\partial t$ of the $j$-th Cartesian equation of $M$. 
One says that $M$ is {\it $k$-nondegenerate at $p$} if
eq.~\thetag{8.9} holds.  Essentially two different classes of CR
manifolds $M$ arise in Proposition~8.8 according to whether
$\ell_{gen}=1$ or $\ell_{gen}=2$. As in the study by Ebenfelt [E] of
normal forms for hypersurfaces in $\C^3$, it holds for $m=d=2$
that $\ell_{gen}=1$ if and only if $M$ is Levi degenerate at a generic
point if and only if $l_{1,gen}=1$.
\endremark

\example{Example~8.10}
The generic manifold $2y_1=w_1\bar w_1$, $2y_2= x_1^2w_2\bar
w_2+w_1^2\bar w_1^2 w_2\bar w_2/4$, or equivalently $ M\: \ z_1=\bar
z_1+i w_1\bar w_1,\ z_2=\bar z_2+i \bar z_1 [iw_1\bar w_1w_2\bar w_2+
\bar z_1 w_2\bar w_2]$,
is minimal at $0$, holomorphically nondegenerate, hence satisfies
$e_{1,gen}=2$ by the above Proposition~8.8 (this can be checked
directly), but here we have $e_1(0)=1$ (true semi-continuity).
Consequently, the conditions {\bf (a)} and {\bf (b)} cannot be
weakened. Inspired by eq.~\thetag{8.28} below, we can further provide
such an example with a {\it rigid} $M$ which is $2$-nondegenerate at
$0$, minimal at $0$ and satisfies $e_{1,gen}=2$ whereas $e_1(0)=1$
only\,:
$$
\left\{
\aligned
&
z_1=\bar z_1+i[w_1\bar w_1+w_1^2\bar w_2+\bar w_1^2w_2]\\
&
z_2=\bar z_2+i[w_1\bar w_1(w_1+\bar w_1)-w_2\bar w_1^2(2w_1+\bar w_1)-
\bar w_2w_1^2(2\bar w_1+w_1)].
\endaligned\right.
\tag 8.11
$$
\endexample

\example{Example~8.12} 
Our goal now is to produce an example of $M$ with $m=d=2$ and
$l_{1,gen}=1$ which is minimal and holomorphically nondegenerate
(exercise\,: verify that simple polynomial rigid examples exist if we
drop one of the last two assumptions).  First, we would like to remind
that, according to the simplification by Freeman of a theorem of
Sommer, the kernel of the Levi form of any $\Cal C^\omega$ (even $\Cal
C^2$) CR-generic $M$ is {\it involutive} locally where it is of
constant rank $e$, from which follows that $M$ is locally {\it foliated} by
complex
manifolds of dimension $e$ ({\it see} [F1], [Chi],~p.152).  An
important example of a real analytic everywhere Levi degenerate but
{\it not} holomorphically degenerate hypersurface is the set of
regular points of the cubic discovered by Freeman [F]\,:
$x_1^3+x_2^3+x_3^3=0$ in $\C_{z_1,z_2,z_3}^3$. Freeman sought
hypersurfaces with everywhere degenerate Levi form which are not
locally biholomorphic to a product by $\Delta$, or equivalently in
modern terms, which are not holomorphically degenerate. In fact, the
same property holds true for the (quite simpler, because of degree
two) hypersurface $x_1^2+x_2^2-x_3^2=0$, called the {\it light cone},
which is simply the tube manifold over the classical real cone in
$\R^3$. The light cone is a homogeneous CR manifold, whereas Freeman's
cubic is not. Both are rigid, because of the existence
of the obvious CR-transversal infinitesimal CR automorphisms $\text{\rm
Re} \, (i\partial / \partial z_j)$, $1\leq j\leq 3$. Ebenfelt [E]
shows that the light cone is finitely nondegenerate (and in
fact $2$-nondegenerate) at every point with exactly one nonzero 
eigenvalue for the Levi form and shows that its equations at 
a generic point (in fact at every point) of the form~\thetag{5.7}
can be written as\,: $z=\bar z+i [w_1\bar w_1+w_1^2\bar w_2+\bar
w_1^2w_2+O_{weighted}(4)]$ ([E], p.~318, eq.~\thetag{A.i.2}; the same
property holds true for Freeman's cubic outside the coordinate
planes).  Building on these two
objects, we can construct as follows an example of
$M=\Sigma_{reg}\subset \C^4_{z_1,z_2,z_3,z_4}$,
with $l_{1,gen}=1$ which is minimal and $2$-nondegenerate at every
point, where\,:
$$
\Sigma:
\left\{
\aligned
&
x_1^2+x_2^2-x_3^2=0,\\
&
x_1^3+x_2^3-x_4^3=0.
\endaligned\right.
\tag 8.13
$$
Here, $\Sigma$ has a complex foliation by complex lines described 
in terms of the holomorphic leaf variable $\zeta\in \C$ and of
four real parameters $a,b,c,d$ as follows\,:
$$
(\zeta,a,b,c,d)\mapsto (a(1+\zeta)+ib,1+\zeta,
(1+a^2)^{1/2}(1+\zeta)+ic,
(1+a^3)^{1/3}(1+\zeta)+id).
\tag 8.14
$$
({\it cf.} [F2]). Over the subset $M:=\Sigma_{reg}=
\Sigma\backslash\{0\}$ of regular
points of $\Sigma$, the complex tangent bundle $T^cM$ admits the
following two spanning global sections (notice that $(x_1,x_2)\neq (0,0)$, 
$x_3\neq 0$ and $x_4\neq 0$ over $M$)\,:
$$
\left\{
\aligned
&
L_1={\partial \over \partial z_1}+
{x_1\over x_3} \ {\partial \over \partial z_3} + {x_1^2 \over x_4^2} \
{\partial \over \partial z_4},\\
&
L_2={\partial \over \partial z_2}+
{x_2\over x_3} \ {\partial \over \partial z_3} + {x_2^2 \over x_4^2} \
{\partial \over \partial z_4}.
\endaligned\right.
\tag 8.15
$$
An isotropic (for the $\R^2$-valued Levi form of $M$) 
vector field is simply the radial vector field, which is
nowhere vanishing on $M$\,:
$$
L=x_1{\partial \over \partial z_1}+x_2 {\partial \over \partial z_2}
+x_3{\partial \over \partial z_3}+x_4{\partial \over \partial z_4}, 
\ \ \ \ \ L\v_M=x_1L_1\v_M+x_2L_2\v_M.
\tag 8.16
$$
Obviously, the 
CR automorphism of $M$ generated by $L$ consists of the real dilatation
$z\mapsto \lambda z$.
By direct computations we have\,:
$$
[L,\bar L]=-\bar L+L, \ \ \ \ \
[L,\bar L_1]=-L_1, \ \ \ \ \ 
[L,\bar L_2]=-L_2,
\tag 8.17
$$
which shows that the kernel of the Levi form contains $L$. Then the
computation of $[L_1,\bar L_1]$ and $[L_2,\bar L_2]$ shows that the
Levi-form of $M$ has exactly one zero and one nonzero eigenvalue at
every point. Consequently, $l_{1,gen}=1$. Further, 
it is easy to show that $M$
is $2$-nondegenerate at every point ({\it cf.} [E]) by checking that
the five vectors $\rho_1$, $\bar L_1(\nabla\rho_1)$, 
$\bar L_2(\nabla\rho_1)$, $\bar L_1\bar
L_1(\nabla\rho_1)$ and $\bar L_2\bar L_2(\nabla\rho_1)$ are linearly
independent at every point of $M$, where $\nabla\rho_1=(\rho_{1,z_1},
\rho_{1,z_2},\rho_{1,z_3})$ is the complex gradient of
$\rho_1:=x_1^2+x_2^2-x_3^2$.  Finally, computing higher order Lie
brackets, one sees that $M$ is minimal at every point, and more
precisely that $\mu_1(p)=2$, $l_1(p)=1$, $\mu_2(p)=3$, $l_2(p)=1$ if
$x_{1,p}x_{2,p}\neq 0$ and that $\mu_1(p)=2$, $l_1(p)=1$,
$\mu_2(p)=4$, $l_2(p)=1$ if $x_{1,p}=0$ or $x_{2,p}=0$
(but $(x_{1,p},x_{2,p})\neq (0,0)$).
\endexample

\demo{Proof of Proposition~8.8}
We have already understood the implications {\bf (a)}
$\Longleftrightarrow$ {\bf (b)} and {\bf (c)} $\Longrightarrow$
{\bf (a)}. Thus, it remains to establish {\bf (b)} $\Longrightarrow$
{\bf (c)}. We shall proceed to another implication, which we now
know to be the contraposition of {\bf (b)} $\Longrightarrow$ {\bf
(c)}, namely we shall prove that if $M$ is holomorphically
nondegenerate and minimal everywhere, then $e_{2,gen}=2$. In fact, we
shall see that the two hand-sides of the last equivalence stated in
Proposition~8.8 come
together in our discussion. To do so, we analyze first a generic point
of minimality of $M$ in terms of Lie brackets. By [BER1,2], the
Levi-number of $M$ is $\ell_M=2$, whence a generic point of
holomorphic nondegeneracy of $M$ is {\it $2$-nondegenerate} and
morevover, it has the property that, in regular coordinates, the
homogeneous Taylor approximation of $M$ of order $\leq 3$ with respect
to $(w,\bar w)$ is already $2$-nondegenerate. We assume that such
a choice of generic point for this property has been made for the
reference point $0\in M$ and we also assume that the H\"ormander
numbers of $M$ at $0$ are the generic ones.  For combinatorial
reasons, there are exactly two possibilities for the generic
H\"ormander numbers of $M$. Either they are equal to $\mu_{1,gen}=2$,
$\mu_{2,gen}=3$, with multiplicities $l_{1,gen}=1$ and $l_{2,gen}=1$ or
$\mu_{1,gen}=2$ with multiplicity $l_{1,gen}=2$. In the latter case,
it is easy to deduce from $l_1(p)=2$ that the two Levi forms of the
two defining equations of $M$ are linearly independent and that the
intersection of their kernel is null.  This assumption corresponds
exactly to the {\it CR-generic manifolds with nondegenerate vector-valued
Levi-form}, in the sense of the Russian school. It is known that such
manifolds can be reduced to three (inequivalent) types ({\it see}
[Lo], [Belo], [EIS]), the {\it elliptic} type\,:
$$
M:
\left\{
\aligned
&
z_1=\bar z_1+i[w_1\bar w_1+O_w(3)+O_{\bar z}(1)]\\
&
z_2=\bar z_2+i[w_2\bar w_2+O_w(3)+O_{\bar z}(1)],
\endaligned\right.
\tag 8.18
$$
or the {\it parabolic} type\,:
$$
M:
\left\{
\aligned
&
z_1=\bar z_1+i[w_1\bar w_1+O_w(3)+O_{\bar z}(1)]\\
&
z_2=\bar z_2+i[w_1\bar w_2+\bar w_1 w_2+O_w(3)+O_{\bar z}(1)],
\endaligned\right.
\tag 8.19
$$
or the {\it hyperbolic} type\,:
$$
M:
\left\{
\aligned
&
z_1=\bar z_1+i[w_1\bar w_1-w_2\bar w_2+O_w(3)+O_{\bar z}(1)]\\
&
z_2=\bar z_2+i[w_1\bar w_2+\bar w_1 w_2+O_w(3)+O_{\bar z}(1)].
\endaligned\right.
\tag 8.20
$$
Direct computation shows that such manifolds are minimal,
holomorphically nondegenerate (in fact their Taylor approximation of
order $2$ is already $1$-nondegenerate) and
that $e_{1,gen}=e_1(0)=2$. Using the representation eq.~\thetag{7.7},
one sees that $e_1(0)=2$ if and only if the following determinant
of order one partial derivatives of the two defining functions
$\Theta_1$, $\Theta_2$, restricted to the second Segre chain 
$S_0^2$\,:
$$
\text{\rm det} \left( \! \! \! \! \!
\matrix
& \Theta_{1,w_1}(\zeta_1,\zeta_2,w_1,w_2,0,0) & 
\Theta_{1,w_2}(\zeta_1,\zeta_2,w_1,w_2,0,0)\\
& \Theta_{2,w_1}(\zeta_1,\zeta_2,w_1,w_2,0,0) & 
\Theta_{2,w_2}(\zeta_1,\zeta_2,w_1,w_2,0,0)
\endmatrix
\right)\not\equiv 0,
\tag 8.21
$$
does not vanish identically in $\C\{w_1,w_2,\zeta_1,\zeta_2\}$.  In
this case, the desired property $e_2(0)=2$ comes immediately for each one
of the three types \thetag{8.18-19-20} (just compute~\thetag{8.21}). 
It remains therefore to study the
(more degenerate) case where $\mu_{1,gen}=2$, $l_{1,gen}=1$,
$\mu_{2,gen}=3$, $l_{2,gen}=1$. Centering the
study at a generic point, this case can be reduced {\it in
regular coordinates} to equations of the form
$$
\left\{
\aligned
& 
z_1=\bar z_1+i[w_1\bar w_1+\varepsilon_2 w_2\bar w_2+H_1^3(w,\bar w)+
O(\v w\v^4)+O(\v \bar z\v \, \v w\v^2)]\\
& 
z_2=\bar z_2+i[H_2^3(w,\bar w)+O(\v w\v^4)+O(\v \bar z\v \,  \v w\v^2)],
\endaligned\right.
\tag 8.22
$$
where $H_1^3$ and $H_2^3$ are homogeneous polynomial of degree $3$ in
$(w,\bar w)$ containing no pluriharmonic terms and
where $\varepsilon_2=0, 1$ or $-1$, depending on
the rank of the Levi-form of $\rho_1$ at $0$. Because of regular
coordinates, using eq.~\thetag{5.8} one can see that $H_1^3$ and $H_2^3$ are {\it real}
(however, the next terms homogeneous of order $4$, $H_j^4(w,\bar w)$,
are not in general real). Clearly, we have $l_1(0)=1$ here. Now,
computing multiple Lie brackets of length equal to $3$, one sees that
$l_2(0)=1$ (which we assumed) if and only if $H_2^3(w,\bar
w)\not\equiv 0$. So in the sequel, we will assume $H_2^3(w,\bar
w)\not\equiv 0$. Consequently, we can assign weight $1$ to 
$w$ and to $\bar w$, weight $2$ to $z_1$ and to $\bar z_1$ and
finally weight $3$ to $z_2$ and to $\bar z_2$. In the sequel, 
we shall work modulo $O_{weighted}(k)$, $k=2,3,4,\ldots$, say
for short $O_{wt}(k)$. Now, since
we have assumed $l_{1,gen}=1$, the following $2\times 4$ matrix 
of the vector-valued Levi form of $M$ modulo $O_{wt}(2)$ must be of
rank $1$ everywhere\,:
$$ 
\left( \! \! \! \! \!
\matrix
& 1+H_{1,w_1\bar w_1}^3 & \varepsilon_2+H_{1,w_2\bar w_2}^3 &
H_{1,w_1\bar w_2}^3 & H_{1,\bar w_1w_2}^3 \\
& H_{2,w_1\bar w_1}^3 & H_{2,w_2\bar w_2}^3 &
H_{2,w_1\bar w_2}^3 & H_{2,\bar w_1w_2}^3
\endmatrix\right),
\tag 8.23
$$
again modulo $O_{wt}(2)$. If $\varepsilon_2\neq 0$, we deduce
$H_{2,w_1\bar w_1}^3 \equiv H_{2,w_1\bar w_2}^3 \equiv H_{2,\bar w_1
w_2}^3 \equiv H_{2,w_2\bar w_2}^3\equiv 0$, whence $H_{2}^3(w,\bar
w)\equiv 0$, a contradiction. Therefore, we must have
$\varepsilon_2=0$ and we then deduce $H_{2,w_1\bar w_2}^3 \equiv
H_{2,w_1\bar w_2}^3 \equiv H_{2,\bar w_1w_2}^3\equiv 0$, whence
$H_3^2(w,\bar w)=Aw_1\bar w_1(w_1+\bar w_1)$, with $A\neq 0$.  After a
complex scaling, we can assume $A=1$ and the equations of $M$ are thus
in the form
$$
\left\{
\aligned
& 
z_1=\bar z_1+i[w_1\bar w_1+H_1^3(w,\bar w)+O_{wt}(4)]\\
& 
z_2=\bar z_2+i[w_1\bar w_1(w_1+\bar w_1)+O_{wt}(4)].
\endaligned\right.
\tag 8.24
$$
Now, we pay attention to the supplementary information coming from the
Levi type $\ell_{gen}$ of $M$.  The case $\ell_{gen}=1$ {\it with $M$
minimal} again implies Levi-nondegeneracy (in the sense of Beloshapka
and Loboda) at a generic point (exercise), {\it i.e.} one of the three
forms eq.~\thetag{8.18-19-20}. Thus, we can assume that $\ell_{gen}=2$
and then that $M$ is $2$-nondegenerate at $0$ after a small shift
(delocalization) of the origin. Then if we develope the homogeneous
polynomial $H_1^3$ in the form
$$
H_1^3=2\text{\rm Re} \, \left(Aw_1^2\bar w_1+ Bw_1^2\bar w_2+ Cw_1w_2\bar
w_2+Dw_1w_2\bar w_1+Ew_2^2\bar w_1+
Fw_2^2\bar w_2\right),
\tag 8.25
$$
we observe that $M$ is $2$-nondegenerate at $0$ {\it if and only if}
$B\neq 0$ or $C\neq 0$ or $F\neq 0$. Then coming back to the
assumption $l_{1,gen}=1$, we observe that this implies that
the rank of the Levi-form of the first defining equation $\rho_1$ of 
$M$ must be equal to $1$ at every point of $M$ in a neighborhood of 
$0$ (otherwise, we come again to the case
$\varepsilon_2\neq 0$, which is excluded). Looking at the first row of
matrix~\thetag{8.23} (with now $\varepsilon_2=0$) and inspecting
the Levi-form of $\rho_1$ modulo $O_{wt}(2)$, we see that necessarily 
$H_{1,w_2\bar w_2}^3(w,\bar w)=Cw_1+\bar C\bar w_1+2Fw_2+
2\bar F \bar w_2\equiv 0$, whence $C=F=0$. After a scaling and a
linear change of coordinates in the $z$-space, we can 
assume that $A=0$ and $B=1$. Then replacing $w_1$ by 
$w_1+Dw_1w_2+Ew_2^2$, we can assume that $D=E=0$ and in 
conclusion, the equations of $M$ are reduced to the following
form\,:
$$
M:
\left\{
\aligned
&
z_1=\bar z_1+i[w_1\bar w_1+w_1^2\bar w_2+\bar w_1^2w_2+
O_{wt}(4)]\\
&
z_2=\bar z_2+i[w_1\bar
 w_1(w_1+\bar w_1)+O_{wt}(4)].
\endaligned\right.
\tag 8.26
$$
We notice that the cubic tangent to eqs.~\thetag{8.26} is
$2$-nondegenerate, and minimal at $0$. Further, it is in fact
Levi-nondegenerate at a generic point (!). However, order 4 terms have
an influence to make possible the property 
$l_{1,gen}=1$, as in Example~8.12, a property which
we are assuming since a while. Thus, continuing the proof of
Proposition~8.8, we have to explicit the order 4 terms in
eq.~\thetag{8.26}, which we now write as follows\,:
$$
M:
\left\{
\aligned
&
z_1=\bar z_1+i[w_1\bar w_1+
w_1^2\bar w_2+\bar w_1^2 w_2+\bar H_1^4(w,\bar w)+\\
&
\ \ \ \ \ \ \ \ \ \ \ \ \ \
+\bar z_1 (\alpha_1 w_1\bar w_1+
\beta_1 w_1\bar w_2+\gamma_1 w_2\bar w_1+\delta_1 w_2\bar w_2)+
O_{wt}(5)],\\
&
z_2=\bar z_2+i[w_1\bar w_1(w_1+\bar w_1)+
\bar H_2^4(w,\bar w)+\\
&
\ \ \ \ \ \ \ \ \ \ \ \ \ \
+\bar z_1 (\alpha_2 w_1\bar w_1+
\beta_2 w_1\bar w_2+\gamma_2 w_2\bar w_1+\delta_2 w_2\bar w_2)+
O_{wt}(5)],
\endaligned\right.
\tag 8.27
$$
where $H_1^4$ and $H_2^4$ are homogeneous polynomials of order four in
$(w,\bar w)$ containing no pluriharmonic terms. To finish the proof of
Proposition~8.8, we claim that $e_2(0)=2$ in the obtained form
eq.~\thetag{8.27} if furthermore $l_{1,gen}=1$.
We shall proceed by contradiction.
Thus, we assume that $e_2(0)=1$, {\it i.e.} that the 
determinant in eq.~\thetag{8.21} vanishes identically.
First, we deduce after inspecting the order three
terms in this determinant that we have
$$
(2w_1\zeta_1+\zeta_1^2) \,
\zeta_1-
\bar H_{2,w_2}^4(w,\zeta)\equiv 0.
\tag 8.28
$$
However, we have not yet expressed the assumption $l_{1,gen}=1$ (which
will contradict eq.~\thetag{8.28}). This assumption can be expressed by
saying that the two Lie brackets $[L_1,\bar L_2]$ and $[L_2,\bar L_2]$
are colinear to the Lie bracket $[L_1,\bar L_1]$ at every nearby
point. Let us compute these Lie brackets. 
We shall write them in the form\,:
$$
\left\{
\aligned
& 
[L_1,\bar L_1]:= P_1 {\partial \over \partial z_1}+ Q_1 {\partial
\over \partial z_2}+
\bar P_1 {\partial \over \partial \bar z_1}+
\bar Q_1 {\partial \over \partial \bar z_2}
\\
& 
[L_1,\bar L_2]:=R {\partial \over \partial z_1}+ S {\partial \over
\partial z_2}+ T {\partial \over \partial \bar z_1} + U {\partial
\over \partial \bar z_2}\\ 
& 
[L_2,\bar L_2]:= P_2 {\partial \over
\partial z_1}+ Q_2 {\partial \over \partial z_2}+
\bar P_2 {\partial \over \partial \bar z_1}+
\bar Q_2 {\partial \over \partial \bar z_2}
\endaligned\right.
\tag 8.29
$$
From the colinearity of all these three Lie brackets at every point,
we shall extract only one equation, namely\,: $\text{\rm det} \left(
\! \! \! \!\matrix & P_1 & Q_1 \\ & R & S \endmatrix
\right)\equiv 0$. After inspection of eq.~\thetag{8.27}, one sees that
$P_1=-i+O_{wt}(1)$, $Q_1=O_{wt}(1)$, $R=O_{wt}(1)$ and $S=O_{wt}(2)$.
Since we want to extract from this equation only the order two terms,
we need only to compute $Q_1$ and $R$ to order $1$ and $S$ to order
$2$. We have the following complete expression for our two $(1,0)$ CR vector
fields\,:
$$
\left\{
\aligned
& L_1={\partial \over \partial w_1}+ i\left[\bar 
w_1+2w_1\bar w_2+\bar H_{1,w_1}^4(w,\bar w)
+\alpha_1\bar z_1\bar w_1+\beta_1\bar z_1\bar w_2+O_{wt}(4)\right]{\partial
\over \partial z_1}+\\
&
\ \ \ \ +i\left[2w_1\bar w_1+\bar w_1^2+
\bar H_{2,w_1}^4(w,\bar w)+
\alpha_2\bar z_1\bar w_1+\beta_2\bar z_1\bar w_2+O_{wt}(4)\right]
{\partial
\over \partial z_2},\\
& 
L_2={\partial \over \partial w_2}+i\left[\bar w_1^2+\bar H_{1,w_2}^4
(w,\bar w)+
\gamma_1\bar z_1\bar w_1+\delta_1\bar z_1\bar w_2+O_{wt}(4)\right]{\partial
\over \partial z_1}+\\
&
\ \ \ \ +i\left[\bar H_{2,w_2}^4(w,\bar w)+
\gamma_2\bar z_1\bar w_1+\delta_2\bar z_1
\bar w_2+O_{wt}(4)\right]{\partial
\over \partial z_2}
\endaligned\right.
\tag 8.30
$$
Now, we can compute 
$$
\aligned
&
P_1=-i+O_{wt}(1), \ \ \ \ \ \ \ \ \ \ 
Q_1=-i[2w_1+2\bar w_1+O_{wt}(2)],\\
&
R=-i[2w_1+O_{wt}(2)], \ \ \ \ \
S=-i[\bar H_{2,w_1\bar w_2}^4(w,\bar w)+
\beta_2\bar z_1+O_{wt}(3)],
\endaligned
\tag 8.31
$$
and the vanishing of the above $2\times 2$ determinant yields the 
following supplementary partial differential equation for $\bar H_2^4$\,:
$$
\bar H_{2,w_1\zeta_2}^4(w,\zeta)+\beta_2 \bar z_1-
(2w_1+2\zeta_1)(2w_1)\equiv 0.
\tag 8.32
$$
On one hand, we deduce from eq.~\thetag{8.28} that the homogeneous
polynomial $\bar H_2^4(w,\zeta)$ necessarily incorporates the monomial
$w_2\zeta_1^3$. On the other hand, we deduce from eq.~\thetag{8.32},
that $\bar H_2^4(w,\zeta)$ contains necessarily the monomial ${4\over
3}\zeta_2 w_1^3$.  The polynomial $\bar H_2^4(w,\zeta)$ is not real,
but nevertheless, it satisfies the relation
$$
\bar H_2^4(w,\zeta)-H_2^4(\zeta,w)-iw_1\zeta_1
[\bar \alpha_2 \zeta_1w_1+\bar \beta_2 \zeta_1w_2+
\bar \gamma_2 \zeta_2 w_1+\bar\delta_2 \zeta_2w_2] \equiv 0,
\tag 8.33
$$
which follows from an inspection of the terms of weighted oder $\leq
4$ in eqs.~\thetag{8.27} and~\thetag{5.8}.  Because of
eq.~\thetag{8.8}, the monomial ${4\over 3}\zeta_2 w_1^3$ should be the
conjugate of the monomial $w_2\zeta_1^3$, but ${4\over 3}\neq 1$. This
is the desired contradiction, which completes the proof of
Proposition~8.8 is complete.
\qed
\enddemo

Along with other results, Ebenfelt in [E] obtains essentially three
forms {\it at a Zariski-generic point} for a minimal and
holomorphically nondegenerate hypersurface in $\C^3$\,: $z_1=\bar
z_1+i[w_1\bar w_1+w_2\bar w_2+O_{wt}(3)]$ or $z_1=\bar z_1+i[w_1\bar
w_1-w_2\bar w_2+O_{wt}(3)]$ in the classical Levi-nondegenerate cases,
and the last form $z_1=\bar z_1+i[w_1\bar w_1+w_1^2\bar w_2+
\bar w_1^2w_2+O_{wt}(4)]$ in the everywhere Levi-degenerate case 
$l_{1,gen}=1$. Incidentally, our analysis has provided a similar list.

\proclaim{Corollary~8.34}
There are exactly four nonequivalent local Taylor representations 
at a Zariski-generic point for an arbitrary real analytic CR-generic 
manifold of codimension $2$ in $\C^4$. If $M$ is Levi-nondegenerate, 
then it is either elliptic~\thetag{8.18} or parabolic~\thetag{8.19}
or hyperbolic~\thetag{8.20}. If $M$ is everywhere Levi-degenerate, 
then it is representable in the form~\thetag{8.26} with further 
conditions as {\it e.g.}~eq.~\thetag{8.32}.
\endproclaim

A quite general problem in CR geometry would be to devise a complete
classification of local Taylor-approximated representations for arbitrary
real analytic CR manifolds in any CR dimension and in any
codimension. The author ignores whether or not this question if
pure Utopia.

\head \S9. Orbits of systems of holomorphic vector fields and 
a refinement of Sussmann's theorem\endhead

This appendix paragraph exhibits the central notion of {\it orbits of
bundles of holomorphic vector fields} in a self-contained way.  As we
have been inspired by the ge\-ne\-ral constructions in [Sus], we give
here a brief generalization of Theorem~7.5.  The notion of orbits of
families of vector fields goes back to a well written paper of
Sussmann [Sus], so we will refer the reader to it for background and
further information. Sussmann considered only the $\Cal C^{\infty}$
case, but his construction works as well in $\Cal C^2$, $\Cal C^k$,
$\Cal C^{\omega}$ or in the complex analytic category. Accordingly, we
will in this section give a proof of the Orbit Theorem about integral
submanifolds of bundles of holomorphic vector fields in the spirit of
[Sus], but with a supplementary important simplification due to the
principle of analytic continuation.

\subhead 9.1.~Flows of vector fields\endsubhead \!
Let $\Delta$ be the unit disc in $\C$ and $r\Delta=\{|z|<r\}$. Let
$\SS=\{L_{\alpha}\}_{\alpha\in A}$, $1\leq
\alpha \leq a$, $A=[\![1,a]\!]$, $a\in \N$, $a\geq 1$, be a {\it
finite} system of nonzero $m$-vectorial holomorphic vector fields over
$\Delta^n$, $m\in \N$, $m\geq 1$. By $m$-vectorial, we mean that each
$L_{\alpha}$ is a collection $(L_{\alpha 1},\ldots,L_{\alpha m})$ of
$m$ commuting linearly independent over $\Delta^n$ vector fields.
Hence, considering the {\it multiple flow mapping} $\C^m\times
\Delta^n \ni(s_1,\ldots,s_m,z) \mapsto 
\exp(s_mL_{\alpha m}) \circ \cdots \circ
\exp(s_1L_{\alpha 1})(z)\in \Delta^n$ (which is
defined on a certain domain), we have for every permutation $\varpi :
[\![1,m]\!]\to [\![1,m]\!]$\,:
$$
\exp(s_{\varpi(m)} L_{\alpha \varpi(m)})\circ
\cdots \circ \exp (s_{\varpi(1)} L_{\alpha \varpi(1)})(z)=
\tag 9.2
$$
$$ \ \ \ \ \ \ \ \ \ \ \ \ \ \ \ \ \ \
=\exp(s_mL_{\alpha m}) \circ \cdots \circ
\exp(s_1L_{\alpha 1})(z).
$$
We shall simply denote this multiple flow map by $(s,z)\mapsto
L_{\alpha s}(z)$ and we will work with multiple vector fields {\it
formally} as if they were {\it usual} vector fields, {\it i.e.} as if
$m=1$. But $s$ will be called an $m$-time and $s\mapsto L_{\alpha
s}(z)$ an $m$-curve. We recall the defining properties of the flow
map\,: $L_{\alpha 0}(z)=z$ and $\frac{d}{ds} (L_{\alpha
s}(z))=L_{\alpha}( L_{\alpha s}(z))$, where $L_{\alpha}(z')$ denotes
the value of $L_{\alpha}$ at $z'$, an $m$-vector in
$T_{z'}\Delta^n$. Now, we assume that $\text{\rm rk}_{\C}
(L_1(p),\ldots,L_a(p))=am$ all over $\Delta^n$, in other words, we
assume that the span over $\Cal O(\Delta^n)$ of the $L_{\alpha}$'s
generates a trivial $am$-dimensional holomorphic vector bundle on
$\Delta^n$. Put $d=n-am$, the {\it codimension} of $\SS$. Let $0<r\leq
1/2$ and $(r\Delta)^n \subset \subset
\Delta^n$. Our aim is to define finite concatenation of
flow mappings of such $m$-vector fields.  If $k\in \N_*(:=\N\backslash
\{0\}$), $\L^k=(L^1,\ldots,L^k)\in \SS^k$,
$t_{(k)}=(t_1,\ldots,t_k)\in
\C^{mk}$, $z\in (r\Delta)^n$, we use the notation
$\L_{t_{(k)}}^k(z)=L_{t_k}^k\circ\cdots\circ L_{t_1}^1(z)$ whenever
the composition is defined. Anyway, after bounding $k\leq 3n$, it is
clear that there exists $\delta >0$ such that all maps
$(t_{(k)},z)\mapsto
\L_{t_{(k)}}^k(z)$ are well-defined
for $t_{(k)}\in (2\delta \Delta^m)^k$, $z\in(\frac{r}{2}\Delta)^n$ and
satisfy $\L_{t_{(k)}}^k(z)\in (r\Delta)^n$, for all $t_{(k)}\in
(2\delta\Delta^m)^k$, $z\in (\frac{r}{2}\Delta)^n$, $k\leq 3n$,
$\L^k\in \SS^k$.
The {\it geometric interpretation} is that the maps $t\mapsto L_{\alpha
t}(z)$ are integral $m$-curves of the element $L_{\alpha}\in \SS$ and
more generally, that the
maps $t_{(k)}\mapsto \L_{t_{(k)}}^k(z)$ can be visualized as multiple,
{\it i.e.} composed, $\L$-integral $m$-curves with source the point
$z$. Also, the point $z'=L_{t_k}^k\circ \cdots \circ L_{t_1}^1(z)$ is
the endpoint of a piecewise smooth $m$-curve with origin $z$: follow
$L^1$ during $m$-time $t_1$, $\ldots$, follow $L^k$ during $m$-time
$t_k$. Now, following a well-established terminology, we shall say that a
manifold $\Lambda\subset \Delta^n$ is called {\it $\SS$-integral} if
$T_z\Lambda \supset \text{\rm Span}_{\C}
\SS(z)$ for all $z\in \Lambda$. Then for
each $L\in \SS$, $L|_{\Lambda}$ is tangent to $\Lambda$. Thus, it is
clear that any integral $m$-curve of an element $L\in \SS$ with origin
a point $z$ of an $\SS$-integral manifold $\Lambda$ stays in
$\Lambda$. Now, we introduce the following definitions.  Let $z\in
(\frac{r}{2}\Delta)^n$. The {\it $\SS$-orbit of $z$} in $\Delta^n$,
$\Cal O_{\SS}(\Delta^n,z)$ is the set of all points
$\L_{t_{(k)}}^k(z)\in (r\Delta)^n$ for all $t_{(k)}\in
(\delta\Delta^m)^k$, $k\leq 3n$.  Finally, As we will be interested
only in the $\SS$-orbit of $0$, we can localize at $z=0$. We shall say
that the open set $(r\Delta)^n$ is {\it $\SS$-minimal at $0$}
if $\Cal O_{\SS}(\Delta^n,0)$ contains a polydisc $(\varepsilon
\Delta)^n$, $\varepsilon>0$.

\subhead 9.3.~The orbit Theorem \endsubhead
Now, we propose a self-contained proof, inspired by 
[Sus] and which does not use
Lie brackets, of the following special case of Nagano's theorem [N].

\proclaim{Theorem 9.4} \text{\rm ([N], [Sus])}
There exists $\varepsilon >0$ such that the $\SS$-orbit of $0$
consists of a closed complex $\SS$-integral manifold-piece $\Cal O_0$
through $0$ in the polydisc $(\varepsilon
\Delta)^n$.
\endproclaim

\noindent
It is easy to check that any $\SS$-integral manifold $\Lambda$ passing
through the origin must contain the manifold-piece $\Cal O_0$, so 
Theorem~9.4 explains what is the smallest such $\SS$-integral
manifold-piece. Also,
as a coorollary, we see that the open set $(r\Delta)^n$ is
$\SS$-minimal at $0$ if and only if the dimension of the $\SS$-orbit is
maximal, {\it i.e.} $\text{\rm dim}_{\C} \Cal O_0 =n$. Let us
introduce the following notation. By $\Cal V_X(p)$, we shall denote a
small open polydisc neighborhood of the point $p$ in the complex
manifold $X$.

\demo{Proof of Theorem~9.4}
Notice that if $a=1$, $\Cal O_0$ is just the $m$-curve through $0$ of
the single element of $\SS$, so we assume $a\geq 2$ in the sequel.
The following definitions will generalize the notion of Segre
multitype given in \S7 above. At first, we need some preliminary.  If
$ \L^k\in \SS^k$ ($k \leq 3n$), we shall denote by $\Gamma_{\L^k}$ the
holomorphic map $(\delta \Delta^m)^k\ni t_{(k)} \mapsto
\L_{t_{(k)}}^k(0)$. By the definition of $\SS$-orbits, it is clear that
$\L_{t_{(k)}}^k(0)\in \Cal O_{\SS}(\Delta^n,0)$. Also, let us recall
that given $f:X\to Y$ a holomorphic map of {\it complex connected
manifolds}, there exists a proper complex subvariety $Z\subset X$ with
$\text{\rm dim}_{\C} Z < \text{\rm dim}_{\C} X$ such that $\text{\rm
rk}_{\C,p}(f)=\max_{q\in X} \text{\rm rk}_{\C,q}(f)$ for all $p\in
X\backslash Z$. This integer is called the {\it generic rank of $f$}
and we shall denote it by $\text{\rm gen-rk}_{\C}(f)$. Of course, if
$U$ is an arbitrary open subset of $X$ then we have $\text{\rm
gen-rk}_{\C}(f|_U)= \text{\rm gen-rk}_{\C}(f)$, thanks to the
principle of analytic continuation. Such properties of the generic
rank of holomorphic maps will be of crucial importance in the
construction of orbits and in the definition of the following
integers.  We construct indeed by induction a special sequence
$\L^{*k}:=(L^{*1},\ldots,L^{*k})$, $k\in \N_*$, as follows.  First,
let us define $\L^{*a}=(L^{*1},\ldots,L^{*a})=(L_1,\ldots,L_a)$, {\it
i.e.}  $L^{*\alpha}= L_{\alpha}$ for $1\leq \alpha\leq a$.  By our
assumption on $\SS$ (shrinking $\delta$ if necessary), we have
$\text{\rm rk}_{\C,0}(\Gamma_{\L^{*a}})= \text{\rm
gen-rk}_{\C}(\Gamma_{\L^{*a}})=
\text{\rm rk}_{\C,t_{(a)}}(\Gamma_{\L^{*a}})=am$, for all $t_{(a)}\in
(\delta\Delta^m)^a$. Let $\alpha\in [\![1,a]\!]$. Given $\L^{*a}$ so
defined, we define $\L^{*a}L_{\alpha}:=(\L^{*a},L_{\alpha})$, an
$(a+1)$-tuple of
elements of $\SS$, and we denote by
$\Gamma_{\L^{*a}L_{\alpha}}$ the map $(\delta\Delta^m)^a\times
(\delta\Delta^m)\ni (t_{(a)},t_{a+1})\mapsto L_{\alpha t_{a+1}}
\circ \L_{t_{(a)}}^{*a}(0)$, which is consistent with our
previous notations. Because
$\Gamma_{\L^{*a}L_{\alpha}}((t_{(a)},0)\equiv
\Gamma_{\L^{*a}}(t_{(a)})$, it is clear that
there exist well-defined
integers $0\leq e_1(\alpha)\leq e_1\leq n-am$ satisfying
$$
\text{\rm gen-rk}_{\C}
(\Gamma_{\L^{*a}L_{\alpha}}):= e_1(\alpha)+am, \ \ \ \ \
\text{\rm and we set} \ \ \ \ \
e_1:=\sup_{1\leq \alpha\leq a} e_1(\alpha).
\tag 9.5
$$ 
If $e_1=0$, our construction stops. If $e_1>0$, we choose an $\alpha$
with $e_1(\alpha)=e_1$ and we define the $(a+1)$-th vector field to be
$L^{*a+1}:=L_{\alpha}$ with this $\alpha$.  So we have completed the
choice of $\L^{*a+1}:=(\L^{*a},L^{*a+1})$. Inductively now, we assume
that $\L^{*k}$ is already defined, with corresponding integers
$e_1\geq 1,\ldots,e_k\geq 1$ satisfying $\text{\rm gen-rk}_{\C}
(\Gamma_{\L^{*k}})=am+e_1+\cdots+ e_k$. We shall write for short
$e_{\{k\}}:=e_1+\cdots+e_k$. Let $\alpha\in [\![1,a]\!]$. To construct
the $(k+1)$-th $m$-vector field $L^{*k+1}$, we define
$\L^{*k}L_{\alpha}:=(\L^{*k},L_{\alpha})$ and
$\Gamma_{\L^{*k}L_{\alpha}}: (\delta\Delta^m)^k\times (\delta\Delta^m)
\ni (t_{(k)},t_{k+1})\mapsto L_{\alpha t_{k+1}}\circ
\L_{t_{(k)}}^{*k}(0)$. Again, because
$\Gamma_{\L^{*k}L_{\alpha}}(t_{(k)},0)\equiv
\Gamma_{\L^{*k}}(t_{(k)})$, it is clear
that there are well-defined integers $0\leq e_{k+1}(\alpha)\leq
e_{k+1}\leq n-am-e_{\{k\}}\leq n-am-k$ satisfying
$$
\text{\rm gen-rk}_{\C}
(\Gamma_{\L^{*k}L_{\alpha}}):= e_{k+1}(\alpha)+e_{\{k\}}+am, \ \ \ \ \
\text{\rm and we set} \ \ \ \ \
e_{k+1}:=\sup_{1\leq \alpha\leq a} e_{k+1}(\alpha).
\tag 9.6
$$ 
If $e_{k+1}=0$, our construction stops. If $e_{k+1}>0$, we choose
$\alpha$ with $e_{k+1}(\alpha)=e_{k+1}$ and we define the $(k+1)$-th
vector field to be $L^{*k+1}=L_{\alpha}$, with this $\alpha$. We
define $e_{\{k+1\}}= e_{\{k\}}+e_{k+1}$. So we have completed the
choice of $\L^{*k+1}:=(\L^{*k},L^{*k+1})$, hence also 
the choice of all the $\L^{*k+1}$. We can now define $\kappa_0$ to be
the smallest integer such that
$e_{\kappa_0+1}=0$ and we set $\mu_0:=a+\kappa_0$.
After $\Gamma_{\L^{*\mu_0}}$, {\it i.e.} for $k>\mu_0$, the 
generic ranks of the maps $\Gamma_{\L^k}$ cease to increase.
Of course, we have the inequalities
$\kappa_0\leq n-am$ and $\mu_0\leq n$ (since $m\geq 1$).
Paralleling the definitions in \S7, we shall call 
the integer $\mu_0$ is a {\it minimality type} of $\SS$ at $0$.
$\mu_0$ depends on the subsequent choices of the $\alpha$'s.
Although the integers $\mu_0$ and $e_k$ may very well depend on the
subsequent choices of the $\alpha$'s which maximizes the subsequent
generic ranks of the maps $\Gamma_{\L^{*k}L_\alpha}$, it will follow
in the end of the proof of Theorem~9.6 that the integer $e=\text{\rm
dim}_{\R} \Cal O_{\SS}(\Delta^n,0)=am+e_{\{\kappa_0\}}$ will not 
depend on the choices of the $\alpha$'s.

To summarize what we have done so far, we have constructed a $\mu_0$-tuple
of $m$-vector fields
$\L^{*\mu_0}=(L_1,\ldots,L_a,L^{*a+1},\ldots,L^{*\mu_0})=
(L^{*1},\ldots,L^{*\mu_0})$ which satisfies
$$
m<2m<\cdots<a(m-1)<am<am+e_1<\cdots<am+e_{\{\kappa_0\}}.
\tag 9.7
$$ 
We define $e:=am+e_{\{\kappa_0\}}$ and we call the $\mu_0$-tuple
$(\overbrace{m,m,\ldots,m}^{a\ \text{\rm times}},e_1,\ldots,
e_{\kappa_0})\in \N^{\mu_0}$ will be called a {\it minimality
multitype} of $\SS$ at $0$, where $\mu_0=a+\kappa_0$.  We remark that
if $a=2$ as in the particular application to CR geometry given in \S7,
if we denote the doubleton $\SS:=\{L, \underline{L}\}$, then
$\L^k_{t_{(k)}}(0)$ is written $\cdots\underline{L}_{t_2}\circ
L_{t_1}(0)$ or $\cdots L_{t_2}
\circ \underline{L}_{t_1}(0)$ and there is no
other possibility. Therefore, if $a=2$, there are exactly at most
minimality multitypes.  Furthermore, if the pair of vector fields
$\{L, \underline{L}\}$ satisfies a particular symmetry condition
$\sigma_*(L)=\underline{L}$, where $\sigma$ is a biholomorphism, as
the pair $\{\Cal L,\underline{\Cal L}\}$ in \S5-7 does, then the two
minimality multitypes must infact coincide. We believe that this explains
a particular feature of the geometry of complexified Segre 
varieties. Now, we can state a slightly more precise version of 
Theorem~9.4 similar to that of Theorem~7.5 as follows.

\proclaim{Theorem~9.8}
Let $\SS=\{L_{\alpha}\}_{1\leq \alpha \leq a}$ be as above be a system
of $m$-vector fields over $\C^n$, with $n=am+d$, let $\mu_0$ be a
minimality type of $\SS$ at 0, with $\mu_0\leq a+d$ and let
$(m,\ldots,m,e_1,\ldots,e_{\kappa_0})$ be the associated multitype,
where $\mu_0=a+\kappa_0$. Then there exists a $\mu_0$-tuple of
$m$-vector fields $\L^{*\mu_0}= (L^1,\ldots,L^{*\mu_0})\in
\SS^{\mu_0}$, there exists an element $t_{(\mu_0)}^*\in
(\delta\Delta^m)^{\mu_0}$ arbitrarily close to the origin of the form
$(t_1^*,\ldots,t_{\mu_0-1}^*,0)$, and there exists a neighborhood
$\Cal W^*$ of $t_{(\mu_0)}^*$ in $(\delta \Delta^m)^{\mu_0}$ such that
after putting
$\underline{\L}^{*\mu_0-1}:=(L^{*\mu_0-1},\ldots,L^{*1})$ and
$\underline{t}_{(\mu_0-1)}^*:=(-t_{\mu_0-1}^*,\ldots,-t_1^*)$, then we
have\text{\rm \,:}
\roster
\item"{\bf 1)}"
The map $\Gamma_{\L^{*\mu_0}}$ is of rank $am+e_1+\ldots+e_{\kappa_0}$
at $t_{(\mu_p)}^*$.
\item"{\bf 2)}" 
$\underline{\L}^{*\mu_0-1}_{\underline{t}_{(\mu_0-1)}^*}
\circ \L_{t_{\mu_0}^*}^{*\mu_0}(0)=0$.
\item"{\bf 3)}"
The map $\underline{\L}_{\underline{t}_{(\mu_0-1)}^*}^{*
\mu_0-1}\circ \Gamma_{\L^{*\mu_0}}\: \Cal W^* \times
\underline{t}_{(\mu_0-1)}^*\to \Delta^n$ is of constant rank 
equal to $am+e_1+\cdots+e_{\kappa_0}$.
\item"{\bf 4)}"
Its image $\underline{\L}_{\underline{t}_{(\mu_0-1)}^*}^{*
\mu_0-1}\circ \Gamma_{\L^{*\mu_0}} (\Cal W^*)
= \Cal O_0$ is a manifold-piece through $0$ of the $\SS$ orbit $\Cal
O_\SS(\Delta^n,0)$. This means that every element of $\SS$ is tangent
to $\Cal O_0$ and that every $\SS$-integral manifold-piece $\Lambda_0$ 
through $0$ must \text{\rm contain} $\Cal O_\SS(\Delta^n,0)$.
\item"{\bf 5)}"
$am+e_1+\cdots+e_{\kappa_0}=\dim_\C \, \Cal O_0$.
\endroster
\endproclaim

\noindent
As the proof of Theorem~9.8 now goes essentially the same way as 
the proof of Theorem~7.5, we shall omit to repeat the detailed
arguments here. The proof of Theorem~9.6 is 
complete.
\qed
\enddemo 

Now, a final remark. As the above considerations are of course
valuable in the $\R$-analytic category, we can deduce from Theorem~9.6
the theorem of Nagano about existence of CR orbits of a CR-generic
$\Cal C^\omega$ manifold $M$ as in \S5-7, which we recovered in
Proposition~6.5 after passing to the extrinsic complexification.

\head \S10. Segre geometry of formal CR manifolds \endhead

Some (straightforward) adaptations are needed to develope the theory
of Segre chains in the {\it real and complex algebraic category}, in
particular to establish the algebraicity of CR-orbits, {\it etc.} The
reader may consult {\it e.g.} [M2] for a complete account about such
modifications. Also, it is easy to observe that most of the theory of
Segre chains extends to the category of formal CR manifolds. Such
objects are given by the ring $\C\dl t,\bar t\dr/((\rho_j(t,\bar
t))_{1\leq j\leq d})$, where the $d$ formal power series 
$\rho_j(t,\bar t)\in
\C\dl t,\bar t\dr$ satisfy $\rho_j(0)=0$ and $\partial \rho_1
(0)\wedge \cdots \wedge \partial \rho_d(0)\neq 0$ and the theory
builds up from the $\rho_j$. In a recent preprint, entitled {\it
Dynamics of the Segre varieties of a real submanifold in complex
space} (arXiv.org/abs/math/0008112), 
Baouendi-Ebenfelt-Rothschild endeavour such a
study (but without the exposition of our geometric viewpoint). In fact,
almost each one of the geometric concepts we have encountered in our
route has an obvious {\it formal counterpart}\,: the {\it formal} extrinsic
complexification, the {\it formal} (conjugate) Segre varieties, the
antiholomorphic involution, the {\it formal} regular coordinates, the
{\it formal} vector fields as in~\thetag{5.10} and their {\it formal}
flow as in~\thetag{6.2} (one can even speak in a rigorous sense of
formal local foliations$\ldots$). Indeed, the definition in
eq.~\thetag{6.2} obviously has a formal meaning. The important fact is
that in the definition of Segre chains, we indeed define inductively
{\it formal} power series, as show
eqs.~\thetag{6.3},~\thetag{6.4}. The formal counterpart of the
notion of generic rank is simply the nonvanishing of a suitable minor
in the ring of formal power series. Thus, all our considerations
extend straightforwardly to the formal category, except perharps for
the proof of the orbit Theorem~7.5, because we use a {\it nonzero}
point $w_{(\mu_p)}^*$ at which the ``value'' of a formal power series
would be senseless. However, the reader may check that if well interpreted
in the formal category, the important map written in {\bf 5)} of
Theorem~7.5 is in fact a formal power series vanishing at $p^c=0$,
whence the proof we give can easily be generalized, as desired.

\remark{Important remark}
Does this theory have deep applications to the smooth
category\,?  We believe not. Indeed, the reader should be aware that
the formal orbits {\it at a point} of system of {\it smooth} vector
fields {\it does not} represent appropriately the {\it smooth} orbit
of that point, as shows the example in $\R^2$ with $L={\partial \over
\partial x}$ and $L'=e^{{-1\over x^2}} {\partial \over \partial
y}$. Here, the smooth orbit of $0$ is a neighborhood of $0$ in $\R^2$
whereas the formal orbit of the Taylor expansions of $L$ and $L'$ at
$0$ is the (formal) $x$-line at $0$. To pursue our informal review of
these questions, let us note that it is easy to show that in general,
the formal orbit is ``contained'' in the smooth orbit (exercise\,:
find the rigorous sense of this assertion and prove it). As
expectable, this shows that the formal theory is much closer to the
analytic one than to the smooth one. Further, it also explains why the
naive theory of formal CR orbits should have very poor applications to
solve the fine questions of regularity about {\it smooth} CR mappings.
\endremark

\head \S11. Application of the formalism to the regularity 
of CR mappings\endhead

We provide here a brief summary of what are the main properties to understand
the CR mapping regularity problems in terms of our formalism and especially 
in terms of {\it the flows of $\Cal L$ and of $\underline{\Cal L}$}.
The reader is referred to [BER2,3] and to [M3] for a complete
account of how the methods work and what are the main results.

\subhead 11.1.~Holomorphic and formal maps of analytic CR manifolds
\endsubhead
We consider a formal or holomorphic invertible mapping $h\: (M,p)\to
(M',p')$ between two real analytic CR manifolds. In coordinates $t\in
\C^n$ vanishing at $p$ and $t'\in \C^n$ vanishing at $p'$, this map
$h(t)=(h_1(t),\ldots,h_n(t))$ is an $n$-tuple of power series with
$h_j(t)\in \C\{t\}$ or $h_j(t)\in \C\dl t\dr$ vanishing at $0$,
$h_j(0)=0$, and its jacobian determinant $\hbox{det} \, ({\partial
h_j\over \partial t_k}_{1\leq j,k\leq n}(0))\neq 0$ is nonzero at $0$.
We denote by $z_j=\bar Q_j(w,\tau)$, $1\leq j\leq d$ and $z_j'=\bar
Q_j'(w',\tau')$ some real analytic equations of $\Cal M$ and of $\Cal
M'$.  According to the splitting of coordinates $(w',z')$, we split
the map $h$ in $h:=(g,f)$. Then the assumption that $h$ maps $\Cal M$
into $\Cal M'$ can be simply restated by saying that the $d$-vectorial
power series (formal or converging) $\bar f(\tau)-Q'(\bar g(\tau),
h(t))\in \C\{t,\tau\}^d$ (or $\in
\C\dl t,\tau\dr^d$) is identically zero in $\C\{w,\zeta,\xi\}^d$
(or in $\C\dl w,\zeta,\xi\dr^d$) after replacing $z$ by $\bar Q
(w,\tau)$ (or after replacing $\xi$ by $Q(\zeta, t)$ in
$\C\{w,z,\zeta\}^d$ or in $\C\dl w,z,\zeta\dr^d$). Of course, we
insist on the parallelism about the two (equivalent) possibilities of
formulating this restatement in coordinates. For short, we shall say
that $\bar f(\tau) \equiv\bar\Theta'(\bar g(\tau),h(t))$ on $\Cal M$.
We denote by $h^c:=(h,\bar h)$ the holomorphic (or formal) map $\Cal
M\to \Cal M'$. Finally, we introduce the $m$-vector fields $\Cal L$
and $\underline{\Cal L}$ of \S5 and we abbreviate their multiple
concatenated flow maps (presented in \S2.2 and in \S6) by $\Gamma_k$
and by $\underline{\Gamma}_k$.

\subhead 11.2.~Interest of flows of CR vector fields \endsubhead
In terms of the flows of the CR vector fields, the recipe for
understanding the results of [BER2,3] and of [M3] is the following.
One of the main assumption in the results therein is that {\it $(M,p)$ is
minimal}. Without entering into all considerations, we shall explain what is
the central role played by the maps $\Gamma_k$ and we shall put in 
perspective the interest of considering {\it simultaneously} some
{\it reflection identities} together with some {\it conjugate reflection 
identities}, as  

\smallskip

\roster
\item"{\bf (I)}"
The first easy remark is the following.
As $\Gamma_k(w_{(k)})\in \Cal M$ for all $k\in \N$, we have
$h^c(\Gamma_k(w_{(k)}))\in \Cal M'$ for all $k\in \N$. At the formal
power series level, this property is expressed by the following
equation\,:
$$
\bar f(\Gamma_k(w_{(k)}))\equiv \Theta'(\bar g(\Gamma_k(w_{(k)})),
h(\Gamma_k(w_{(k)}))) \ \ \ \text{\rm in} \ \ \ \C\dl w_{(k)}\dr^d.
\tag 11.3
$$

\smallskip

\item"{\bf (II)}"
By Theorem~7.5, the minimality of $(M,p)$ is equivalent to the fact
that $\Gamma_{2\mu_p-1}$ is a submersion onto $(\Cal M,p^c)$. Then 
the relation \thetag{11.3} for $k\geq 2\mu_p-1$ becomes equivalent to 
the fundamental identity $\bar f(\tau) \equiv  Q'(\bar g(\tau), h(t))$ 
on $\Cal M$. Incidentally, this shows how minimality plays its first role.

\smallskip

\item"{\bf (III)}"
Now, we come to the most important step. For all $\beta\in \N^m$, we
consider the derivations $\underline{\Cal L}^\beta:=
\underline{\Cal L}^{1,\beta_1}\cdots \underline{\Cal L}^{m,\beta_m}$
and we apply them to the fundamental equation $\bar f(\tau) \equiv
\sum_{\gamma\in \N^m} \bar g(\tau)^\gamma \, Q_\gamma'(h(t))$. This
process is very classical. We get an infinite family of equations,
called {\it reflection identities}, of the form
$$
\underline{\Cal L}^\beta(\bar f)=\sum_{\gamma\in \N^m}
\underline{\Cal L}^\beta (\bar g(\tau)^\gamma) \  Q_\gamma'(h(t)), 
\ \ \ \ \ \forall \ \beta\in\N^m.
\tag 11.4
$$
However, we could also have applied these derivations to the {\it
conjugate} equation $f(t)\equiv \sumg g(t)^\gamma
\, \bar Q_\gamma'(h(\tau)))$ on $\Cal
M$ {\it and since there is no reason why not to do it, we do it.}

\smallskip

\item"{\bf (IV)}"
Consequently, we obtain two {\it a priori} different families
of equations\,:
$$
\left\{
\aligned
(*): \ \ f\equiv \bar Q'(g,\bar h), \ \ \ \ \ 0 \ \ \equiv & \sumg
g^\gamma \ \underline{\Cal L}^\beta (\bar{Q}_\gamma'(\bar h)), \
\ \ \ \ \forall \ \beta \in\N_*^m.\\ (\underline{*}): \ \
\bar{f}=Q'(\bar{g},h), \ \ \ \underline{\Cal L}^{\beta} \bar{f} \
\equiv & \sumg \underline{\Cal L}^\beta (\bar{g}^\gamma) \
Q_\gamma'(h), \ \ \ \ \ \forall \ \beta \in\N_*^m.
\endaligned
\right.
\tag 11.5
$$
These equations should be understood with $(t,\tau)\in \Cal M$.
Their respective bar-conjugates are the following\,:
$$
\left\{
\aligned
(\bar *): \ \ \bar f\equiv Q'(\bar g,h), \ \ \ \ \ 0 \ \ \equiv & \sumg
\bar g^\gamma \ {\Cal L}^\beta ({Q}_\gamma'(h)), \
\ \ \ \ \forall \ \beta \in\N_*^m.\\ (\overline{\underline{*}}): \ \
f=\bar Q'(g,\bar h), \ \ \ {\Cal L}^{\beta} f \
\equiv & \sumg {\Cal L}^\beta (g^\gamma) \
\bar Q_\gamma'(\bar h), \ \ \ \ \ \forall \ \beta \in\N_*^m.
\endaligned
\right.
\tag 11.6
$$
In summary, we obtain {\it four families of (conjugate) reflection
identities}.  To our knowledge, the system $(\underline{*})$ (or its
bar-conjugate $(\overline{\underline{*}})$) is nowhere considered in
the previous works on the subject, whereas $(*)$ and its bar-conjugate
$(\overline{*})$ are always used. The interest of the system
$(\underline{*})$ and the important information that it adds to 
study the regularity of the formal CR-reflection mapping is argued
in [M3], but we cannot enter into all the details here. We just want to 
say here that according to our heuristic principle of studying everything
{\it as pairs} in analytic CR geometry, we have discovered new reflection 
identities which are {\it not completely equivalent} to the known ones.

\smallskip

\item"{\bf (V)}"
Finally, we would like to make an last important observation which
shows how the flow maps $(w,q)\mapsto \Cal L_w(q)$ and the reflection
identities are linked together. Recall that if $q(x)\in \C\{ x\}^{2n}$
is a series vanishing at $0$ with $q(x)\in \Cal M$, for $x\in \C^\nu$,
then the derivative $[\partial_\zeta^\beta (\underline{\Cal
L}_\zeta(h^c(q(x)))]\v_{w=0}$ is equal to $[\underline{\Cal L}^\beta
h^c](q(x))$ (just by definition of flows). Then to obtain
eqs.~\thetag{11.4}, it is equivalent to differentiate the fundamental
equation $\bar f=Q'(\bar g, h)$ with respect to 
$\underline{\Cal L}^\beta$ as usual or to take the
differentiations
$$
\partial_\zeta^\beta\v_{\zeta=0} [\bar f(\underline{\Cal L}_\zeta(\tau))]=
\sumg \partial_\zeta^\beta\v_{\zeta=0} 
[\bar g(\underline{\Cal L}_\zeta(\tau))^\gamma \
Q_\gamma'(h(\underline{\Cal L}_\zeta(t)))], \ \ \ \ \forall \ \beta\in
\N^m.
\tag 11.7
$$
In particular, when differentiating eqs.~\thetag{11.3} with respect to
$w_k$ at $w_k=0$, we obtain either the system $(\underline{*})$ (if
$k$ is even) or the system $(\overline{\underline{*}})$ at the point
$(t,\tau)=\Gamma_{k-1}(w_{(k-1)})$ (if $k$ is odd). In conclusion to
this discussion and to this article, we would like to say that such an
identity~\thetag{11.7} is one of the key fact which explains the
mystery of the CR regularity properties in the minimal case, because
{\it this identity exhibits a natural relation between the reflection
identities and the flows generating the Segre chains, some two a priori
different objects which reveal therefore to be in fact intimately related
with each other}.
\endroster

\footnotetext[1]{pdf file~: 
cmi.univ-mrs.fr/$\sim$merker/index.html.}


\Refs\widestnumber\key{MM55}

\ref \key AG \by A. A. Agrachev and R. V. Gamkredlidze \paper
The exponential representation of flows and the chronological
calculus \jour Mat. Sbornik (N.S.) {\bf 107} (149) (1978), 467--532, 639.
English transl: Math USSR Sbornik \vol 35 \yr 1979 \pages 727--785\endref

\ref \key Ar \by V. I. Arnold \book Mathematical 
methods of Classical Mechanics,
\publ Graduate Texts in Mathematics,
\vol 60 \publaddr Springer, Heidelberg\yr 1989 \endref

\ref\key BER1 \by M.S. Baouendi, P. Ebenfelt and
L.P. Rothschild\paper Algebraicity of holomorphic mappings between
real algebraic sets in $\C^n$\jour Acta Mathematica \vol 177
\yr1996\pages 225--273
\endref

\ref\key BER2 \by M.S. Baouendi, P. Ebenfelt and
L.P. Rothschild \book Real submanifolds in complex space and their
mappings\publ Princeton Mathematical Series \vol 47\publaddr Princeton
University Press, Princeton, NJ\yr 1999
\endref

\ref\key BER3 \by M.S. Baouendi, P. Ebenfelt and
L.P. Rothschild \paper Convergence and finite determinacy of formal
CR mappings\jour Preprint 1999 
\endref

\ref\key Bell \by A. Bella\"{\i}che
\book SubRiemannian Geometry \publ Progress in Mathematices
\vol 144 \publaddr Brkh\"auser Verlag, Basel/Switzerland
\yr 1996 \pages 1--78
\endref 

\ref\key Belo \by V. K. Beloshapka \paper
On holomorphic transformation of a quadric \jour 
Mat. Sb. {\bf 182} (1991), no2, 203--219. English Transl. in 
Math. USSR Sbornik \vol 72 \yr 1992, no1, 189--205\endref

\ref\key BlGr \by T. Bloom and I. Graham\paper On
``type'' conditions for generic real submanifolds of $\C^n$\jour
Invent. Math \vol 40 \yr1977\pages 217--243
\endref

\ref\key Bo \by A. Boggess
\book CR manifolds and the
tangential Cauchy-Riemann complex \publ CRC Press
\publaddr Boca Raton\yr
1991\endref

\ref\key BuGo \manyby D. Burns and X. Gong \paper
Singular Levi-flat real analytic hypersurfaces\jour Amer. J. Math.
\vol 121 \yr 1999 \pages 23--53 \endref

\ref\key Chi \by E. Chirka \paper An introduction to the geometry of CR
manifolds\jour Russian Math. Surveys (no 1) \vol 46\yr
1991\pages 95--197\endref

\ref\key Cho \by W. L. Chow \paper \"Uber Systeme von linearen 
partiellen Differentialgleichungen erster Ordnung \jour
Math. Ann. \vol 117 \yr 1939 \pages 98--105\endref

\ref\key CMS \manyby B. Coupet, F. Meylan and A.  Sukhov \paper
Holomorphic maps of algebraic CR manifolds\jour International
Mathematics Research Notices\yr 1999\pages no 1, 1--29 \endref

\ref\key DP \manyby
Klas Diederich and Serguei Pinchuk \paper Proper holomorphic mappings
in dimension two etxend\jour Indiana University Mathematics Journal
\vol 44 \yr 1995\pages 1089--1126
\endref

\ref\key DW
\manyby K. Diederich and S.M. Webster 
\paper A reflection principle for degenerate hypersurfaces 
\jour Duke Math. J. \vol 47 \yr 1980 \pages 835--843\endref

\ref\key EIS \by V. V. Ezhov, A. V. Isaev and G. Schmalz \paper
Invariants of elliptic and hyperbolic CR structures of 
codimension $2$ \jour International Journal of Mathematics \vol
10 \yr 1999 \pages 1--52 \endref

\ref\key F1 \by M. Freeman \paper
The Levi form and local complex foliations \jour Proc. Amer. Math. Soc.
\vol 57 \yr 1976 \pages 369--370 \endref

\ref\key F2 \by M. Freeman \paper Real submanifolds with 
degenerate Levi forms \jour Proceedings of Symposia in Pure 
Mathematics XXX, part I, Several Complex Variables, Amer. 
Math. Soc., Providence, RI, 1977\endref

\ref\key Gro \by M. Gromov \paper 
Carnot-Carath\'eodory spaces seen from within  
\publ Progress in Mathematics
\vol 144 \publaddr Birkh\"auser Verlag, Basel/Switzerland
\yr 1996 \pages 79--323\endref

\ref\key Gru \by V. V. Gru$\check{\text{\rm s}}$in \paper
On a class of hypoelliptic operators \jour Mat. Sbornik (N.S.) {\bf 83}
(125) (1970), 456--473.
English transl: Math USSR Sbornik \vol 12 \yr 1970 \pages 458--476 \endref

\ref\key HaTr
\manyby N. Hanges and F. Treves \paper Propagation of holomorphic
extendability of CR functions \jour Math. Ann \vol 263 \yr 1983 \pages
157-177 \endref

\ref\key H\"o \by L. H\"ormander \paper Hypo-elliptic second order 
differential operators \jour Acta Math. \vol 119 \yr 1967 \pages 
147--171 \endref

\ref\key Hu \by X. Huang \paper Schwarz reflection principle in 
complex spaces of dimension two \jour Comm. Partial DIff. Equations
\vol 21 \yr 1996  \pages 1781--1828\endref

\ref\key J \by B. J\"oricke\paper Deformation of CR
manifolds, minimal points and CR-manifolds with the microlocal
analytic extension property\jour J. Geom. Anal \yr 1996\pages
555--611.
\endref

\ref\key K \by I. Kupka\paper G\'eom\'etrie sous-riemannienne \jour
S\'eminaire Bourbaki, Vol 1995/96. Ast\'erique $n^o$ 24 (1997),
Exp. $n^o$ 817, 5 \pages 352--380 \endref

\ref\key Le \by H. Lewy\paper On the boundary behaviour of 
holomorphic mappings\jour Rend. Acad. Naz. Lincei \vol 35
\yr1977\pages 1--8
\endref

\ref\key Lo \by A. Loboda \paper
Real analytic generating manifolds of codimension 2 in $\C^4$ and
their biholomorphic invariants
{\rm (translated from Russian)} \jour Math. USSR Izv. \vol 33 \yr 1989
\pages 295--315
\endref

\ref\key M1 \by J. Merker\paper Global minimality of
generic manifolds and holomorphic extendibility of CR functions\jour
Int. Math. Res. Not. \vol 8 \yr 1994\pages 329--342
\endref

\ref\key M2 \by J. Merker\paper On the partial algebraicity
of holomorphic mappings between real algebraic sets\jour
Preprint 1999 \toappear \endref

\ref\key M3 \by J. Merker\paper
On the convergence of formal CR mappings\footnotemark[1] 
\jour 3 manuscripts \yr 1999-2000 \endref

\ref\key Met \by G. M\'etivier \paper Fonctions spectrale et valeurs
propres d'une classe d'op\'erateurs non elliptiques \jour Comm. PArtial 
Diff. Equations \vol 1 \yr 1976 \pages 467--519 \endref

\ref\key Mit \by J. Mitchell \paper On Carnot-Carath\'eodory metrics \jour 
J. Diff. Geom. \vol 21 \yr 1985 \pages 35--45\endref

\ref\key MW \by J. Moser and S. Webster\paper Normal
forms for real surfaces in $\C^2$ near complex tangents and hyperbolic
surface transformation\jour Acta Mathematica \vol 150 \yr1983\pages
255--296 \endref

\ref\key N \by T. Nagano \paper Linear differential systems with 
singularities and applications to transitive Lie algebras
\jour J. Math. Soc. Japan \vol 18 \yr 1966 \pages 398--404 \endref

\ref\key P \by S. Pinchuk\paper On the analytic continuation
of holomorphic mappings, Mat. Sb. {\bf 98} (140) (1975),
416-435. English Transl. in Math. USSR Sb \vol 27 \yr1975\pages
375--392
\endref

\ref\key RS \by L.P. Rothschild and E. Stein \paper
Hypoelliptic differential operators
and nilpotent groups \jour Acta Mathematica \vol 137 \yr 1976 
\pages 247--320 \endref

\ref\key Se \by B. Segre\paper Intorno al problema di Poincar\'e
della rappresentazione pseudoconforme\jour Rendi. Acc. Lincei \vol 13
\yr 1931 \pages676--683
\endref

\ref\key SS \by R. Sharipov and A. Sukhov\paper
On CR mappings between algebraic Cauchy-Riemann manifolds and separate
algebraicity for holomorphic functions\jour Trans. Amer. Math. Soc
\vol 348 \yr 1996\pages 767--780
\endref

\ref\key Sp \by M. Spivak \book A comprehensive introduction to 
differential geometry \publ Vol 1, second edition, Publish or
Perish \publaddr Inc., Houston, Texas \yr 1979 \endref

\ref\key Sta
\by N. Stanton 
\paper Infinitesimal CR automorphisms of real hypersurfaces
\jour Amer. J. Math.\yr 1996\vol 118\pages 209--233\endref

\ref\key Ste \by E. Stein \paper Some problems in harmonic analysis
suggested by symmetric spaces and semi-simple groups, in \jour
Actes du Congr\`es international des math\'ematiciens, Nice, 1970 \vol
1 \pages 173--189\endref

\ref\key Suk2 \by A. Sukhov \paper Segre varieties, CR geometry
and Lie symmetries of second order PDE systems \jour Preprint
2000, arXiv.org/abs/math/0002197 \endref

\ref\key Suk2 \by A. Sukhov \paper Segre varieties and 
Lie symmetries \jour Math. Z. 
\toappear \endref

\ref\key SJ \by H. J. Sussmann and V. Jurdjevic \paper
Controllability of nonlinear systems \jour J. Differential Equations
\vol 12 \yr 1972 \pages 95--166\endref

\ref\key Sus \by H. J. Sussmann\paper Orbits of families of
vector fields and integrability of distributions\jour Trans.
Amer. Math. Soc \vol 180 \yr 1973\pages 171--188
\endref

\ref\key Trp \by J.-M. Tr\'epreau\paper Sur la
propagation des singularit\'es dans les vari\'et\'es CR\jour
Bull. Soc. Math. Fr \vol 118 \yr 1990 \pages 403--450
\endref
 
\ref\key Trv \by F. Treves \book Hypoanalytic
structures. Local theory \publ Princeton Mathematical Series \vol
40\publaddr Princeton University Press, Princeton, NJ \yr 1992
\endref
 
\ref\key Tu1 \by A. E. Tumanov\paper Extending CR
functions on a manifold of finite type over a wedge, Mat. Sbornik {\bf
136} (1988), 129--140. English transl. in \jour Math. of the Ussr
Sbornik \vol 64 \yr 1989 \pages 129--140
\endref

\ref\key Tu2 \by A. E. Tumanov\paper
Connections and propagation of Analyticity for CR Functions, Duke
Math. J \vol 73 \yr 1994\pages 1--24
\endref

\ref\key W \by S.M. Webster\paper On the mapping problem
for algebraic real hypersurfaces\jour Invent. Math \vol 43
\yr 1977\pages 53--68
\endref

\ref\key Z \by D. Zaitsev \paper Germs of local
automorphisms of real analytic CR structures and analytic dependence
on $k$-jets\jour Math. Res. Letters \vol 4 \yr1997\pages 823--842
\endref

\endRefs
\enddocument